\documentclass[10pt]{amsart}

\usepackage{amssymb, amsmath}
\usepackage{mathrsfs}
\usepackage{amscd}
\usepackage{verbatim}

\theoremstyle{plain}
\newtheorem{theorem}{Theorem}[section]
\theoremstyle{remark}
\newtheorem{remark}[theorem]{Remark}
\newtheorem{example}[theorem]{Example}
\theoremstyle{plain}
\newtheorem{corollary}[theorem]{Corollary}
\newtheorem{lemma}[theorem]{Lemma}
\newtheorem{proposition}[theorem]{Proposition}

\numberwithin{equation}{section}

\def\R{{\mathbb R}}

\newcommand{\E}{{\mathbb E}}
\renewcommand{\P}{{\mathbb P}}
\newcommand{\F}{{\mathcal F}}
\newcommand{\G}{{\mathcal G}}

\newcommand{\bbF}{{\mathbb F}}

\renewcommand{\a}{\alpha}
\renewcommand{\b}{\beta}
\newcommand{\g}{\gamma}
\renewcommand{\d}{\delta}
\newcommand{\e}{\varepsilon}
\renewcommand{\l}{\lambda}
\newcommand{\om}{\omega}
\renewcommand{\O}{\Omega}

\newcommand{\D}{{\mathcal D}}
\newcommand{\calL}{{\mathcal L}}
\newcommand{\n}{\Vert}
\newcommand{\one}{{{\bf 1}}}
\newcommand{\embed}{\hookrightarrow}
\newcommand{\s}{^*}
\newcommand{\lb}{\langle}
\newcommand{\rb}{\rangle}

\newcommand{\limn}{\lim_{n\to\infty}}

\newcommand{\sumn}{\sum_{n\ge 1}}

\newcommand{\Tend}{T_{0}}
\newcommand{\gL}{L^2_{\g}}
\newcommand{\Vinf}{V^{p}_{\alpha, \infty}([0,T]\times \O;E_\eta)}
\newcommand{\V}{\Vinf}
\newcommand{\VVinf}{V^{p}_{\alpha,\infty}([0,T_0]\times \O;E_\eta)}
\newcommand{\VV}{\VVinf}
\newcommand{\Vo}{V^{0}_{\a,p}([0,T]\times\O;E_{\eta})}
\newcommand{\VVo}{V^{0}_{\a,p}([0,T_0]\times\O;E_{\eta})}
\newcommand{\Vp}{{\tilde V}_{\a,p}^{p}([0,T]\times\O;E_{\eta})}
\newcommand{\VVp}{{\tilde V}_{\a,p}^{p}([0,T_0]\times\O;E_{\eta})}

\newcommand{\Vpwithouttilde}{{V}_{\a,p}^{p}([0,T]\times\O;E_{\eta})}
\newcommand{\Vadm}{V^{0,{\rm loc}}_{\alpha,p}([0,\varrho)\times \O;E_\eta)}

\newcommand{\Vadmone}{V^{0,{\rm loc}}_{\alpha,p}([0,\varrho_1)\times \O;E_\eta)}
\newcommand{\Vadmtwo}{V^{0,{\rm loc}}_{\alpha,p}([0,\varrho_2)\times \O;E_\eta)}

\begin{document}

\author{J.M.A.M. van Neerven}
\address{Delft Institute of Applied Mathematics\\
Delft University of Technology \\ P.O. Box 5031\\ 2600 GA Delft\\The
Netherlands} \email{J.M.A.M.vanNeerven@tudelft.nl}

\author{M.C. Veraar}\address{Delft Institute of Applied Mathematics\\
Delft University of Technology \\ P.O. Box 5031\\ 2600 GA Delft\\The
Netherlands\footnote{Current address: Mathematisches Institut I,
Technische Universit\"at Karlsruhe, D-76128 Karls\-ruhe, Germany}}
\email{mark@profsonline.nl}

\author{L. Weis}
\address{Mathematisches\, Institut\, I \\
\, Technische\, Universit\"at \, Karlsruhe \\
\, D-76128 \, Karls\-ruhe\\Germany} \email{Lutz.Weis@math.uni-karlsruhe.de}

\title[Stochastic evolution equations]{Stochastic evolution equations \\ in UMD Banach spaces}

\begin{abstract}
We discuss existence, uniqueness, and space-time H\"older regularity
for solutions of the parabolic stochastic evolution equation
\[\left\{\begin{aligned}
dU(t) & = (AU(t) + F(t,U(t)))\,dt + B(t,U(t))\,dW_H(t), \qquad t\in [0,\Tend],\\
 U(0) & = u_0,
\end{aligned}
\right.
\]
where $A$ generates an analytic $C_0$-semigroup on a UMD Banach space $E$
and $W_H$ is a cylindrical Brownian motion with values in a Hilbert space $H$.
We prove that if the
mappings $F:[0,T]\times E\to E$ and $B:[0,T]\times E\to \mathscr{L}(H,E)$
satisfy
suitable Lipschitz conditions and $u_0$ is $\F_0$-measurable and bounded,
then this problem has a unique mild
solution, which has trajectories in $C^\l([0,T];\D((-A)^\theta)$ provided
$\lambda\ge 0$ and $\theta\ge 0$ satisfy $\l+\theta<\frac12$.
Various extensions of this result are given and the results are applied to
parabolic stochastic partial differential equations.
\end{abstract}

\thanks{The first and second named authors are supported by
a `VIDI subsidie' (639.032.201) in the `Vernieuwingsimpuls'
programme of the Netherlands Organization for Scientific Research
(NWO). The second named author is also supported by the Humboldt
Foundation. The third named author is supported by a grant from the
Deutsche Forschungsgemeinschaft (We 2847/1-2).}


\date\today

\maketitle

\section{Introduction and statement of the results\label{intro}}

In this paper we prove existence, uniqueness, and space-time regularity
results for the abstract semilinear stochastic Cauchy problem
\begin{equation}\label{eq:SCP}\tag{SCP}
\left\{\begin{aligned}
dU(t) & = (AU(t) + F(t,U(t)))\,dt + B(t,U(t))\,dW_H(t), \qquad t\in [0,\Tend],\\
 U(0) & = u_0.
\end{aligned}
\right.
\end{equation}
Here $A$ is the generator of an analytic $C_0$-semigroup
$(S(t))_{t\ge 0}$ on a UMD Banach space $E$, $H$ is a separable Hilbert
space, and for suitable $\eta\ge 0$ the functions
$F:[0,T]\times \D((-A)^\eta)\to E$ and $B:[0,T]\times \D((-A)^\eta)\to \calL(H,E)$ enjoy
suitable Lipschitz continuity properties. The driving process $W_H$
is an $H$-cylindrical Brownian motion
adapted to a filtration $(\F_t)_{t\ge 0}$. In fact we shall allow
considerably less restrictive assumptions on $F$ and $B$; both
functions may be unbounded and may depend on the underlying probability
space.

A Hilbert space theory for stochastic evolution equations of the
above type has been developed since the 1980s by the schools of Da
Prato and Zabczyk \cite{DPZ92}. Much of this theory has been
extended to martingale type $2$-spaces \cite{Brz1,Brz2}; see also
the earlier work \cite{Nei}. This class of Banach spaces covers the
$L^p$-spaces in the range $2\leq p<\infty$, which is enough for many
practical applications to stochastic partial differential equations.
Let us also mention an alternative approach to the $L^p$-theory of
stochastic partial differential equations has been developed by
Krylov \cite{Kry}.

Extending earlier work of McConnell \cite{MCC}, the present authors
have developed a theory of stochastic integration in UMD spaces
\cite{NVWco,NVW} based on decoupling inequalities for UMD-valued
martingale difference sequences due to Garling \cite{Ga1, Ga2}. This
work is devoted to the application of this theory to stochastic
evolution equations in UMD spaces. In this introduction we will
sketch in an informal way the main ideas of our approach. For the
simplicity of presentation we shall consider the special case $H=\R$
and make the identifications $\calL(\R,E)=E$ and $W_\R=W$, where $W$
is a standard Brownian motion. For precise definitions and
statements of the results we refer to the main body of the paper.

A solution of equation \eqref{eq:SCP} is defined as
an $E$-valued adapted
process $U$ which satisfies the variation of constants formula
\[ U(t) = S(t)u_0 + \int_0^t S(t-s)F(s,U(s))\,ds +  \int_0^t
S(t-s)B(s,U(s))\,dW(s).\] The relation of this solution concept
with other type of solutions is considered in \cite{VThesis}. The
principal difficulty to be overcome for the construction of a
solution, is to find an appropriate space of processes which is
suitable for applying the Banach fixed point theorem. Any such space
$V$ should have the property that $U\in V$ implies that the
deterministic convolution
\[  t\mapsto  \int_0^t S(t-s)F(s,U(s))\,ds\]
and the stochastic convolution
\[t\mapsto  \int_0^t S(t-s)B(s,U(s))\,dW(s)\]
belong to $V$ again. To indicate why this such a space is difficult
to construct we recall a result from \cite{NVgammaLip} which states,
loosely speaking, that if $E$ is a Banach space which has the
property that $f(u)$ is stochastically integrable for every
$E$-valued stochastically integrable function $u$ and every
Lipschitz function $f:E\to E$, then $E$ is isomorphic to a Hilbert
space. Our way out of this apparent difficulty is by strengthening
the definition of Lipschitz continuity to $L_\g^2$-Lipschitz
continuity, which can be thought of as a Gaussian version of
Lipschitz continuity. From the point of view of stochastic PDEs,
this strengthening does not restrict the range of applications of
our abstract theory. Indeed, we shall prove that under standard
measurability and growth assumptions, Nemytskii operators are
$L_\g^2$-Lipschitz continuous in $L^p$. Furthermore, in type $2$
spaces the notion of $L_\g^2$-Lipschitz continuity coincides with
the usual notion of Lipschitz continuity.

Under the assumption that $F$ is Lipschitz continuous in the second
variable and $B$ is $L_\g^2$-Lipschitz continuous in the second
variable, uniformly with respect to bounded time intervals in their
first variables, the difficulty described above is essentially
reduced to finding a space of processes $V$ having the property that
 $\phi\in V$ implies that the pathwise deterministic convolutions
 \[  t\mapsto  \int_0^t S(t-s)\phi(s)\,ds\]
and the stochastic convolution integral
\begin{equation}\label{eq:stoch-conv1}
t\mapsto  \int_0^t S(t-s)\phi(s)\,dW(s)
\end{equation}
define processes which again belong to $V$.
The main tool for obtaining estimates for this stochastic integral is
$\g$-boundedness. This is the Gaussian version of the notion of
$R$-boundedness which in the past years has established itself as
a natural generalization to Banach spaces of the notion of
uniform boundedness in the Hilbert space context and which played an
essential role in much recent progress in the area of parabolic evolution equations.
The power of both notions derives from the fact that they connect probability in
Banach spaces with harmonic analysis.

From the point of view of stochastic integration, the importance of
$\g$-bounded families of operators is explained by the fact that they
act as pointwise multipliers in spaces of stochastically integrable processes.
This would still not be very useful if it were not the case that
one can associate $\g$-bounded families of operators with an analytic
$C_0$-semigroup $(S(t))_{t\ge 0}$ with generator $A$.
In fact, for all $\eta>0$ and $\e>0$, families such as
\[\big\{t^{\eta+\e} (-A)^{\eta}S(t): \ t\in (0,T_0)\big\}\]
are $\g$-bounded. Here, for simplicity, we are assuming that the fractional
powers of $A$ exist; in general one has to consider translates of $A$.
This suggests to rewrite the stochastic convolution \eqref{eq:stoch-conv1}
as
\begin{equation}\label{eq:stoch-conv2}
t\mapsto  \int_0^t
\big[(t-s)^{\eta+\e}(-A)^{\eta}S(t-s)\big](t-s)^{-\eta-\e}(-A)^{-\eta}\phi(s)\,dW(s).
\end{equation}
By $\g$-boundedness we can estimate the $L^p$-moments of this integral by the
$L^p$-moments of the simpler integral
\begin{equation}\label{eq:stoch-conv3}
t\mapsto  \int_0^t (t-s)^{-\eta-\e}(-A)^{-\eta}\phi(s)\,dW(s).
\end{equation}
Thus we are led to define $V_{\a,\infty}^{p}([0,T_0]\times \O;\D((-A)^{\eta}))$
as the space of all continuous adapted processes
$\phi:(0,T_0)\times\O\to \D((-A)^{\eta})$ for which the norm
\[\begin{aligned}
\ & \|\phi\|_{V_{\a,\infty}^{p}([0,T_0]\times\O;\D((-A)^{\eta}))}
\\ & := \Big(\E\|\phi\|_{C([0,T_0];\D((-A)^{\eta}))}^p\Big)^\frac1p +
\sup_{t\in [0,T_0]}\Big(\E \|
(t-\cdot)^{-\alpha}\phi(\cdot)\|_{\g(L^2(0,t),\D((-A)^{\eta}))}^p\Big)^{\frac{1}{p}}
\end{aligned}\]
is finite. Here, $\g(L^2(0,t),F)$ denotes the Banach space of
$\g$-radonifying operators from $L^2(0,t)$ into the Banach space
$F$; by the results of \cite{NW1}, a function $f:(0,t)\to F$ is
stochastically integrable on $(0,t)$ with respect to $W$ if and only
if it is the kernel of an integral operator belonging to $\g(L^2(0,t),F)$.

Now we are ready to formulate a special case of one of the main
results (see Theorems \ref{thm:mainexistenceL}, \ref{thm:Holdercont},
\ref{thm:Holdercontloc}).

\begin{theorem}\label{thm:main}
Let $E$ be a UMD space and let $\eta\ge 0$ and $p>2$ satisfy $\eta+\frac1p<\frac12$.
Assume that:
\begin{itemize}
\item[(i)]  $A$ generates an analytic $C_0$-semigroup on $E$;
\item[(ii)]  $F:[0,T_0]\times \D((-A)^\eta)\to E$ is Lipschitz continuous and of
linear growth in the
second variable, uniformly on $[0,T_0]$;
\item[(iii)]  $B:[0,T_0]\times \D((-A)^\eta)\to \calL(H,E)$ is
$L_\g^2$-Lipschitz continuous and of linear growth in the
second variable, uniformly on $[0,T_0]$;
\item[(iv)] $u_0\in L^p(\O,\F_0;\D((-A)^\eta))$.
\end{itemize} Then:
\begin{enumerate}
\item {\rm (Existence and uniqueness)}\
For all $\a> 0$ such that
$\eta+\frac1p<\a<\frac12$ the problem
\eqref{eq:SCP} admits a unique solution $U$ in
$V_{\a,\infty}^{p}([0,T_0]\times\O;\D((-A)^\eta))$.
\item {\rm (H\"older regularity)}\
For all $\l\ge 0$ and $\delta\ge \eta$ such that
$\l+\delta<\frac12$
the process $U-S(\cdot)u_0$ has a version with paths in $C^\lambda([0,T_0];\D((-A)^{\delta}))$.
\end{enumerate}
\end{theorem}

For martingale type $2$ spaces $E$, Theorem 1.1 was proved by
Brze\'zniak \cite{Brz2}; in this setting the $L^2_\gamma$-Lipschitz
assumption in (iii) reduces to a standard Lipschitz assumption. As
has already been pointed out, the class of martingale type $2$
spaces includes the spaces $L^p$ for $2\leq p<\infty$, whereas the
UMD spaces include $L^p$ for $1<p<\infty$. The UMD assumption in
Theorem \ref{thm:main} can actually be weakened so as to include
$L^1$-spaces as well; see Section \ref{sec:gen}. The assumptions on
$F$ and $B$ as well as the integrability assumption on $u_0$ can be
substantially weakened; we shall prove versions of Theorem
\ref{thm:main} assuming that $F$ and $B$ are merely locally
Lipschitz continuous and locally $L_\g^2$-Lipschitz continuous,
respectively, and $u_0$ is $\F_0$-measurable.

Let us now briefly discuss the organization of the paper.
Preliminary material on $\g$-radonifying operators, stochastic
integration in UMD spaces, and $\g$-boundedness of families of
operators, is collected in Section \ref{sec:prelimin}. In Sections
\ref{sec:convol} and \ref{sec:stoch_convol} we prove estimates for
deterministic and stochastic convolutions. After introducing the
notion of $L_\g^2$-Lipschitz continuity in Section
\ref{sec:g-Lipschitz} we take up the study of problem \eqref{eq:SCP}
in Section \ref{sec:main_results}, where we prove Theorem
\ref{thm:main}. The next two sections are concerned with refinements
of this theorem. In Section \ref{sec:meas_initial} we consider
arbitrary $\F_0$-measurable initial values, still assuming that the
functions $F$ and $B$ are globally Lipschitz continuous and
$L_\g^2$-Lipschitz continuous respectively. In Section
\ref{sec:loc_Lipschitz} we consider the locally Lipschitz case and
prove existence and uniqueness of solutions up to an explosion time.
In Section \ref{sec:gen} we discuss how the results of this paper can be extended
to a larger class of Banach spaces including the  UMD spaces as well as the spaces
$L^1$.

The final Section \ref{sec:applications}
is concerned with applications to stochastic partial differential equations.
On bounded smooth domains $S\subseteq \R^d$ we consider the parabolic problem
\begin{eqnarray*}
\nonumber \frac{\partial u}{\partial t}(t,s) &=& A(s,D) u(t,s) +
f(t,s,u(t,s))
\\ & & \qquad +\, g(t,s,u(t,s)) \, \frac{\partial w}{\partial t}(t,s), \quad s\in S, \ t\in
(0,T],
\\  B_j(s,D) u(t,s) &=& 0, \quad s\in \partial S, \ t\in (0,T],
\\  u(0,s) &=& u_0(s), \quad s\in S.
\end{eqnarray*}
Here $A$ is of the form
\[
A(s,D) = \sum_{|\alpha|\leq 2m} a_{\alpha}(s) D^{\alpha}
\]
with $D = -i(\partial_1, \ldots, \partial_d)$ and for $j=1, \ldots,
m$,
\[B_j(s,D) = \sum_{|\beta|\leq m_j} b_{j\beta}(s) D^{\beta}\]
where $1\leq m_j<2m$ is an integer. As a sample existence result, we
prove that if $f$ and $g$ satisfy standard measurability assumptions
and are locally Lipschitz and of linear growth in the third
variable, uniformly with respect to the first and second variables,
and if $u\in H^{2 m\eta,p}_{\{B_j\}}(S)$, then the above problem
admits a solution with paths in $C^{\lambda}([0,T];H^{2m
\delta,p}_{\{B_j\}}(S))$ for all $\d> \frac{d}{2mp}$ and $\lambda>0$
that satisfy $\delta+\lambda<\frac{1}{2}-\frac{d}{4m}$ and $2m
\delta - \frac1p \neq m_j$, for all $j=1, \ldots, m$. Uniqueness
results are obtained as well.

\medskip
All vector spaces in this paper are real.
Throughout the paper, $H$ and $E$ denote a separable Hilbert space and a Banach space,
respectively. We study the problem \eqref{eq:SCP} on a time interval
$[0,T_0]$ which is always considered to be fixed.
In many estimates below we are interested on bounds on
sub-intervals $[0,T]$ of
$[0,T_0]$ and it will be important to keep track of the dependence upon $T$ of
the constants appearing in these bounds. For this purpose we shall use the convention that
the letter $C$ is used
for generic constants which are independent of $T$
but which may depend on $T_0$ and all other relevant data in the estimates.
The numerical value of $C$ may vary from line to line.

We write $Q_1 \lesssim_A Q_2$ to express that there exists a
constant $c$, only depending on $A$, such that $Q_1\leq cQ_2.$ We
write $Q_1\eqsim_A Q_2$ to express that $Q_1 \lesssim_A Q_2$ and
$Q_2 \lesssim_A Q_1$.

\section{Preliminaries}\label{sec:prelimin}

The purpose of this section is to collect the basic stochastic tools
used in this paper. For proofs and further details we refer the
reader to our previous papers \cite{NVW, NW1}, where also references
to the literature can be found.

Throughout this paper, $(\O, \F, \P)$ always
denotes a complete probability space with a
filtration $(\F_t)_{t\geq 0}$. For a Banach space
$F$ and a finite measure space $(S,\Sigma, \mu)$, $L^0(S;F)$ denotes
the vector space of strongly measurable functions $\phi:S\to F$,
identifying functions which are equal almost everywhere.
Endowed with the topology induced by convergence
in measure, $L^0(S;F)$ is a complete metric space.

\medskip
\paragraph{\bf $\g$-Radonifying operators}
A linear operator $R:H\to E$ from a
separable Hilbert space $H$ into a Banach space
$E$ is called {\em $\g$-radonifying}
if for some (and then for every) orthonormal basis $(h_n)_{n\ge 1}$ of $H$ the Gaussian
sum
$ \sumn \g_n Rh_n$
converges in $L^2(\O;E)$. Here, and in the rest of the paper, $(\g_n)_{n\ge 1}$ is a
{\em Gaussian sequence}, i.e., a sequence of independent standard real-valued
Gaussian random variables. The space $\g(H,E)$ of all $\g$-radonifying
operators from $H$ to $E$ is a Banach space with respect to the norm
\[ \n R\n_{\g(H,E)} := \Big( \E \Big\n \sumn \g_n Rh_n
\Big\n^2\Big)^\frac12.\] This norm is independent of the orthonormal
basis $(h_n)_{n\ge 1}$. Moreover, $\g(H,E)$ is an operator ideal in
the sense that if $S_1:H'\to H$ and $ S_2:E\to E'$ are bounded
operators, then $R\in \g(H,E)$ implies $S_2RS_1\in \g(H',E')$ and
\begin{equation}\label{eq:ideal}
 \n S_2RS_1\n_{ \g(H',E')} \leq \n S_2\n \n R\n_{\g(H,E)} \n S_1\n.
\end{equation}
We will be mainly interested in the case where $H =
L^2(0,T;\mathcal{H})$, where $\mathcal{H}$ is another separable
Hilbert space.

The following lemma gives necessary and sufficient conditions for an
operator from $H$ to an $L^p$-space to be $\g$-radonifying.
It unifies various special cases in the literature, cf. \cite{BrzvN03, VTC}
and the references given therein.
In passing we note that by using the techniques of \cite{LiTz} the lemma can be generalized to arbitrary
Banach function spaces with finite cotype.

\begin{lemma}\label{thm:sq-fc-Lp}
Let $(S,\Sigma,\mu)$ be a $\sigma$-finite measure space and let
$1\leq p<\infty$. For an  operator $T\in\calL(H,L^p(S))$ the
following assertions are equivalent:

\begin{enumerate}
\item[\rm(1)]  $T\in \g(H,L^p(S))$;
\item[\rm(2)]  For some orthonormal basis $(h_n)_{n=1}^\infty$ of $H$ the function
$\big(\sumn |T h_n|^2\big)^{\frac12}$ belongs to $L^p(S)$;
\item[\rm(3)]  For all orthonormal bases $(h_n)_{n=1}^\infty$ of $H$ the function
$\big(\sum_{n=1}^\infty |T h_n|^2\big)^{\frac12}$ belongs to $L^p(S)$;
\item[\rm(4)]  There exists a function $g\in L^p(S)$ such that for all $h\in H$ we have
$|T h|\leq \|h\|_H\cdot g$ $\mu$-almost everywhere;
\item[\rm(5)]  There exists a function $k\in L^p(S;H)$ such that
$T h = [k(\cdot),h]_H$ $\mu$-almost everywhere.
\end{enumerate}
Moreover, in this situation we may take $k = \big(\sum_{n=1}^\infty |T
h_n|^2\big)^{\frac12}$ and have
\begin{equation}\label{eq:Rsquarekappa}
\|T\|_{\g(H,L^p(S))} \eqsim_p \Big\|\Big(\sum_{n=1}^\infty |T
h_n|^2\Big)^{\frac12}\Big\| \leq \|g\|_{L^p(S)}.
\end{equation}
\end{lemma}
\begin{proof}
By the
Kahane-Khintchine inequalities and Fubini's theorem we have,
for all $f_1,\dots,f_N\in L^p(S)$,
\[
\begin{aligned}
\ & \Big\|\Big(\sum_{n=1}^N |f_n|^2\Big)^{\frac12}\Big\|_{L^p(S)}
\\ & \qquad = \Big\|\Big(\E\Big|\sum_{n=1}^N \g_n
f_n\Big|^2\Big)^{\frac12}\Big\|_{L^p(S)}
\eqsim_p\Big\|\Big(\E\Big|\sum_{n=1}^N \g_n
f_n\Big|^p\Big)^{\frac1p}\Big\|_{L^p(S)}
\\ & \qquad = \Big(\E\Big\n\sum_{n=1}^N \g_n f_n\Big\|_{L^p(S)}^p\Big)^{\frac1p}
\eqsim_p
\Big(\E\Big\|\sum_{n=1}^N \g_n f_n\Big\|_{L^p(S)}^2\Big)^{\frac12}.
\end{aligned}\]
The equivalences (1)$\Leftrightarrow$(2)$\Leftrightarrow$(3) follow by
taking $f_n := Th_n$, $n=1,\dots,N$. This also gives the first part of
\eqref{eq:Rsquarekappa}.

(2)$\Rightarrow$(4): \ Let $g\in L^p(S)$ be defined as $g= \big(\sum_{n=1}^\infty |T
h_n|^2\big)^{\frac12}$. For $h = \sum_{n=1}^N a_n h_n$ we have,
for $\mu$-almost all $s\in S$,
\[|T h(s)| = \Big|\sum_{n=1}^N a_n T h_n(s) \Big|\leq \Big(\sum_{n=1}^N
|a_n|^2\Big)^{\frac12} \Big(\sum_{n=1}^N |T
h_n(s)|^2\Big)^{\frac12}\leq g(s)\|h\|_H.\] The case of a general
$h\in H$ follows by an approximation argument.

(4)$\Rightarrow$(5): \ Let $H_0$ be a countable dense set in $H$ which is
closed under taking $\mathbb{Q}$-linear combinations.
Let $N\in \Sigma$ be a $\mu$-null set such that for all $s\in \complement N$ and for
all $h\in H_0$, $|T h(s)|\le g(s)\|f\|_H$ and $h\mapsto T
h(s)$ is $\mathbb{Q}$-linear on $H_0$. By the Riesz representation theorem,
applied for each fixed
$s\in \complement N$, the mapping $h \to T h(s)$ has a unique extension
to an element $k(s)\in H$ with $T h(s) = [h, k(s)]_H$ for all $h\in
H_0$. By an approximation argument we obtain
that for all $h\in H$ we have $T h(s) = [h, k(s)]_H$ for $\mu$-almost all
$s\in S$. For all $s\in \complement N$,
\[\|k(s)\|_H = \sup_{\|h\|_H\le 1, h\in H_0} |[h,k(s)]| = \sup_{\|h\|_H\le 1, h\in H_0} |T h(s)| \le g(s).\]
Putting $k(s) = 0$ for $s\in N$, we obtain (5) and the last inequality
in \eqref{eq:Rsquarekappa}.

(5)$\Rightarrow$(3): \ Let $(h_n)_{n= 1}^\infty$ be an orthonormal basis for $H$.
Let $N\in \Sigma$ be a $\mu$-null set such that for all $s\in \complement N$ and all
$n\geq 1$ we have $T h_n(s) = [h_n,k(s)]$. Then for $s\in \complement N$,
\[\Big(\sum_{n=1}^\infty |T h_n(s)|^2\Big)^{\frac12} = \Big(\sum_{n=1}^\infty |[h_n, k(s)]|^2\Big)^{\frac12} = \|k(s)\|_H.\]
This gives (3) and the middle equality of \eqref{eq:Rsquarekappa}.
\end{proof}

Recall that for domains $S\subseteq \R^d$ and $\lambda>\frac{d}{2}$
one has $H^{\lambda, 2}(S) \hookrightarrow C_{\rm b}(\overline{S})$ (cf.
\cite[Theorem 4.6.1]{Tr}). Applying Lemma \ref{thm:sq-fc-Lp} with $g
\equiv C \cdot 1_S$ we obtain the following result.
\begin{corollary}\label{cor:gammaSobolevembed}
Assume $S\subseteq \R^d$ is a bounded domain. If
$\lambda>\frac{d}{2}$, then for all $p\in [1, \infty)$, the
embedding $I:H^{\lambda,2}(S)\to L^p(S)$ is $\g$-radonifying.
\end{corollary}

From the lemma we obtain an isomorphism of Banach spaces
\[L^p(S;H)\simeq \g(H, L^p(S)),\]
which is given by $ f\mapsto (h\mapsto [f(\cdot),h]_H)$.
The next result generalizes this observation:

\begin{lemma}[\cite{NVW}]\label{lem:gamma-Fubini}
Let $(S,\Sigma,\mu)$ be a  $\sigma$-finite measure space and let
$p\in[1,\infty)$ be fixed. Then $ f\mapsto (h\mapsto f(\cdot)h)$
defines an iso\-morphism of Banach spaces \[L^p(S;\g(H,E))\simeq
\g(H,L^p(S;E)).\]
\end{lemma}

\paragraph{\bf Stochastic integration}
In this section we recall some aspects of stochastic integration
in UMD Banach spaces. For proofs and more details we refer to our paper
\cite{NVW}, whose terminology we follow.

A Banach space $E$ is called a
{\em UMD} space if for some (equivalently, for
all) $p\in (1,\infty)$ there exists a constant $\beta_{p,E}\ge 1$ such that
for all $L^p$-integrable $E$-valued martingale difference sequences
$(d_j)_{j=1}^n$ and all $\{-1, 1\}$-valued sequence $(\e_j)_{j=1}^n$ we have
\begin{equation}\label{eq:UMD}
\Bigl(\E\Bigl\|\sum_{j=1}^n \e_j d_j \Bigr\|^p\Bigr)^\frac1p
\le \beta_{p,E} \,
\Bigl(\E\Bigl\|\sum_{j=1}^n d_j\Bigr\|^p\Bigr)^\frac1p.
\end{equation}
The class of UMD spaces was introduced in the 1970s by Maurey and Burkholder and has
been studied by many authors. For more information and references to the
literature we refer the reader to the review articles
\cite{BU3, RF}.
Examples of UMD spaces are all Hilbert spaces and the spaces $L^p(S)$
for $1<p<\infty$ and $\sigma$-finite measure spaces $(S,\Sigma,\mu)$.
If $E$ is a UMD space, then $L^p(S;E)$ is a UMD space for $1<p<\infty$.

Let $H$ be a separable Hilbert space.
An {\em $H$-cylindrical Brownian motion} is
family $W_H=(W_H(t))_{t\in [0,T]}$ of bounded linear
operators from $H$ to $L^2(\O)$ with the following two properties:
\begin{enumerate}
\item $W_H h = (W_H(t)h)_{t\in [0,T]}$ is real-valued Brownian motion for each $h\in H$,
\item $ \E (W_H(s)g \cdot W_H(t)h) = (s\wedge t)\,[g,h]_{H}$ for all $s,t\in [0,T], \ g,h\in H.$
\end{enumerate}
The stochastic integral of the indicator process
$1_{(a,b]\times A}\otimes (h\otimes x)$, where $0\le a<b<T$ and the
subset $A$ of $\O$ is $\F_a$-measurable, is defined as
\[\int_0^T 1_{(a,b]\times A}\otimes (h\otimes x)\,dW_H := 1_A(W_H(b)h - W_H(a)h)
x.\] By linearity, this definition extends to adapted step processes
$\Phi:(0,T)\times\O\to \calL(H,E)$ whose values are finite rank
operators.

In order to extend this definition to a more general class of
processes we introduce the following terminology. A process
$\Phi:(0,T)\times\O\to \calL(H,E)$ is called {\em $H$-strongly
measurable} if $\Phi h$ is strongly measurable for all $h\in H$.
Here, $(\Phi h)(t,\omega):= \Phi(t,\omega)h$. Such a process is
called {\em stochastically integrable} with respect to $W_H$ if it
is adapted and there exists a sequence of adapted step processes
$\Phi_n:(0,T)\times \O \to\calL(H,E)$ with values in the finite rank
operators from $H$ to $E$ and a pathwise continuous process
$\zeta:[0,T]\times\O\to E$, such that the following two conditions
are satisfied:
\begin{enumerate}
\item $\limn \Phi_n h = \Phi h$ in $L^0((0,T)\times\O;E)$ for all $h\in H$;
\item $\displaystyle \limn \int_0^\cdot \Phi_n\,dW_H = \zeta$ in
$L^0(\O;C([0,T];E))$.
\end{enumerate} In this situation,
$\zeta$ is determined uniquely as an element of $L^0(\O;C([0,T];E))$
and is called the {\em stochastic integral} of $\Phi$ with respect
to $W_H$, notation:
\[ \zeta = \int_0^\cdot \Phi\,dW_H.\]
The process $\zeta$ is a continuous local martingale starting at
zero. The following result from \cite{NVWco,NVW} states necessary and sufficient
conditions for stochastic integrability.

\begin{proposition}\label{prop:NVW}
Let $E$ be a UMD space. For an adapted $H$-strongly measurable
process $\Phi:(0,T)\times\O\to\calL(H,E)$ the
following assertions are equivalent: \begin{enumerate}
\item the process $\Phi$ is stochastically integrable with respect to $W_H$;
\item for all $x\s\in
E\s$ the process $\Phi\s x\s$ belongs to $L^0(\O;L^2(0,T;H))$, and there exists a pathwise continuous process $\zeta:[0,T]\times\O\to E$
such that for all $x\s\in E\s$ we have
\[\lb \zeta, x\s\rb =  \int_0^\cdot \Phi\s x\s\,dW_H \ \ \hbox{ in
$L^0(\O;C([0,T])$;}\]
\item\label{NVW3} for all $x\s\in
E\s$ the process $\Phi\s x\s$ belongs to $L^0(\O;L^2(0,T;H))$, and there exists an operator-valued random variable
$R:\O\to \g(L^2(0,T;H),E))$ such that for all $f\in L^2(0,T;H)$ and
$x\s\in E\s$ we have
\[ \lb Rf,x\s\rb = \int_0^T [f(t), \Phi\s(t)x\s]_H\,dt  \ \ \hbox{ in
$L^0(\O)$.}\] \end{enumerate} In this situation we have $\zeta =
\int_0^{\cdot} \Phi\,dW_H$ in $L^0(\O;C([0,T];E))$. Furthermore, for
all $p\in (1, \infty)$,
\[ \E \sup_{t\in [0,T]}\Big\n \int_0^t \Phi\,dW_H\Big\n^p \eqsim_{p,E}
\E\n R\n_{\g(L^2(0,T;H),E)}^p.\]
\end{proposition}

In the situation of (3) we shall say that $R$ is {\em represented} by $\Phi$.
Since $\Phi$ is uniquely determined almost everywhere on $(0,T)\times\O$ by
$R$ and vise versa (this readily follows from
\cite[Lemma 2.7 and Remark 2.8]{NVW}), in what follows we shall frequently
identify $R$ and $\Phi$.

The next lemma will be useful in Section \ref{sec:meas_initial}.

\begin{lemma}
\label{lem:localityloc} Let $\Phi:(0,T)\times\O\to\calL(H,E)$ be
stochastically integrable with respect to $W_H$. Suppose $A\in\F$ is
a measurable set such that for all $x\s\in E\s$ we have
\[\Phi\s(t,\omega)x\s = 0\ \ \hbox{for almost all $\,(t,\omega)\in (0,T)\times
A$}.\] Then almost surely in $A$, for all $t\in [0,T]$ we have
$ \int_0^t \Phi\,dW_H = 0.$
\end{lemma}
\begin{proof}
Let $x^*\in E^*$ be arbitrary. By strong measurability it suffices
to show that, almost surely in $A$, for all $t\in [0,T]$ we have
\[M_t:=\int_0^t \Phi^*x^*\,dW_H = 0.\]
For the quadratic variation of the continuous local martingale $M$
we have
\[[M]_T = \int_0^T \|\Phi^*(s)x^*\|^2 \, d s=0 \ \ \text{a.s. on $A$}.\]
Therefore, $M = 0$ a.s. on $A$. Indeed, let
\[\tau:= \inf\{t\in [0,T]: [M]_t
>0\},\]
where we take $\tau=T$ if the infimum is taken over the empty set.
Then $M^{\tau}$ is a continuous local martingale with
quadratic variation $[M^{\tau}] = [M]^{\tau}=0$. Hence $M^\tau=0$
a.s. This implies the result.
\end{proof}

\medskip
\paragraph{\bf $R$-Boundedness and $\g$-boundedness}
Let $E_1$ and $E_2$ be Banach spaces and let $(r_n)_{n\ge 1}$
be a {\em Rademacher sequence}, i.e., a sequence of independent random variables
satisfying $\P\{r_n = -1\} = \P\{r_n=1\}=\frac12$.
A family $\mathscr{T}$ of bounded linear
operators from $E_1$ to $E_2$ is called {\em $R$-bounded} if there exists
a constant $C\ge 0$ such that for all finite sequences $(x_n)_{n=1}^N$ in $E_1$
and $(T_n)_{n=1}^N$ in ${\mathscr {T}}$ we have
\[ \E \Big\n \sum_{n=1}^N r_n T_n x_n\Big\n^2 \le C^2\E \Big\n \sum_{n=1}^N r_n
x_n\Big\n^2.
\]
The least admissible constant $C$ is called the {\em $R$-bound} of
$\mathscr {T}$, notation $R(\mathscr{T})$. By the Kahane-Khintchine
inequalities the exponent $2$ may be replaced by any $p\in
[1,\infty)$. This only affects the value of the $R$-bound; we shall
use the notation $R_p(\mathscr{T})$ for the $R$-bound of
$\mathscr{T}$ relative to exponent $p$.

Upon replacing the Rademacher sequence by a Gaussian sequence we arrive at the
notion of a {\em $\g$-bounded} family of operators, whose
{\em $\g$-bound} will be denoted by $\g(\mathscr{T})$.
A standard randomization argument shows that every $R$-bounded family is
$\g$-bounded, and both notions are equivalent if the range space has finite
cotype (the definitions of type and cotype are recalled in the next section).

The notion of $R$-boundedness has played an important role in recent
progress in the regularity theory of parabolic evolution equations.
Detailed accounts of these developments are presented
in \cite{DHP,KuWe}, where more about the history of this concept and further
references to the literature can be found.

Here we shall need various examples of $R$-bounded families, which are stated in
the form of lemmas.

\begin{lemma}[\cite{We}]\label{lem:int-der} If $\Phi:(0,T)\to \calL(E_1,E_2)$ is differentiable
with integrable derivative, the family \[{\mathscr{T}_\Phi} =
\big\{\Phi(t): \ t\in (0,T)\big\}\] is $R$-bounded in
$\calL(E_1,E_2)$, with
\[R({\mathscr{T}_\Phi})\le \n \Phi(0+)\n + \int_0^T\n \Phi'(t)\n\,dt.\]
\end{lemma}

We continue with a lemma which connects the notions of
$R$-boundedness and $\g$-radonification. Let $H$ be a Hilbert space
and $E$ a Banach space. For each $h\in H$ we obtain a linear
operator $T_h: E\to \g(H,E)$ by putting \[T_hx:= h\otimes x, \qquad
x\in E.\]

\begin{lemma}[\cite{KaiWe}]\label{lem:KaiWe}
If $E$ has finite cotype, the family
\[\mathscr{T}=\big\{T_h:\ \n h\n_H\le 1\big\}\]
is $R$-bounded in $\calL(E,\g(H,E))$.
\end{lemma}

Following \cite{KWcalc}, a Banach space $E$ is said to have {\em property $(\Delta)$} if there exists
a constant $C_\Delta$ such that if $(r_n')_{n=1}^N$ and $(r_n'')_{n=1}^N$
are Rademacher sequences on probability spaces $(\O',\P')$ and $(\O'',\P'')$
respectively, and $(x_{mn})_{m,n=1}^N$ is a doubly indexed sequence of  elements of $E$, then
\[
\E'\E''\Big\n\sum_{n=1}^N \sum_{m=1}^n r_m' r_n'' x_{mn}\Big\n^2 \le C_\Delta^2
\E'\E''\Big\n\sum_{n=1}^N \sum_{m=1}^N r_m' r_n'' x_{mn}\Big\n^2 .
\]
Every UMD space has property $(\Delta)$ \cite{CPSW} and every Banach
space with property $(\Delta)$ has finite cotype. Furthermore the
spaces $L^1(S)$ with $(S,\Sigma,\mu)$ $\sigma$-finite have
property $(\Delta)$. The space of trace class operators does
not have property $(\Delta)$ (see \cite{KWcalc}).

The next lemma is a variation of Bourgain's
vector-valued Stein inequality for UMD spaces \cite{Bou, CPSW}
and  was kindly communicated to us
by Tuomas Hyt\"onen.

\begin{lemma} \label{lem:Stein}
Let $W_H$ be an $H$-cylindrical Brownian motion,  adapted to a
filtration $(\F_t)_{t\in [0,T]}$, on a probability space $(\O,P)$.
If $E$ is a Banach space enjoying property $(\Delta)$, then for all
$1\le p< \infty$ the family of conditional expectation operators
\[\mathscr{E}_p = \big\{\E(\cdot|{\F_t}): \ t\in [0,T]\big\}\] is
$R$-bounded, with $R$-bound $C_\Delta$, on the closed linear
subspace $G^p(\O;E)$ of $L^p(\O;E)$ spanned by all random variables of the
form $\int_0^T \Phi\,dW_H$ with $\Phi\in \g(L^2(0,T;H),E)$.
\end{lemma}
\begin{proof}
Let $1\le p<\infty$ be fixed and choose $\E_1,\dots,\E_N \in \mathscr{E}_p$,
say $\E_n = \E(\cdot|\F_{t_n})$ with $0\le t_n\le T$. By relabeling the indices we may assume that
$t_1\le \dots\le t_N$.
We must show that for all $F_1,\dots,F_N \in L^p(\O;E)$ of the form
$F_n = \int_0^T \Phi_n\,dW_H$ we have
\[ \E' \big\n\sum_{n=1}^N r_n' \E_n F_n \Big\n^2
\le C_\Delta^2  \E' \big\n\sum_{n=1}^N r_n'
F_n \Big\n^2.
\]
We write $\E_n = \sum_{j=1}^n D_j$, where $D_j := \E_j - \E_{j-1}$
with the convention that $\E_0 = 0$. The important point to observe
is that if $\Psi_j\in \g(L^2(0,T;H),E)$ and $G_j:= \int_0^T
\Psi_j\,dW_H$, the random variables $D_j G_j$ are symmetric and
independent. Hence, by a standard randomization argument,
\[\begin{aligned}
\ & \E' \big\n\sum_{n=1}^N r_n' \E_n F_n \Big\n_{G^p(\O;E)}^2
 =  \E' \big\n\sum_{n=1}^N\sum_{j=1}^n r_n' D_j F_n \Big\n_{G^p(\O;E)}^2
\\ & \quad =  \E' \big\n\sum_{j=1}^N D_j \sum_{n=j}^N r_n' F_n \Big\n_{G^p(\O;E)}^2
 =  \E'\E'' \big\n\sum_{j=1}^N r_j'' D_j \sum_{n=j}^N r_n'F_n
\Big\n_{G^p(\O;E)}^2
\\ & \quad \le C_\Delta^2  \E'\E'' \big\n\sum_{j=1}^N r_j'' D_j\! \sum_{n=1}^N r_n'  F_n \Big\n_{G^p(\O;E)}^2
 = C_\Delta^2  \E' \big\n\sum_{j=1}^N D_j\! \sum_{n=1}^N r_n'F_n
\Big\n_{G^p(\O;E)}^2
\\ & \quad = C_\Delta^2  \E' \big\n\E_N \sum_{n=1}^N r_n' F_n \Big\n_{G^p(\O;E)}^2
 \le C_\Delta^2  \E' \big\n \sum_{n=1}^N r_n' F_n \Big\n_{G^p(\O;E)}^2.
\end{aligned}
\]
\end{proof}

The next lemma, obtained in \cite{KaWe} for the case
$H=\R$, states that $\g$-bounded families act boundedly as pointwise
multipliers on spaces of $\g$-radonifying operators.
The proof of the general case is entirely similar.

\begin{lemma}\label{lem:KW}
Let $E_1, E_2$ be Banach spaces and let $H$ be a separable Hilbert
space. Let $T>0$. Let $M:(0,T)\to \calL(E_1,E_2)$ be function with
the following properties:
\begin{enumerate}
\item for all $x\in E_1$
the function $M(\cdot)x$ is strongly measurable in $E_2$;
\item the
range $\mathscr{M} = \{M(t): \ t\in (0,T)\}$ is $\gamma$-bounded in
$\calL(E_1,E_2)$.
\end{enumerate}
Then for all step functions $\Phi:(0,T)\to \calL(H,E_1)$ with values
in the finite rank operators from $H$ to $E_1$ we have
\begin{equation}\label{eq:KW}
\|M\Phi\|_{\g(L^2(0,T;H), E_2)}\le \gamma(\mathscr{M})\|\Phi\|_{\g(L^2(0,T;H),
E_1)}.
\end{equation}
Here, $(M\Phi)(t):= M(t)\Phi(t)$.
As a consequence, the mapping $\Phi\mapsto M\Phi$ has a unique extension to a bounded operator
from $\g(L^2(0,T;H), E_1)$ to $\g(L^2(0,T;H), E_2)$ of norm at most $\gamma(\mathscr{M})$.
\end{lemma}

In \cite{KaWe} it is shown that under slight regularity assumptions
on $M$, the $\g$-bounded\-ness is also a necessary condition.

\section{Deterministic convolutions}
\label{sec:convol}

After these preliminaries we take up our main line of study and begin
with some estimates for deterministic convolutions. The main tool will be
a multiplier lemma for vector-valued
Besov spaces, Lemma \ref{lem:Besov}, to which we turn first.

Let $E$ be a Banach space, let $I=(a,b]$ with $-\infty\le a<b\le
\infty$ be a (possibly unbounded) interval, and let $s\in (0,1)$ and
$1\le p, q\le \infty$ be fixed. Following \cite[Section 3.b]{Ko},
the Besov space $B_{p,q}^s(I;E)$ is defined as follows. For $h\in\R$
and a function $f:I\to E$, we define $T(h)f:I\to E$ as the translate
of $f$ by $h$, i.e.,
\[(T (h) f)(t) :=
  \begin{cases}
    f(t+h) & \text{if $t+h\in I$}, \\
    0 & \text{otherwise}.
  \end{cases}
\]
Put
\[ I[h] : =  \{t\in I:\ t+  h \in I\}\] and, for $f\in L^p(I;E)$
and $t>0$,
\[ \varrho_p(f,t) := \sup_{|h|\le t} \|T(h) f -
f\|_{L^p(I[h];E)}.
\]
Now define
\[B^s_{p, q}(I;E) := \{f\in L^p(I;E): \|f\|_{B^s_{p, q}(I;E)}<\infty\},\]
where
\begin{equation}\label{besov}
\|f\|_{B^s_{p, q}(I;E)} := \|f\|_{L^p(I;E)} + \Big(\int_0^1
\big(t^{-s}\varrho_p(f,t)\big)^q\, \frac{dt}{t}\Big)^\frac1q
\end{equation}
with the obvious modification for $q=\infty$. Endowed with the norm
$\|\cdot\|_{B^s_{p, q}(I;E)}$, $B^s_{p,q}(I;E)$ is a Banach space.

The following continuous inclusions hold
for all $s,s_1, s_2\in (0,1)$, $p, q, q_1,
q_2\in [1, \infty]$ with $q_1\le q_2$, $s_2 \le s_1$:
\[\ B_{p, q_1}^s(I; E) \embed B_{p, q_2}^s(I;E), \ B_{p, q}^{s_1}(I;
E) \embed B_{p, q}^{s_2}(I;E).\]
If $I$ is bounded, then also
\[B_{p_1, q}^s(I; E) \embed B_{p_2, q}^s(I;E)\]
for $1\le p_2\le p_1\le \infty$.

The next lemma will play an important role in setting up our basic
framework. We remind the reader of the convention, made at the end
of Section \ref{intro}, that constants appearing in estimates may
depend upon the number $T_0$ which is kept fixed throughout the
paper.

\begin{lemma}\label{lem:Besov}
Let $1\le q < p< \infty$, $s>0$ and $\a\ge 0$  satisfy
$s<\frac{1}{q} - \frac{1}{p}$ and $\a <\frac{1}{q} - \frac{1}{p}-s$, and let
$1\le r<\infty$.
For all $T\in [0,T_0]$ and
$\phi\in B^s_{p,r}(0,T;E)$ the function $t\mapsto t^{-\alpha}
\phi(t)\one_{(0,T)}(t)$ belongs to $B^{s}_{q,r}(0,T_0;E)$ and there
exists a constant $C\ge 0$, independent of $T\in [0,T_0]$, such that
\[\|t\mapsto t^{-\alpha} \phi(t)\one_{(0,T)}(t)\|_{B^{s}_{q,r}(0,T_0;E)} \le
C T^{\frac1q-\frac1p-s-\alpha} \|\phi\|_{B^s_{p,r}(0,T;E)}.
\]
\end{lemma}
\begin{proof}
We prove the lemma under the additional assumption that $\a>0$; the
proof simplifies for case $\a=0$. We shall actually prove the
following stronger result
\[\|t\mapsto t^{-\alpha} \phi(t)\one_{(0,T)}(t)\|_{B^{s}_{q,r}(\R;E)} \le
C T^{\frac1q-\frac1p-s-\alpha} \|\phi\|_{B^s_{p,r}(0,T;E)}
\]
with a constant $C$ independent of $T\in [0,T_0]$.

Fix $u\in [0,T]$ and $|h|\le u$. First assume that $h\ge 0$. Then
$I[h]=[0,T-h]$ and, by H\"older's inequality,
\[\begin{aligned}
\Big(\int_{\R} &\Big\|\frac{\phi(t+h)\one_{(0,T)}(t+h)-\phi(t)
\one_{(0,T)}(t)}{(t+h)^{\alpha}}\Big\|^q \,dt\Big)^{\frac{1}{q}} \\
& \leq \Big(\int_{-h}^0
\Big\|\frac{\phi(t+h)}{(t+h)^{\alpha}}\Big\|^q
 \,dt\Big)^{\frac{1}{q}}
+ \Big(\int_{0}^{T-h}
\Big\|\frac{\phi(t+h)-\phi(t)}{(t+h)^{\alpha}}\Big\|^q
\,dt\Big)^{\frac{1}{q}} \\ & \qquad + \Big(\int_{T-h}^T
\Big\|\frac{\phi(t)}{(t+h)^{\alpha}}\Big\|^q \,dt\Big)^{\frac{1}{q}}
\\ & \leq C u^{\frac1q-\frac1p-\alpha} \|\phi\|_{L^p(0,T;E)}  + C T^{\frac1q-\frac1p-\alpha}  \Big(\int_{I[h]} \|\phi(t+h)-\phi(t)\|^p \,
dt\Big)^{\frac1p}.
\end{aligned}\]
Again by H\"older's inequality,
\[\begin{aligned}
\Big(\int_{\R} \Big\|\frac{\phi(t) \one_{(0,T)}(t)}{(t+h)^{\alpha}}
-&\frac{\phi(t)\one_{(0,T)}(t)}{t^{\alpha}}\Big\|^q
\,dt\Big)^{\frac{1}{q}}
\\ & \leq \Big(\int_{0}^{T} \big|(t+h)^{-\alpha} -
t^{-\alpha}\big|^{\frac{pq}{p-q}} \,dt\Big)^{\frac{p-q}{pq}}
\|\phi\|_{L^p(0,T;E)}
\end{aligned}\]
with
\[
 \int_{0}^{T} \big|(t+h)^{-\alpha} -
t^{-\alpha}\big|^{\frac{pq}{p-q}} \,dt
\le \int_{0}^{T} t^{-\frac{\alpha
pq}{p-q}} - (t+h)^{-\frac{\alpha pq}{p-q}}  \,dt \le C
h^{1-\frac{\alpha pq}{p-q}}\le C
u^{1-\frac{\alpha pq}{p-q}}.
\]
Combining these estimates with the triangle inequality we obtain
\[
\begin{aligned} \ & \Big(\int_{\R}
\Big\|\frac{\phi(t+h)\one_{(0,T)}(t+h)}{(t+h)^{\alpha}}
-\frac{\phi(t) \one_{(0,T)}(t) }{t^{\alpha}}\Big\|^q
\,dt\Big)^{\frac{1}{q}}
\\ & \qquad \qquad\le C
u^{\frac1q-\frac1p-\alpha } \|\phi\|_{L^p(0,T;E)}
+   CT^{\frac1q-\frac1p-\alpha}
\Big(\int_{I[h]} \|\phi(t+h)-\phi(t)\|^p \,
dt\Big)^{\frac1p}.
\end{aligned}
\]
A similar estimate holds for $h\le 0$.

Next we split $[0,1]=[0,T\wedge 1]\cup [T\wedge 1,1]$ and estimate
the integral in \eqref{besov}. For the first we have
\[\begin{aligned}
\Big(&\int_0^{T\wedge 1} u^{-sr} \sup_{|h|\le u}\Big\|t\mapsto
\frac{\phi(t+h)\one_{(0,T)}(t+h)}{(t+h)^{\alpha}} -\frac{\phi(t)
\one_{(0,T)}(t)}{t^{\alpha}}\Big\|_{L^q(\R;E)}^r   \,
\frac{du}{u}\Big)^{\frac{1}{r}}
\\ & \le C
\Big(\int_0^{T\wedge 1}  u^{-sr} \Big[ T^{\frac1q-\frac1p-\alpha}
\sup_{|h|\le u} \|\phi(\cdot+h)-\phi(\cdot)\|_{L^p(I[h];E)}
\\ & \hskip7.6cm  + u^{\frac{p-q}{pq} -\alpha} \|\phi\|_{L^p(0,T;E)} \Big]^{r}\, \frac{du}{u}\Big)^{\frac{1}{r}}
\\ & \stackrel{\rm(i)}{\le} C
T^{\frac1q-\frac1p-\alpha} \Big(\int_0^1 u^{-sr} \Big[ \sup_{|h|\le u}
\|\phi(\cdot+h)-\phi(\cdot)\|_{L^p(I[h];E)}\Big]^{r}\, \frac{du}{u}\Big)^{\frac{1}{r}}
\\ & \hskip5cm + C
\Big(\int_0^{T} u^{-sr} u^{\frac{(p-q)r}{pq} -\alpha r} \,
\frac{du}{u}\Big)^{\frac{1}{r}}\|\phi\|_{L^p(0,T;E)}
\\ &\stackrel{\rm(ii)}{\le} CT^{\frac1q-\frac1p-\alpha}
\|\phi\|_{B^s_{p,r}(0,T;E)} + C T^{\frac1q-\frac1p-s-\alpha}
\|\phi\|_{L^p(0,T;E)}.
\end{aligned}\]
In (i) we used the triangle inequality in $L^r(0, T\wedge 1,\frac{du}{u})$ and in
(ii) we noted that $\alpha<\frac{1}{q} - \frac{1}{p}-s$.

Next,
\[\begin{aligned}
\Big(\int_{\R}
\Big\|\frac{\phi(t+h)\one_{(0,T)}(t+h)}{(t+h)^{\alpha}}
-\frac{\phi(t)\one_{(0,T)}(t)}{t^{\alpha}}\Big\|^q
\,dt\Big)^{\frac{1}{q}}
 &\le 2\Big(\int_{0}^T \Big\|\frac{\phi(t) }{t^{\alpha}}\Big\|^q
\,dt\Big)^{\frac{1}{q}}\\ & \le CT^{\frac1q-\frac1p -\alpha}
\|\phi\|_{L^p(0,T;E)}.
\end{aligned}\]
Using this we estimate the second part:
\[\begin{aligned}
\ & \Big(\int_{T\wedge 1}^1 u^{-sr} \sup_{|h|\le
u}\Big\|\frac{\phi(t+h)\one_{(0,T)}(t+h)}{(t+h)^{\alpha}}
-\frac{\phi(t) \one_{(0,T)}(t)}{t^{\alpha}}\Big\|_{L^q(I[h];E)}^r \,
\frac{du}{u}\Big)^{\frac{1}{r}}
\\ & \qquad \qquad\qquad \le CT^{\frac1q-\frac1p -\alpha}
\|\phi\|_{L^p(0,T;E)} \Big(\int_{T\wedge 1}^1
u^{-sr}\frac{du}{u}\Big)^{\frac{1}{r}}
\\ &  \qquad \qquad\qquad \le CT^{\frac1q-\frac1p-s-\alpha}
\|\phi\|_{L^p(0,T;E)}.
\end{aligned}\]
Putting everything together and using
 H\"older's inequality to estimate the $L^q$-norm
of $t^{-\alpha}\phi(t)$ we obtain
\[\begin{aligned}
\ & \|t\mapsto t^{-\alpha} \phi(t)\|_{B^s_{q,r}(0,T;E)}
\\ &  =
\|t\mapsto t^{-\alpha}\phi(t)\|_{L^q(0,T;E)}
\\ & \qquad  + \Big(\int_0^1
u^{-sr} \sup_{|h|\le
u}\Big\|\frac{\phi(t+h)\one_{(0,T)}(t+h)}{(t+h)^{\alpha}}
-\frac{\phi(t)\one_{(0,T)}(t)}{t^{\alpha}}\Big\|_{L^q(\R;E)}^r \,
\frac{du}{u}\Big)^{\frac{1}{r}}
\\ & \ \le  CT^{\frac1q-\frac1p -\alpha}
\|\phi\|_{L^p(0,T;E)} +  CT^{\frac1q-\frac1p-\alpha}
\|\phi\|_{B^s_{p,r}(0,T;E)} + CT^{\frac1q-\frac1p-s-\alpha}
\|\phi\|_{L^p(0,T;E)}.
\end{aligned}\]
\end{proof}

A Banach space $E$ has {\em type $p$}, where $p\in
[1,2]$, if there exists a constant $C\ge 0$ such that for all
$x_1,\dots,x_n\in E$ we have
\[ \Big(\E \Big\n \sum_{j=1}^n r_j\, x_j\Big\n^2\Big)^\frac12\le C
\Big(\sum_{j=1}^n \n x_j\n^p\Big)^\frac1p.\] Here $(r_j)_{j\ge 1}$
is a Rademacher sequence. Similarly $E$ has  {\em cotype $q$}, where
$q\in [2,\infty]$, if there exists a constant $C\ge 0$ such that for
all $x_1,\dots,x_n\in E$ we have
\[\Big(\sum_{j=1}^n \n x_j\n^q\Big)^\frac1q
\le C\Big(\E \Big\n \sum_{j=1}^n r_j\, x_j\Big\n^2\Big)^\frac12.\]
In these definitions the Rademacher variables may be replaced by
Gaussian variables without changing the definitions; for a proof and more
information see \cite{DJT}.
Every Banach space has type $1$ and cotype $\infty$, the spaces
$L^p(S)$, $1\le p<\infty$, have type $\min\{p,2\}$ and cotype
$\max\{p,2\}$, and Hilbert spaces have type $2$ and cotype $2$.
Every UMD space has nontrivial type, i.e., type $p$ for some $p\in
(1,2]$.

In view of the basic role of the space $\g(L^2(0,T;H),E)$
in the theory of vector-valued stochastic
integration, it is natural to look for conditions
on a function $\Phi:(0,T)\to \calL(H,E)$
ensuring that the associated integral operator $I_\Phi: L^2(0,T;H)\to E$,
\[
 I_\Phi f := \int_0^T \Phi(t)f(t)\,dt, \qquad f\in L^2(0,T;H),
\]
is well-defined and belongs to $\g(L^2(0,T;H),E)$. The next
proposition, taken from \cite{NVWco}, states such a condition for functions
$\Phi$ belonging to
suitable Besov spaces of $\g(H,E)$-valued functions.

\begin{lemma}\label{lem:Besovemb}
If $E$ has type $\tau \in [1,2)$, then $\Phi\mapsto I_\Phi$
defines a continuous embedding
\[B_{\tau,\tau}^{\frac{1}\tau-\frac{1}2}(0,T_0;\g(H,E))\embed  \g(L^2(0,T_0;H),E),\]
where the constant of the embedding depends on $T_0$ and the type
$\tau$ constant of $E$.
\end{lemma}

Conversely, if $\Phi\mapsto I_\Phi$ defines a continuous embedding
$B_{\tau,\tau}^{\frac{1}\tau-\frac{1}2}(0,T_0;\g(H,E))\embed
\g(L^2(0,T_0;H),E)$, then $E$ has type $\tau$ (see \cite{KNVW}); we
will not need this result.

\begin{lemma}\label{lem:gammaHolder}
Let $E$ be a Banach space with type $\tau\in [1, 2)$. Let $\alpha\ge
0$ and $q>2$ be such that $\alpha<\frac12-\frac1q$. There
exists a constant $C\ge 0$ such that for all $T\in [0, T_0]$
and $\Phi\in B_{q,\tau}^{\frac1\tau-\frac12}(0,T;\g(H,E))$ we have
\[
\begin{aligned}
&\sup_{t\in (0,T)}\|s\mapsto
(t-s)^{-\alpha}\Phi(s)\|_{\g(L^2(0,t;H),E)} \le C
T^{\frac{1}{2}-\frac1q-\alpha}
\|\Phi\|_{B_{q,\tau}^{\frac1\tau-\frac12}(0,T;\g(H,E))}.
\end{aligned}
\]
\end{lemma}

\begin{proof}
Fix $T\in [0,T_0]$ and $t\in [0,T]$. Then,
\[\begin{aligned}
 \| s\mapsto (t-s)^{-\alpha}\Phi(s)\|_{\g(L^2(0,t;H),E)}
  & = \n s\mapsto s^{-\alpha}\Phi(t-s)\|_{\g(L^2(0,t;H),E)}
\\ & = \| s\mapsto s^{-\alpha}\Phi(t-s)\one_{(0,t)}(s)\|_{\g(L^2(0,T_0;H),E)}
\\ & \stackrel{\rm(i)}{\le}
C\|s\mapsto s^{-\alpha}\Phi(t-s) \one_{(0,t)}(s)
\|_{B^{\frac1\tau-\frac12}_{\tau,\tau}(0,T_0;\g(H,E))}
\\ & \stackrel{\rm(ii)}{\le} C
t^{\frac12-\frac1q-\alpha} \|s\mapsto
\Phi(t-s)\|_{B^{\frac1\tau-\frac12}_{q,\tau}(0,t;\g(H,E))}
\\ &\le C
T^{\frac12-\frac1q-\alpha}
\|\Phi\|_{B^{\frac1\tau-\frac12}_{q,\tau}(0,T;\g(H,E))}.
\end{aligned}\]
In (i) we used Lemma \ref{lem:Besovemb} and (ii) follows from Lemma
\ref{lem:Besov}.
\end{proof}

In the remainder of this section we assume that $A$ is the
infinitesimal generator of an analytic $C_0$-semigroup $S = (S(t))_{t\ge 0}$
on $E$.
We fix an arbitrary  number $w\in\R$ such that
the semigroup generated by $A-w$ is uniformly exponentially stable.
The fractional powers $(w-A)^\eta$ are then well-defined, and for $\eta>0$
we put
\[ E_\eta := \D((w-A)^\eta).\]
This is a Banach space with respect to the norm
\[\n x\n_{E_\eta}:= \n x\n + \n (w-A)^\eta x\n.\]
As is well known, up to an equivalent norm this definition is independent of the
choice of $w$. The basic estimate
\begin{equation}\label{eq:Ja}
 \n S(t)\n_{\calL(E,E_\eta)} \le C t^{-\eta}, \qquad t\in [0,T_0],
\end{equation}
valid for $\eta>0$ with $C$ depending on $\eta$,
will be used frequently.

The extrapolation spaces $E_{-\eta}$ are defined, for $\eta>0$, as
the completion of $E$ with respect to the norm \[\n x\n_{E_{-\eta}}
:= \n (w-A)^{-\eta}x\n.\] Up to an equivalent norm, this space is
independent of the choice of $w$.

We observe at this point that the spaces $E_\eta$ and $E_{-\eta}$ inherit all
isomorphic Banach space properties of $E$, such as (co)type, the UMD property,
and property $(\Delta)$, via the isomorphisms $(w-A)^\eta:
E_\eta\simeq E$ and $(w-A)^{-\eta}:
E_{-\eta}\simeq E$.

The following lemma is well-known; a sketch of a proof is
included for the convenience of the reader.

\begin{lemma}\label{lem:detconv2}
Let  $q\in [1, \infty)$ and  $\tau\in [1,2)$ be given, and let $\eta\ge
0$ and $\theta\ge 0$ satisfy $\eta+\theta<\frac32-\frac1\tau$. There
exists a constant $C\ge 0$ such that for all $T\in [0,T_0]$ and
$\phi\in L^\infty(0,T;E_{-\theta})$ we have $S*\phi\in
B_{q,\tau}^{\frac1\tau-\frac12}(0,T;E_\eta)$ and
\[\|S*\phi\|_{B_{q,\tau}^{\frac1\tau-\frac12}(0,T;E_\eta)} \le C T^{\frac1q}
\|\phi\|_{L^\infty(0,T;E_{-\theta})}.
\]
\end{lemma}
\begin{proof}
Without loss of generality we may assume that $\eta,\theta>0$.
Let $\e>0$ be such that
$\eta+\theta < \frac32 - \frac1\tau-\e$. Then
\[\|S*\phi\|_{B_{q,\tau}^{\frac1\tau-\frac12}(0,T;E_\eta)}\leq C
T^{\frac1q}\|S*\phi\|_{C^{\frac1\tau-\frac12-\e}([0,T];E_\eta)} \leq C
T^{\frac1q} \|\phi\|_{L^\infty(0,T;E_{-\theta})}.\]
The first estimate is a direct consequence of the definition
of the Besov norm, and the second follows from \cite[Proposition 4.2.1]{Lun}.
\end{proof}

From the previous two lemmas we deduce the next convolution estimate.

\begin{proposition}\label{prop:detminconvgamma}
Let $E$ be a Banach space with type $\tau\in [1,2]$ and let  $0\le \alpha <
\frac12$. Let $\eta\ge 0$ and $\theta\geq 0$ satisfy
$\eta+\theta <\frac{3}{2}-\frac1\tau$. Then there is a constant
$C\ge 0$ such that for all $0\le t\le T\le T_0$ and $\phi\in
L^\infty(0,T;E)$,
\[
\|s\mapsto  (t-s)^{-\alpha} (S*\phi)(s)\|_{\g(L^2(0,t),E_\eta)}
 \le C T^{\frac12-\alpha}
\|\phi\|_{L^\infty(0,T;E_{-\theta})}.\]
\end{proposition}
\begin{proof}
First assume that $1\le \tau<2$. It follows from Lemmas
\ref{lem:gammaHolder} and \ref{lem:detconv2} that for any $q>2$ such
that $\a<\frac12-\frac1q$,
\[
\begin{aligned} \|s\mapsto  (t-s)^{-\alpha}S*\phi(s)\|_{\g(L^2(0,t),E_\eta)} &
\le CT^{\frac{1}{2}-\frac1q-\alpha}
\|S*\phi\|_{B_{q,\tau}^{\frac1\tau-\frac12}(0,T;E_\eta)}
\\ & \le C T^{\frac12-\alpha}
\|\phi\|_{L^\infty(0,T;E_{-\theta})}. \end{aligned}
\]

For $\tau=2$ we argue as follows. Since $E_\eta$ has type $2$, we
have a continuous embedding
$L^2(0,t;E_\eta)\hookrightarrow\g(L^2(0,t),E_\eta)$; see \cite{RS}. Therefore,
using \eqref{eq:Ja},
\[\begin{aligned}
 \|s\mapsto   (t-s)^{-\alpha}S*\phi(s)\|_{\g(L^2(0,t),E_\eta)}
&  \le C\|s\mapsto
(t-s)^{-\alpha} S*\phi(s)\|_{L^2(0,t;E_\eta)}
\\ & \le C
\|s\mapsto (t-s)^{-\alpha} \|_{L^{2}(0,t)} \|
S*\phi\|_{L^\infty(0,T;E_{\eta})}
\\ & \le C T^{\frac12-\alpha}
T^{1-\eta-\theta}\|\phi\|_{L^\infty(0,T;E_{-\theta})}.
\end{aligned}\]
\end{proof}

The following lemma, due to Da Prato, Kwapie{\'n} and Zabczyk \cite[Lemma
2]{DPKZ} in the Hilbert space case, gives a H\"older estimate for the convolution
\[R_{\alpha} \phi(t) := \frac{1}{\Gamma(\alpha)}\int_0^t (t-s)^{\a-1} S(t-s) \phi(s) \, ds. \]
The proof carries over to Banach spaces without change.

\begin{lemma}[\cite{DPKZ}]\label{lem:detconv}\label{fact2}
Let $0<\a\le 1$, $1<p< \infty$, $\lambda\ge 0$, $\eta\ge 0$,
and $\theta\geq 0$ satisfy
$\lambda+\eta+\theta < \a-\frac{1}{p}$.
Then there exist a constant $C\ge 0$ and an $\e>0$
such that for all $\phi\in L^p(0,T;E)$ and $T\in [0,T_0]$,
\[
\|R_\a \phi\|_{C^\lambda([0,T];E_{\eta})}
\le CT^\e\| \phi\|_{L^p(0,T;E_{-\theta})}.
\]
\end{lemma}

\section{Stochastic convolutions}
\label{sec:stoch_convol}

We now turn to the problem of estimating stochastic convolution integrals.
We start with a lemma which, in combination with Lemma \ref{lem:KW},
can be used to estimate stochastic convolutions involving
analytic semigroups.

\begin{lemma}\label{lem:derL1}
Let $S$ be an analytic $C_0$-semigroup on a Banach space $E$.
For all $0\le a <1$ and $\e>0$ the family
\[\big\{t^{a+\e} S(t)\in \calL(E, E_a): \ t\in [0,T]\big\}\]
is $R$-bounded in $ \calL(E, E_a)$,  with
$R$-bound of order $O(T^\e)$ as $T\downarrow 0$.
\end{lemma}
\begin{proof}
Let $N:[0,T]\to
\calL(E, E_a)$ be defined as $N(t) = t^{a+\e} S(t)$. Then $N$ is continuously
differentiable on $(0,T)$ and $N'(t) = (a+\e) t^{a+\e-1} S(t) + t^{a+\e} A
S(t)$, where $A$ is the generator of $S$.
Hence, by \eqref{eq:Ja},
\[\|N'(t)\|_{\calL(E,E_a)} \le C
t^{\e-1} \ \ \text{for $t\in (0,T)$}.\] By Lemma \ref{lem:int-der}
the $R$-bound on $[0,T]$ can now be bounded from above by
\[ \int_0^T \n N'(t)\n _{\calL(E,E_a)}\,dt \le C T^\e.\]
\end{proof}

We continue with an extension of the Da Prato-Kwapie\'n-Zabczyk factorization
method \cite{DPKZ} for Hilbert spaces to UMD spaces.
For deterministic $\Phi$, the assumption that $E$ is UMD can be
dropped. A related regularity result for arbitrary $C_0$-semigroups
is due to Millet and Smole\'nski \cite{MilSmo}.

It will be convenient to introduce the notation
\[S \diamond \Phi(t) := \int_0^t S(t-s) \Phi(s) \,  d W_H(s)\]
for the stochastic convolution with respect to $W_H$ of $S$ and
$\Phi$, where $W_H$ is an $H$-cylindrical Brownian motion.

\begin{proposition}\label{prop:path-hcont2}
Let $0<\alpha<\frac12$, $\lambda\ge 0$, $\eta\geq 0$, $\theta\ge 0$, and
$p>2$ satisfy $\lambda+\eta+\theta<\alpha-\frac1p$.
Let $A$ be the generator of an analytic $C_0$-semigroup $S$ on a UMD
space $E$ and let $\Phi:(0,T)\times\O\to \calL(H,E_{-\theta})$ be
$H$-strongly measurable and adapted. Then there exist
$\e>0$ and $C\ge 0$ such that
\[
\E\|S\diamond \Phi \|^p_{C^\lambda([0,T];E_{\eta})} \le C^p T^{\e p}
\int_0^T \E \|s\mapsto (t-s)^{-\alpha}
\Phi(s)\|_{\g(L^2(0,t;H),E_{-\theta})}^p \, dt.
\]
\end{proposition}

Here, and in similar formulations below, it is part of the
assumptions that the right-hand side is well-defined and finite.
In particular it follows from the proposition there exist
$\e>0$ and $C\ge 0$ such that
\[
\E\|S\diamond \Phi \|^p_{C^\lambda([0,T];E_{\eta})} \le C^p T^{\e p}
\sup_{t\in [0,T]} \E \|s\mapsto (t-s)^{-\alpha}
\Phi(s)\|_{\g(L^2(0,t;H),E_{-\theta})}^p \]
provided the right-hand
side is finite.

\begin{proof}
The idea of the proof is the same as in \cite{DPKZ}, but there are
some technical subtleties which justify us to outline the main
steps.

Let $\beta\in (0,\frac12)$ be such that $\lambda+ \eta<\beta-\frac1p
< \alpha-\theta-\frac1p$.  It follows from Lemmas \ref{lem:KW} and
\ref{lem:derL1} that, for almost all $t\in [0,T]$, almost surely we have
\begin{equation}\label{eq:convol-est}
\begin{aligned}
\ & \|s\mapsto (t-s)^{-\beta} S(t-s)\Phi(s)\|_{\g(L^2(0,t;H),E)}
\\ & \qquad\qquad
\le Ct^{\alpha-\beta-\theta} \|s\mapsto (t-s)^{-\alpha}
\Phi(s)\|_{\g(L^2(0,t;H),E_{-\theta})}.
\end{aligned}
\end{equation}
By Proposition \ref{prop:NVW}, the process
$\zeta_{\beta}:[0,T]\times\O\to E$,
\[\zeta_{\beta}(t) := \frac{1}{\Gamma(1-\beta)}
\int_0^t (t-s)^{-\beta} S(t-s)\Phi(s) \, dW_H(s),\] is well-defined
for almost all $t\in [0,T]$ and satisfies
\[
\big(\E\|\zeta_{\beta}(t) \|^p \big)^\frac1p \le
Ct^{\alpha-\beta-\theta}
\big(\E \|s\mapsto(t-s)^{-\alpha}
\Phi(s)\|_{\g(L^2(0,t;H),E_{-\theta})}^p\big)^\frac1p.\]
By Proposition \ref{prop:convprogr} the process $\zeta_\b$ is strongly
measurable. Therefore, by Fubini's theorem,
\[\begin{aligned}
\|\zeta_{\beta}\|_{L^p(\O;L^p(0,T;E))}
 & \le C T^{\alpha-\beta-\theta}\int_0^T \E \|s\mapsto
(t-s)^{-\alpha} \Phi(s)\|_{\g(L^2(0,t;H),E_{-\theta})}^p \, dt.
\end{aligned}
\]
By Lemma \ref{fact2}, the paths of $R_\beta \zeta_{\beta}$
belong to
$C^\lambda([0,T];E_{\eta})$
almost surely, and for some $\e'>0$ independent of $T\in [0,T_0]$ we have
\begin{equation}\label{eq:estRzeta1}
\begin{aligned}
\ & \|R_\beta \zeta_{\beta}\|_{L^p(\O;C^\lambda([0,T];E_{\eta}))}
\\ & \qquad\qquad \le CT^{\e'} \|\zeta_{\beta}\|_{L^p(\O;L^p(0,T;E))}
\\ & \qquad\qquad \le C T^{\alpha-\beta-\theta+\e'} \Big(\int_0^T \E \|s\mapsto
(t-s)^{-\alpha} \Phi(s)\|_{\g(L^2(0,t;H),E_{-\theta})}^p \,
dt\Big)^{\frac1p}.
\end{aligned}
\end{equation}

The right ideal property \eqref{eq:ideal}, \eqref{eq:convol-est},
and Proposition \ref{prop:NVW} imply the stochastic integrability of
$s\mapsto S(t-s) \Phi(s)$ for almost all $t\in [0,T]$. The proof will
be finished (with $\e = \a-\b-\theta+\e'$) by showing that
almost surely on $(0,T)\times \O$,
\[
S\diamond \Phi = R_\beta \zeta_{\beta}.
\]
It suffices to check that for almost all $t\in [0,T]$ and $x^* \in
E^*$ we have, almost surely,
\begin{equation}\label{eq:weakstochfub}
\lb S\diamond \Phi(t), x^*\rb = \frac{1}{\Gamma(\beta)}\int_0^t
(t-s)^{\beta-1} \lb S(t-s) \zeta_{\beta}(s), x^*\rb \, ds.
\end{equation}
This follows from a standard argument via the stochastic Fubini
theorem, cf. \cite{DPKZ}, which can be applied here since almost
surely we have, writing $\lb \Phi(r),x\s\rb := \Phi\s(r)x\s$,
\[\begin{aligned}
\int_0^t& \big\n\big<(t-s)^{\beta-1} S(t-s) (s-\cdot)^{-\beta}
S(s-\cdot)\Phi(\cdot) \one_{[0,s]}(\cdot), x^*\big>\big\|_{L^2(0,t;H)}\,ds
\\ &= \int_0^t\big\|\big<(s-\cdot)^{-\beta} S(s-\cdot)\Phi(\cdot),
(t-s)^{\beta-1} S^*(t-s) x^*\big>\big\|_{L^2(0,t;H)}\,ds
\\ &\le \int_0^t \|(s-\cdot)^{-\beta} S(s-\cdot)\Phi(\cdot)\|_{\g(L^2(0,t;H),E)}
\|(t-s)^{\beta-1} S^*(t-s) x^*\| \, ds,
\end{aligned}\]
which is finite for almost all $t\in [0,T]$ by H\"older's
inequality.
\end{proof}

\begin{remark}\label{rem:path-hcont2}
The stochastic integral $S\diamond \Phi$ in Proposition
\ref{prop:path-hcont2} may be defined only for almost all $t\in
[0,T]$. If in addition one assumes that $\Phi\in
L^p((0,T)\times\O;\g(H,E_{-\theta}))$,
then $S\diamond
\Phi(t)$ is well-defined in $E_{\eta}$ for all $t\in [0,T]$. This
follows readily from \eqref{eq:weakstochfub},
\cite[Theorem 3.6(2)]{NVW} and the density of $E^*$ in $(E_{\eta})^*$. Since we
will not need this in the sequel, we leave this to the interested
reader.
\end{remark}

As a consequence  we have the following regularity result of
stochastic convolutions in spaces with type $\tau\in [1, 2)$. We will
not need this result below, but we find it interesting enough to
state it separately.

\begin{corollary}\label{cor:convBesov}
Let $E$ be a UMD space with type $\tau\in [1, 2)$. Let $p>2$, $q>2$,
$\lambda\ge 0$, $\eta\geq 0$, $\theta\ge 0$ be such that
$\lambda+\eta+\theta< \frac12-\frac1p-\frac1q$.
Then there is an $\d>0$ such that for all
$H$-strongly strongly measurable and adapted $\Phi:(0,T)\times\O\to
\calL(H,E_{-\theta})$,
\begin{equation}\label{eq:Cldp2}
\begin{aligned}
\E\|  S\diamond \Phi\|^p_{C^\lambda([0,T];E_{\eta})} & \leq C^p T^{\d p}
\E\|\Phi\|_{B^{\frac{1}{\tau}-\frac12}_{q,\tau}(0,T;\g(H,E_{-\theta}))}^p.
\end{aligned}
\end{equation}
\end{corollary}

\begin{proof} By assumption we may choose $\a\in (0,\frac12)$ such that
$\lambda+\eta+\theta+\frac1p<\a<\frac12-\frac1q$. The result now
follows from Proposition \ref{prop:path-hcont2} and Lemma
\ref{lem:gammaHolder} (noting that $E_{-\theta}$ has type $\tau$):
\[\begin{aligned}
\E\| S\diamond \Phi
\|^p_{C^\lambda([0,T];E_{\eta})}
 & \leq C^p T^{\e p} \sup_{t\in [0,T]}\E\|s\mapsto
(t-s)^{-\alpha}  \Phi(s)\|_{\g(L^2(0,t;H),E_{-\theta})}^p
\\ & \leq C^p T^{(\frac12-\frac1q-\a+\e)p}\E\|\Phi\|_{B^{\frac{1}{\tau}-\frac12}_{q,\tau}(0,T;\g(H,E_{-\theta}))}^p.
\end{aligned}\]
\end{proof}

The main estimate of this section is contained in the next result.

\begin{proposition}\label{prop:gammanormconv}
Let $E$ be a UMD Banach space. Let $\eta\ge 0$, $\theta \ge 0$,
$\a>0$ satisfy $0\le \eta+\theta<\a<\frac12$. Let $\Phi:(0,T)\times
\O\to \calL(H,E_{-\theta})$ be adapted and $H$-strongly measurable. Then for
all $1<p<\infty$ and all $0\le t\le T\le T_0$,
\[
\E \| (t-\cdot)^{-\alpha} S\diamond
\Phi(\cdot)\|_{\g(L^2(0,t;H),E_\eta)}^p \le C^p
T^{(\frac12-\eta-\theta)p} \E \|(t-\cdot)^{-\alpha}
\Phi(\cdot)\|_{\g(L^2(0,t;H),E_{-\theta})}^p.
\]
\end{proposition}

\begin{proof} Fix $0\le t\le T\le T_0$.
As in Proposition \ref{prop:path-hcont2} one shows that the
finiteness of the right-hand side  implies that $s\mapsto S(t-s)
\Phi(s)$ is stochastically integrable on $[0,t]$. We claim that
$s\mapsto  S(t-s) \Phi(s)$ takes values in $E_{\eta}$ almost surely
and is stochastically integrable on $[0,t]$ as an $E_{\eta}$-valued
process. Indeed, let $\e>0$ be such that $\b:=\eta+\theta+\e<\alpha$
and put
\[N_\b(t):= t^\b (\mu-A)^{\eta+\theta} S(t).\]
It follows from Lemmas  \ref{lem:KW} and \ref{lem:derL1} that
\[\begin{aligned}
 \E\|
S(t-\cdot) \Phi(\cdot)\|_{\g(L^2(0,t;H),E_{\eta})}^p
&  \le C \E\|
N_\b(t-\cdot) (t-\cdot)^{-\b}
\Phi(\cdot)\|_{\g(L^2(0,t;H),E_{-\theta})}^p
\\ & \le C
T^{\e p}\E\|(t-\cdot)^{-\b}
\Phi(\cdot)\|_{\g(L^2(0,t;H),E_{-\theta})}^p,
\end{aligned}\]
and the expression on the right-hand side is finite by the
assumption. The stochastic
integrability now follows from Proposition \ref{prop:NVW}.
This proves the claim.
Moreover, by Proposition \ref{prop:convprogr}, the stochastic
convolution process $S\diamond \Phi$
is adapted and strongly measurable as an $E_\eta$-valued process.

Let $G^p(\O;E_{\eta})$ and $G^p(\O\times\widetilde\O;E_\eta)$ denote
the closed subspaces in $L^p(\O;E_{\eta})$ and
$L^p(\O\times\widetilde\O;E_\eta)$ spanned by all elements of the
form $\int_0^T\Psi\,dW_H$ and $\int_0^T\Psi\,d\widetilde W_H$,
respectively, where $\widetilde W_H$ is an independent copy of $W_H$
and $\Psi$ ranges over all adapted elements in
$L^p(\O;\g(L^2(0,T;H),E))$. Since $E_\eta$ is a UMD space, by
Proposition \ref{prop:NVW} the operator
\[D_p
\int_0^T\Psi\,d\widetilde W_H := \int_0^T\Psi\,d W_H,
\]
is well defined and bounded from $G^p(\O\times\widetilde\O;E_\eta)$
to $G^p(\O;E_\eta)$. Using the Fubini isomorphism of Lemma
\ref{lem:gamma-Fubini} twice, we estimate
\[\begin{aligned}
\ & \big\|s\mapsto (t-s)^{-\alpha}
S\diamond \Phi(s)\big\|_{L^p(\O;\g(L^2(0,t),E_{\eta}))}
\\ & \qquad \eqsim  \Big\|s\mapsto
\int_0^s (t-s)^{-\alpha} S(s-r) \Phi(r)\,dW_H(r) \Big\|_{\g(L^2(0,t),G^p(\O;E_\eta))}
\\ & \qquad =  \Big\|s\mapsto
D_p \int_0^t 1_{(0,s)}(r) (t-s)^{-\alpha} S(s-r) \Phi(r)\,d\widetilde W_H(r)
\Big\|_{\g(L^2(0,t),G^p(\O;E_{\eta}))}
\\ & \qquad \lesssim
 \Big\|s\mapsto
 \int_0^t 1_{(0,s)}(r) (t-s)^{-\alpha} S(s-r) \Phi(r)\,d\widetilde W_H(r)
\Big\|_{\g(L^2(0,t),G^p(\O\times \widetilde\O;E_{\eta}))}
\\ & \qquad \eqsim
\Big\|s\mapsto \int_0^s (t-s)^{-\alpha} S(s-r)\Phi(r)\,d\widetilde
W_H(r) \Big\|_{L^p(\O;\g(L^2(0,t),L^p(\widetilde \O;E_{\eta})))}.
\end{aligned}
\]
Rewriting the right-hand side in terms of the function $N_\beta(t) =
t^\b (\mu-A)^{\eta+\theta} S(t)$ introduced above and using the
stochastic Fubini theorem to interchange the Lebesgue integral and
the stochastic integral, the right-hand side can be estimated as
\[\begin{aligned}
\ & \Big\|s\mapsto \int_0^s (t-s)^{-\alpha} S(s-r)\Phi(r)\,d\widetilde
W_H(r) \Big\|_{L^p(\O;\g(L^2(0,t),L^p(\widetilde \O;E_{\eta})))}
\\ & \eqsim
\Big\|s\mapsto \int_0^s (t-s)^{-\alpha}(\mu-A)^{\eta+\theta} S(s-r)\Phi(r)\,d\widetilde
W_H(r) \Big\|_{L^p(\O;\g(L^2(0,t),L^p(\widetilde \O;E_{-\theta})))}
\\ & =
\Big\|s\mapsto \int_0^s (t-s)^{-\alpha}(s-r)^{-\beta} N(s-r)\Phi(r)\,d\widetilde
W_H(r) \Big\|_{L^p(\O;\g(L^2(0,t),L^p(\widetilde \O;E_{-\theta})))}
\\ & = \Big\|s\mapsto
\int_0^s (t-s)^{-\alpha}
\\ & \qquad\quad\times
(s-r)^{-\b} \int_0^{s-r}   N_\b'(w)
 \Phi(r)\,dw\,d\widetilde W_H(r)
\Big\|_{L^p(\O;\g(L^2(0,t),L^p(\widetilde \O;E_{-\theta})))}
\\ & = \Big\|s\mapsto
\int_0^s  N_\b'(w)
\\ & \qquad\quad\times \int_0^{s-w} (t-s)^{-\alpha} (s-r)^{-\b}
\Phi(r)\,d\widetilde W_H(r)\,dw
\Big\|_{L^p(\O;\g(L^2(0,t),L^p(\widetilde \O;E_{-\theta})))}
\\ & = \Big\|s\mapsto
\int_0^t  N_\b'(w) 1_{(0,s)}(w)
\\ & \qquad\quad\times   \E_{\widetilde\F_{s-w}}\int_0^{s} (t-s)^{-\alpha} (s-r)^{-\b}
 \Phi(r)\,d\widetilde W_H(r)\,dw
\Big\|_{L^p(\O;\g(L^2(0,t),L^p(\widetilde \O;E_{-\theta})))},
\end{aligned}
\]
where $\E_{\widetilde\F_t}(\xi):= \E(\xi|\widetilde\F_t)$ is the
conditional expectation with respect to $\widetilde \F_t =
\sigma(\widetilde W_H(s)h: \ 0\le s\le t, \ h\in H\}$. Next we note
that \[\int_0^{t} \n N_\b'(w) \n \,dw \lesssim T^\e.
\]
Applying Lemmas \ref{lem:Stein} and \ref{lem:KW} pointwise with
respect to $\om\in\O$, we may estimate the right-hand side above by
\[\begin{aligned}
 &  \int_0^t \n N_\b'(w)\n \Big\|
s\mapsto 1_{(w,t)}(s)
\\ & \qquad\qquad\times \E_{\widetilde\F_{s-w}}
\int_0^{s}(t-s)^{-\alpha}(s-r)^{-\b} \Phi(r)\,d\widetilde W_H(r)
\Big\|_{L^p(\O;\gamma(L^2(0,t),L^p(\widetilde \O;E_{-\theta})))} \,dw
\\ & \lesssim  T^{\e}
 \Big\|s\mapsto \E_{\widetilde\F_{s-w}}\int_0^s (t-s)^{-\alpha}
(s-r)^{-\b} \Phi(r)\,d\widetilde W_H(r)
\Big\|_{L^p(\O;\gamma(L^2(0,t),L^p(\widetilde \O;E_{-\theta})))}
\\ & \lesssim T^\e
 \Big\|s\mapsto \int_0^s (t-s)^{-\alpha}
(s-r)^{-\b} \Phi(r)\,d\widetilde W_H(r)
\Big\|_{L^p(\O;\gamma(L^2(0,t),L^p(\widetilde \O;E_{-\theta})))}
\\ & \lesssim T^\e
\big\|s\mapsto \big[r\mapsto (t-s)^{-\alpha} (s-r)^{-\b}
1_{(0,s)}(r)\Phi(r)\big]
\big\|_{L^p(\O;\g(L^2(0,t),\g(L^2(0,t;H),E_{-\theta})))}.
\end{aligned}
\]
Using the isometry
\[\g(H_1,\g(H_2,F)) \simeq \g(H_2,\g(H_1;F)),\]
and the Fubini isomorphism, the right hand side is equivalent
to
\[
\begin{aligned} \ & \eqsim  T^{\e} \big\|s\mapsto \big[r\mapsto
(t-s)^{-\alpha} (s-r)^{-\b} 1_{(0,s)}(r)\Phi(r)\big]
\big\|_{L^p(\O;\g(L^2(0,t),\g(L^2(0,t;H),E_{-\theta})))}
\\ & \eqsim  T^{\e}
\big\|r\mapsto \!\big[s\mapsto (t-s)^{-\alpha} (s-r)^{-\b}1_{(0,s)}(r)
\Phi(r)\big]
\big\|_{L^p(\O;\g(L^2(0,t;H),\g(L^2(0,t),E_{-\theta})))}.
\end{aligned}\]

To proceed further we want to apply, pointwise with respect to $\O$,
Lemma \ref{lem:KW} to the multiplier
\[M: (0,t)\to \calL(E_{-\theta},\g(L^2(0,t),E_{-\theta}))\]
defined by
\[ M(r)x:= f_{r,t}\otimes x, \quad s\in (0,t), \ \ x\in
E_{-\theta},\] where $f_{r,t}\in L^2(0,t)$ is the function
\[ f_{r,t} (s):= (t-r)^{\alpha}(t-s)^{-\alpha}  (s-r)^{-\b}
\one_{(r,t)}(s).\] We need to check that the range of $M$ is
$\g$-bounded in $\calL(E_{-\theta},\g(L^2(0,t),E_{-\theta}))$. For
this we invoke Lemma \ref{lem:KaiWe}, keeping in mind that
$R$-bounded families are always $\g$-bounded and that UMD spaces
have finite cotype. To apply the lemma we check that functions
$f_{s,t}$ are uniformly bounded in $L^2(0,t)$:
\[\begin{aligned}
\int_0^t |f_{r,t}(s)|^2\,ds & = (t-r)^{2\alpha}  \int_r^t
(t-s)^{-2\alpha}  (s-r)^{-2\b}\, ds
\\ &= (t-r)^{1-2\b}  \int_0^{1} (1-u)^{-2\alpha}  u^{-2\b}\, du
\\ &\le T^{1-2\b}  \int_0^{1} (1-u)^{-2\alpha}  u^{-2\b}\, du.
\end{aligned}\]

It follows from Lemma \ref{lem:KW} that
\[\begin{aligned}
\ & \Big\|s\mapsto (t-s)^{-\alpha} (s-\cdot)^{-\b}1_{(0,s)}(\cdot)
\Phi(\cdot) \Big\|_{L^p(\O;\g(L^2(0,t;H),\g(L^2(0,t),E_{-\theta})))}
\\ & \qquad \le C T^{\frac12-\beta}
\Big\|r\mapsto (t-r)^{-\alpha} \Phi(r)
\Big\|_{L^p(\O;\g(L^2(0,t;H),E_{-\theta}))}
\\ & \qquad = C T^{\frac12 - \eta - \theta-\e}
\Big\|r\mapsto (t-r)^{-\alpha} \Phi(r)
\Big\|_{L^p(\O;\g(L^2(0,t;H),E_{-\theta}))}.
\end{aligned}\]

Combining all estimates we obtain the result.
\end{proof}

\section{$L_\g^2$-Lipschitz functions}
\label{sec:g-Lipschitz}

Let
$(S,\Sigma)$ be a countably generated measurable space and let $\mu$
be a finite measure on $(S,\mu)$. Then $L^2(S,\mu)$ is separable and we may define
 \[L^2_{\g}(S,\mu;E) :=
\g(L^2(S,\mu);E)\cap L^2(S,\mu;E).\] Here, $\g(L^2(S,\mu);E)\cap
L^2(S,\mu;E)$ denotes the Banach space of all strongly $\mu$-measurable functions
$\phi:S\to E$ for which
\[\|\phi\|_{L^2_{\g}(S,\mu;E)} := \|\phi\|_{\g(L^2(S,\mu);E)} + \|\phi\|_{L^2(S,\mu;E)}\]
is finite. One easily checks that the simple functions are dense in
$L^2_{\g}(S,\mu;E)$.

Next let $H$ be a nonzero separable Hilbert space,
let $E_1$ and $E_2$ be Banach spaces, and let
$f:S\times E_1\to \calL(H,E_2)$ be a function such that for all $x\in E_1$
we have $f(\cdot, x)\in \g(L^2(S,\mu;H),E_2)$. For simple functions $\phi:S\to
E_1$ one easily checks that $s\mapsto f(s, \phi(s))\in \g(L^2(S,\mu;H),E_2)$.
We call $f$
{\em $\gL$-Lipschitz function with respect to $\mu$}
if $f$ is strongly continuous in the
second variable and for all simple functions $\phi_1, \phi_2:S\to E_1$,
\begin{equation}\label{eq:glipsimp}
\|f(\cdot,\phi_1) - f(\cdot, \phi_2)\|_{\g(L^2(S,\mu;H),E_2)}\le
C\|\phi_1-\phi_2\|_{L^2_{\g}(S,\mu;E_1)}.
\end{equation}
In this case the mapping $\phi\mapsto S_{\mu,f} \phi := f(\cdot,\phi(\cdot))$
extends uniquely to a Lipschitz
mapping from $L^2_{\g}(S,\mu;E_1)$ into $\g(L^2(S,\mu;H),E_2)$. Its
Lipschitz constant will be denoted by $L^{\g}_{\mu,f}$.

It is evident from the definitions that for simple functions $\phi:S\to E_1$,
the operator $S_f(\phi)\in \g(L^2(S,\mu;H), E_2)$ is
represented by the function $f(\cdot,\phi(\cdot))$
The next lemma extends this to arbitrary functions $\phi\in L^2_\g(S,\mu;E_1)$.

\begin{lemma} If $f:S\times E_1\to \calL(H,E_2)$ is an
$L^2_{\g}$-Lipschitz function, then for all $\phi\in L^2_\g(S,\mu;E_1)$
the operator $S_{\mu,f}\phi\in \g(L^2(S,\mu;H), E_2)$ is
represented by the function $f(\cdot, \phi(\cdot))$.
\end{lemma}
\begin{proof}
Let $(\phi_n)_{n\geq 1}$ be a sequence of simple functions such that
$\phi = \limn \phi_n$ in $L^2_{\g}(S,\mu;E_1)$. We may assume that
$\phi = \limn \phi_n$ $\mu$-almost everywhere. It follows from
\eqref{eq:glipsimp} that $(f(\cdot, \phi_n(\cdot)))_{n\geq 1}$ is a Cauchy
sequence in $\g(L^2(S,\mu;H),E_2)$. Let $R\in \g(L^2(S,\mu;H),E_2)$
be its limit. We must show that $R$ is represented by $f(\cdot,
\phi(\cdot))$. Let $x^*\in E^*_2$ be arbitrary. Since $R^*x^* = \limn
f^*(\cdot, \phi_n(\cdot))x^*$ in $L^2(S,\mu;H)$ we may choose a subsequence
$(n_k)_{k\geq 1}$ such that $R^*x^* = \lim_{k\to \infty} f^*(\cdot,
\phi_{n_k}(\cdot)) x^*$ $\mu$-almost everywhere. On the other hand since $f$ is
strongly continuous in the second variable we have
\[\lim_{k\to \infty} f^*(s, \phi_{n_k}(s)) x^* = f^*(s, \phi(s)) x^*\ \
\text{for $\mu$-almost all $s\in S$}.\]
This proves that for all $h\in H$ we have $R^*x^* = f^*(\cdot,
\phi(\cdot)) x^*$ $\mu$-almost everywhere and the result follows.
\end{proof}

Justified by this lemma, in what follows we shall always identify
$S_{\mu,f}\phi$ with $f(\cdot, \phi(\cdot))$.

If $f$ is $\gL$-Lipschitz with respect to all finite measures $\mu$ on $(S,\Sigma)$ and
\[L^\g_{f} := \sup\{L^{\g}_{\mu,f}: \mu \ \text{is a finite measure on
$(S,\Sigma)$}\}\] is finite, we say that $f$ is a {\em $\gL$-Lipschitz function}.
In type $2$ spaces there is the following easy criterium to check whether a function is
$\gL$-Lipschitz.

\begin{lemma}\label{lem:type2Lipschitz}
Let $E_2$ have type $2$. Let $f:S\times E_1\to \g(H,E_2)$ be such that
for all $x\in E_1$, $f(\cdot,x)$ is strongly measurable. If there is a
constant $C$ such that
\begin{eqnarray}
\label{eq:type2bdd} \|f(s,x)\|_{\g(H,E_2)} & \le & C(1+\|x\|), \ s\in S, \ \ x\in
E_1,
\\ \label{eq:type2Lipschitz} \|f(s,x) - f(s,y)\|_{\g(H,E_2)}&\le& C\|x-y\|, \
s\in S, \ \ x,y\in
E_1,
\end{eqnarray}
then $f$ is a $\gL$-Lipschitz function and $L_{f}^\g\le C_2 C$, where $C_2$
is the Rademacher type $2$ constant of $E_2$. Moreover, it satisfies the
following linear growth condition
\[\|f(\cdot,\phi)\|_{\g(L^2(S,\mu;H),E_2)} \le C_2 C(1+\|\phi\|_{L^2(S,\mu;E_1)}).\]
\end{lemma}
If $f$ does not depend on $S$, one can check that \eqref{eq:glipsimp} implies
\eqref{eq:type2bdd} and \eqref{eq:type2Lipschitz}.
\begin{proof}
Let $\phi_1, \phi_2\in L^2(S,\mu;E_1)$. Via an approximation
argument and \eqref{eq:type2Lipschitz} one easily checks that
$f(\cdot,\phi_1)$ and $f(\cdot,\phi_2)$ are strongly measurable. It
follows from \eqref{eq:type2bdd} that $f(\cdot,\phi_1)$ and
$f(\cdot,\phi_2)$ are in $L^2(S,\mu;\g(H,E_2))$ and from
\eqref{eq:type2Lipschitz} we obtain
\begin{equation}\label{eq:type2Lipf}
\|f(\cdot,\phi_1) - f(\cdot,\phi_2)\|_{L^2(S,\mu;\g(H,E_2))}\le
C\|\phi_1-\phi_2\|_{L^2(S,\mu;E_1)}.
\end{equation}
Recall from \cite{vNWe} that $L^2(S,\mu;\g(H,E_1))\hookrightarrow
\g(L^2(S,\mu;H),E_1)$ where the norm of the embedding equals $C_2$.
From this and \eqref{eq:type2Lipf} we conclude that
\[\|f(\cdot,\phi_1) - f(\cdot,\phi_2)\|_{\g(L^2(S,\mu;H),E_2)}\le C_2 C \|\phi_1-\phi_2\|_{L^2(S,\mu;E_1)}.\]
This clearly implies the result. The second statement follows in the same way.
\end{proof}

A function $f:E_1\to \calL(H,E_2)$ is said to be {\em
$\gL$-Lipschitz} if the induced function $\tilde f:S\times E_1\to
\calL(H,E_2)$, defined by $\tilde f(s,x) = f(x)$, is $\gL$-Lipschitz
for every finite measure space $(S,\Sigma,\mu)$.

\begin{lemma}\label{lem:gammalipf}
For a function $f:E_1\to
\calL(H,E_2)$, the following assertions are equivalent:
\begin{enumerate}
\item $f$ is $\gL$-Lipschitz;
\item There is a constant $C$ such that for some (and then for every) orthonormal basis
$(h_m)_{m\geq 1}$ of $H$ and all finite sequences $(x_n)_{n=1}^N,
(y_n)_{n=1}^N$ in $E_1$ we have
\[
\begin{aligned}
\E\Big\|\sum_{n=1}^N \sum_{m\geq 1} & \g_{nm} (f(x_n) h_m - f(y_n) h_m)
\Big\|^2 \\ & \le C^2 \E\Big\|\sum_{n=1}^N \g_n (x_n - y_n) \Big\|^2 + C^2
\sum_{n=1}^N \|x_n - y_n\|^2.
\end{aligned}
\]
\end{enumerate}
\end{lemma}
\begin{proof} (1) $\Rightarrow$ (2): \ Let $(h_m)_{m\geq 1}$ be an orthonormal basis and let
$(x_n)_{n=1}^N$ and $(y_n)_{n=1}^N$ in $E_1$ be arbitrary. Take
$S=(0,1)$ and $\mu$ the Lebesgue measure and choose disjoint sets
$(S_n)_{n=1}^N$ in $(0,1)$ such that $ \mu(S_n)=\frac1N$ for all
$n=1, \ldots, N$. Now define $\phi_1 := \sum_{n=1}^N \one_{S_n}
\otimes x_n$ and $\phi_2 := \sum_{n=1}^N \one_{S_n} \otimes y_n$.
Then (2) follows from \eqref{eq:glipsimp}.

(2) $\Rightarrow$ (1): \ Since the distribution of Gaussian vectors
is invariant under orthogonal transformations, if (2) holds for one
orthonormal basis $(h_m)_{m\geq 1}$, then it holds for every
orthonormal basis $(h_m)_{n\ge 1}$. By a well-known argument (cf.
\cite[Proposition 1]{Ja}), (2) implies that for all $(a_n)_{n=1}^N$
in $\R$ we have
\[
\begin{aligned}
\E\Big\|\sum_{n=1}^N \sum_{m\geq 1} & a_n \g_{nm} (f(x_n) h_m - f(y_n) h_m)
\Big\|^2 \\ & \le C^2 \E\Big\|\sum_{n=1}^N a_n \g_n (x_n - y_n) \Big\|^2 +
C^2 \sum_{n=1}^N a_n^2 \|x_n - y_n\|^2.
\end{aligned}
\]
Now \eqref{eq:glipsimp} follows for simple functions $\phi$, and the general case follows from this
by an approximation argument.
\end{proof}

Clearly, every $\gL$-Lipschitz function $f:E_1\to \g(H,E_2)$ is a
Lipschitz function. It is a natural question whether Lipschitz
functions are automatically $\gL$-Lipschitz. Unfortunately, this is
not true. It follows from the proof of \cite[Theorem 1]{NVgammaLip}
that if $\text{dim}(H)\geq 1$, then every Lipschitz function
$f:E_1\to \g(H,E_2)$ is $\gL$-Lipschitz if and only if $E_2$ has
type $2$.

A Banach space $E$ has {\em property $(\a)$} if
for all $N\ge 1$ and all sequences $(x_{mn})_{m,n=1}^N$ in $E$ we have
\[
\begin{aligned}
 \E \Bigl\n   \sum_{m,n=1}^N r_{mn} x_{m n} \Bigr\n^2
\eqsim \E'\E'' \Bigl\n  \sum_{m,n=1}^N r_{m}'r_n''x_{m n} \Bigr\n^2.
\end{aligned}
\]
Here, $(r_{mn})_{m,n\ge 1}$, $(r_{m}')_{m\ge 1}$, and
$(r_{n}'')_{n\ge 1}$ are Rademacher sequences, the latter two
independent of each other. By a randomization argument one can show
that the Rademacher random variables can be replaced by Gaussian
random variables. It can be shown using the Kahane-Khintchine
inequalities that the exponent $2$ in the definition can be replaced
by any number $1\le p<\infty$.

Property $(\a)$ has been introduced by Pisier
\cite{Pialpha}. Examples of spaces with this property are the Hilbert spaces
and the spaces $L^p$ for $1\le p<\infty$.

The next lemma follows directly from the definition of property $(\alpha)$ and
Lemma \ref{lem:gammalipf}.

\begin{lemma}\label{lem:alphaGammaLipschitz}
Let $E_2$ be a space with property $(\alpha)$. Then $f:E_1\to \g(H,E_2)$ is
$\gL$-Lipschitz if and only if there exists a constant $C$ such that
for all finite sequences $(x_n)_{n=1}^N$ and $(y_n)_{n=1}^N$ in $E_1$ we have
\[
 \E\Big\|\sum_{n=1}^N \g_{n} (f(x_n) - f(y_n))
\Big\|^2_{\g(H,E_2)}\le C^2 \E\Big\|\sum_{n=1}^N \g_n (x_n - y_n) \Big\|^2 +
C^2 \sum_{n=1}^N \|x_n - y_n\|^2.
\]
In particular, every $f\in \calL(E_1, \g(H,E_2))$ is $\gL$-Lipschitz.
\end{lemma}
When $H$ is finite dimensional, this result remains valid even if $E_2$ fails to
have
property $(\alpha)$.

The next example identifies an important class of $\gL$-Lipschitz continuous
functions.

\begin{example}[Nemytskii maps]
Fix $p\in [1, \infty)$ and
let $(S,\Sigma, \mu)$ be a $\sigma$-finite measure space. Let $b:\R\to \R$ be a
Lipschitz function; in case $\mu(S)=\infty$ we also assume that $b(0)=0$. Define the
Nemytskii map $B:L^p(S)\to L^p(S)$ by
$B(x)(s) := b(x(s)).$ Then $B$ is $\gL$-Lipschitz with respect to $\mu$. Indeed, it
follows from the Kahane-Khintchine inequalities that
\[\begin{aligned}
\Big(\E\Big\|\sum_{n=1}^N \g_n (B(x_n) - B(y_n))\Big\|^2\Big)^{\frac12}
&\eqsim_{p} \Big(\int_S \Big(\sum_{n=1}^N |b(x_n(s)) -
b(y_n(s))|^2\Big)^{\frac{p}{2}} \, d\mu(s)\Big)^{\frac1p}
\\ & \le L_b \Big(\int_S \Big(\sum_{n=1}^N |x_n(s) - y_n(s)|^2\Big)^{\frac{p}{2}} \,
d\mu(s)\Big)^{\frac1p}
\\ & \eqsim_p L_b\Big(\E\Big\|\sum_{n=1}^N \g_n (x_n - y_n)\Big\|^2\Big)^{\frac12}.
\end{aligned}\]
Now we apply Lemma \ref{lem:gammalipf}.
\end{example}

\section{Stochastic evolution equations I: integrable initial values}
\label{sec:main_results}

On the space $E$ we consider the stochastic equation:
\begin{equation}\tag{SCP}\label{SE}
\left\{\begin{aligned}
dU(t) & = (AU(t) + F(t,U(t)))\,dt + B(t,U(t))\,dW_H(t), \qquad t\in [0,\Tend],\\
 U(0) & = u_0,
\end{aligned}
\right.
\end{equation}
where $W_H$ is an $H$-cylindrical Brownian motion.
We make the following assumptions on $A$, $F$, $B$, $u_0$, the numbers
$\eta, \theta_F, \theta_B\ge 0$:
\let\ALTERWERTA\theenumi
\let\ALTERWERTB\labelenumi
\def\theenumi{(A1)}
\def\labelenumi{(A1)}
\begin{enumerate}
\item \label{as:semigroup}
The operator $A$ is the generator of an analytic $C_0$-semigroup $S$ on a UMD
Banach space $E$.
\end{enumerate}
\let\theenumi\ALTERWERTA
\let\labelenumi\ALTERWERTB

\let\ALTERWERTA\theenumi
\let\ALTERWERTB\labelenumi
\def\theenumi{(A2)}
\def\labelenumi{(A2)}
\begin{enumerate}
\item \label{as:LipschitzF}
The function
\[
F:[0,\Tend]\times\O\times E_{\eta}\to E_{-\theta_F}\] is Lipschitz
of linear growth uniformly in $[0,\Tend]\times\O$, i.e., there are
constants $L_F$ and $C_F$ such that for all $t\in [0,\Tend],
\omega\in \O \ \text{and} \ x,y\in E_\eta$,
\begin{eqnarray*}
\|(F(t,\omega,x) - F(t,\omega,y))\|_{E_{-\theta_F}} & \le
L_F\|x-y\|_{E_\eta},
\\
\| F(t,\omega,x)\|_{E_{-\theta_F}} & \le  C_F(1+\|x\|_{E_\eta}).
\end{eqnarray*}
Moreover, for all $x\in E_\eta$, $(t,\omega)\mapsto F(t,\omega, x)$
is strongly measurable and adapted in $E_{-\theta_F}$.
\end{enumerate}
\let\theenumi\ALTERWERTA
\let\labelenumi\ALTERWERTB

\let\ALTERWERTA\theenumi
\let\ALTERWERTB\labelenumi
\def\theenumi{(A3)}
\def\labelenumi{(A3)}
\begin{enumerate}
\item \label{as:LipschitzB}
The function
\[
B:[0,\Tend]\times\O\times E_\eta\to \calL(H,E_{-\theta_B})\] is
$\gL$-Lipschitz of linear growth uniformly in $\O$, i.e.,
there are constants $L_B^{\g}$ and $C_B^{\g}$ such that for all
finite measures $\mu$ on $([0,\Tend], \mathcal{B}_{[0,\Tend]})$, for
all $\omega\in \O$, and all $\phi_1, \phi_2\in
L^2_{\g}((0,\Tend),\mu;E_\eta)$,
\[
\begin{aligned}
\|
(B(\cdot,\omega,\phi_1) -
\ & B(\cdot,\omega,\phi_2))\|_{\g(L^2((0,\Tend),\mu;H),E_{-\theta_B})}
\\ &\ \ \qquad\qquad\qquad \le
L_B^{\g}\|\phi_1-\phi_2\|_{L^2_{\g}((0,\Tend),\mu;E_\eta)},
\end{aligned}
\]
and
\[
B(\cdot,\omega,\phi)\|_{\g(L^2((0,\Tend),\mu;H),E_{-\theta_B})}
\le C_B^{\g}(1+\|\phi\|_{L^2_{\g}((0,\Tend),\mu;E_\eta)}).
\]
Moreover, for all $x\in E_\eta$, $(t,\omega)\mapsto
B(t,\omega,
x)$ is $H$-strongly measurable and adapted in $E_{-\theta_B}$.
\end{enumerate}
\let\theenumi\ALTERWERTA
\let\labelenumi\ALTERWERTB

\let\ALTERWERTA\theenumi
\let\ALTERWERTB\labelenumi
\def\theenumi{(A4)}
\def\labelenumi{(A4)}
\begin{enumerate}
\item \label{as:initial_value}
The initial value $u_0:\O\to E_\eta$ is strongly $\F_0$-measurable.
\end{enumerate}
\let\theenumi\ALTERWERTA
\let\labelenumi\ALTERWERTB

We call a process $(U(t))_{t\in [0,\Tend]}$
a {\em mild $E_\eta$-solution} of \eqref{SE} if
\begin{enumerate}
\item[(i)] $U:[0,T_0]\times\O\to E_\eta$ is strongly measurable and adapted,
\item[(ii)] for all $t\in [0,T_0]$, $s\mapsto S(t-s) F(s,U(s))$ is in $L^0(\O;L^1(0,t;E))$,
\item[(iii)] for all $t\in [0,T_0]$, $s\mapsto S(t-s) B(s,U(s))$ $H$-strongly measurable
and adapted and in $\g(L^2(0,t;H),E)$ almost surely,
\item[(iv)] for all $t\in [0,T_0]$, almost surely
\[U(t) = S(t) u_0 + S*F(\cdot,U)(t) + S\diamond B(\cdot,U)(t).\]
\end{enumerate}

By (ii) the deterministic
convolution is defined pathwise as a Bochner integral,
and since $E$ is a UMD space, by (iii) and Proposition
\ref{prop:NVW} the
stochastic convolutions is well-defined.

We shall prove an existence and uniqueness result
for \eqref{SE} using a fixed point argument in a suitable scale of Banach
spaces of $E$-valued processes introduced next.
Fix $T\in (0,\Tend]$, $p\in [1, \infty)$, $\alpha\in
(0,\frac12)$.
We define $V_{\a,\infty}^{p}([0,T]\times\O;E)$ as the space of all
continuous adapted processes $\phi:[0,T]\times\O\to E$ for which
\[\begin{aligned}
\ & \|\phi\|_{V_{\a,\infty}^{p}([0,T]\times\O;E)}
\\ & \qquad := \big(\E\|\phi\|_{C([0,T];E)}^p\big)^\frac1p +
\sup_{t\in [0,T]}\Big(\E \|s\mapsto
(t-s)^{-\alpha}\phi(s)\|_{\g(L^2(0,t),E)}^p\Big)^{\frac{1}{p}}
\end{aligned}\] is finite.
Similarly we define $V_{\a,p}^{p}([0,T]\times\O;E)$ as the space of
pathwise continuous and adapted processes $\phi:[0,T]\times\O\to E$
for which
\[\begin{aligned}
\ & \|\phi\|_{V_{\a,p}^{p}([0,T]\times\O;E)}
\\ & \qquad := \big(\E\|\phi\|_{C([0,T];E)}^p\big)^\frac1p +
\Big(\int_0^T \E \|s\mapsto
(t-s)^{-\alpha}\phi(s)\|_{\g(L^2(0,t),E)}^p \, dt\Big)^{\frac{1}{p}}
\end{aligned}\] is finite.
 Identifying processes which are indistinguishable,
the above norm on $V_{\a,p}^{p}([0,T]\times\O;E)$ and $V_{\a,\infty}^{p}([0,T]\times\O;E)$
turn these spaces into Banach spaces.

The main result of this section, Theorem \ref{thm:mainexistenceL} below,
establishes existence and
uniqueness of a mild solution of \eqref{SE} with initial value
$u_0\in L^p(\O,\F_0;E_\eta)$ in each of the spaces
$V_{\a,p}^{p}([0,\Tend]\times\O;E)$ and $V_{\a,\infty}^{p}([0,\Tend]\times\O;E)$.
Since we have a continuous embedding
$V_{\a,\infty}^{p}([0,\Tend]\times\O;E)\embed
V_{\a,p}^{p}([0,\Tend]\times\O;E)$, the existence result is stronger
for $V_{\a,\infty}^{p}([0,\Tend]\times\O;E)$ while the uniqueness
result is stronger for $V_{\a,p}^{p}([0,\Tend]\times\O;E)$.

For technical reasons, in the next section we will also need the space
${\tilde V}_{\a,p}^{p}([0,T]\times\O;E)$
which is obtained by `pathwise continuous' replaced by
`pathwise bounded and $\mathcal{B}_{[0,T]}\otimes\F$-measurable' and
 $C([0,T];E)$ replaced by $B_{\rm b}([0,T];E)$ in the definition of
 ${\tilde V}_{\a,p}^{p}([0,T]\times\O;E)$. Here $B_{\rm b}([0,T];E)$
denotes the Banach space of bounded strongly Borel measurable functions on
$[0,T]$ with values in
$E$, endowed with the supremum norm.

Consider the fixed point operator
\[L_T (\phi)= \big[t\mapsto S(t) u_0 + S*F(\cdot,\phi)(t) + S\diamond B(\cdot,\phi)(t)\big].\]
In the next proposition we show that
$L_T$ is well-defined on each of the three spaces introduced above
and that it is a strict contraction for $T$ small enough.

\begin{proposition}\label{prop:contr} Let $E$ be a UMD space with
type $\tau\in [1, 2]$. Suppose that {\rm
\ref{as:semigroup}-\ref{as:initial_value}} are satisfied and assume
that $0\le \eta+\theta_F<\frac32 -\frac1\tau$ and $0\le
\eta+\theta_B<\frac12$. Let $p>2$ and $\alpha\in (0,\frac12)$ be
such that $\eta+\theta_B<\a-\frac1p$. If $u_0\in L^p(\O;E_\eta)$,
then the operator $L_T$ is well-defined and bounded on each of the
spaces
\[V \in \big\{ V_{\a,\infty}^p([0,T]\times\O;E_\eta), \
V_{\a,p}^p([0,T]\times\O;E_\eta), \ \tilde
V_{\a,p}^p([0,T]\times\O;E_\eta)\big\},\] and there exist a constant
$C_T$, with $\lim_{T\downarrow 0} C_T =0$, such that for all
$\phi_1, \phi_2\in V$,
\begin{equation}\label{eq:fixpointest}
\|L_T(\phi_1) - L_T(\phi_2)\|_{V} \le C_T \|\phi_1 - \phi_2\|_{V}.
\end{equation} Moreover, there is a constant
$C\ge 0$, independent of $u_0$, such that for all $\phi\in V$,
\begin{equation}\label{eq:operatorest}
\|L_T(\phi)\|_{V}\le C(1+ (\E\|u_0\|_{E_\eta}^p)^{\frac1p}) + C_T \n \phi\n_{V}.
\end{equation}
\end{proposition}
\begin{proof}
We give a detailed proof for the space
$V_{\a,\infty}^p([0,T]\times\O;E_\eta)$. The proof for
$V_{\a,p}^p([0,T]\times\O;E_\eta)$ is entirely similar. For the
proof for $\tilde V_{\a,p}^p([0,T]\times\O;E_\eta)$ one replaces
$C([0,T];E)$ by $B_{\rm b}((0,T);E)$.

{\em Step 1: Estimating the initial value part.}\ Let $\e\in (0,\frac12)$. From
Lemmas  \ref{lem:KW} and
\ref{lem:derL1} we infer that
\[\begin{aligned}
\|s\mapsto (t-s)^{-\alpha} S(s) u_0\|_{\g(L^2(0,t),E_\eta)} &\le C
\|s\mapsto (t-s)^{-\alpha} s^{-\e} u_0\|_{\g(L^2(0,t),E_\eta)}
\\ & = C \|s\mapsto
(t-s)^{-\alpha}s^{-\e}\|_{L^2(0,t)} \|u_0\|_{E_\eta}
\\ &  \leq C \|u_0\|_{E_\eta}.
\end{aligned}\]
For the other part of the $\V$-norm we note that
\[\|Su_0\|_{C([0,T];E_\eta)} \le C
\|u_0\|_{E_\eta}.\]
It follows that
\[\|S u_0\|_{\V} \le C
\|u_0\|_{L^p(\O;E_\eta)}.\]

{\em Step 2: Estimating the deterministic convolution.}\
We proceed in two steps.

(a): \ For $\psi\in C([0,T] ;E_{-\theta_F})$ we estimate the
$\V$-norm of $S*\psi$.

By Lemma \ref{lem:detconv} (applied with $\a=1$ and $\l=0$) $S*\psi$ is
continuous in $E_\eta$. Using \eqref{eq:Ja} we estimate:
\begin{equation}\label{eq:detphi2}
\begin{aligned}
\|S*\psi\|_{C([0,T];E_\eta)} &\le C \int_0^t
(t-s)^{-\eta-\theta_F}\, ds \, \n\psi\n_{C([0,T];E_{-\theta_F})}
\\ & \le CT^{1-\eta-\theta_F} \|\psi\|_{C([0,T];E_{-\theta_F})}.
\end{aligned}
\end{equation}

Also, since $E$ has type $\tau$, it follows from Proposition
\ref{prop:detminconvgamma} that
\begin{equation}\label{eq:detphi1}
\begin{aligned} \|s\mapsto
(t-s)^{-\alpha} S*\psi(s)\|_{\g(L^2(0,t),E_\eta)}
 & \le T^{\frac12-\a}\|\psi\|_{C([0,T];E_{-\theta_F})}.
\end{aligned}
\end{equation}

Now let $\Psi\in L^p(\O;C([0,T];E_{-\theta_F}))$. By applying \eqref{eq:detphi2}  and \eqref{eq:detphi1} to the paths $\Psi(\cdot, \omega)$
one obtains that $S*\Psi\in \V$ and
\begin{equation}\label{eq:detLphi}
\|S*\Psi\|_{\V} \le C T^{\min\{\frac12-\a, 1-\eta-\theta_F\}} \|
\Psi\|_{L^p(\O;C([0,T];E_{-\theta_F}))}.
\end{equation}

(b): \ Let $\phi_1, \phi_2\in \V$. Since $F$ is of linear growth, $
F(\cdot,\phi_1)$ and $F(\cdot,\phi_2)$ belong to
$L^p(\O;C([0,T];E_{-\theta_F}))$. From \eqref{eq:detLphi} and the
fact that $F$ is Lipschitz continuous in its $E_\eta$-variable we
deduce that $S*(F(\cdot,\phi_1)), S*(F(\cdot,\phi_2))\in \V$ and
\begin{equation}\label{eq:detLFphi}
\begin{aligned}
\ & \|S*(F(\cdot,\phi_1)-F(\cdot,\phi_2))\|_{\V}
\\ & \quad \le C
T^{\min\{\frac12-\a, 1-\eta-\theta_F\}} \|
(F(\cdot,\phi_1)-F(\cdot,\phi_2))\|_{L^p(\O;C([0,T];E_{-\theta_F}))}
\\ & \quad
\le C T^{\min\{\frac12-\a, 1-\eta-\theta_F\}} L_{F}
\|\phi_1-\phi_2\|_{\V}.
\end{aligned}
\end{equation}

{\em Step 3: Estimating the stochastic convolution.}\ Again we proceed in two
steps.

(a): \ Let $\Psi:[0,T]\times\O\to \calL(H,E_{-\theta_B})$ be $H$-strongly
measurable and adapted and suppose that
\begin{equation}\label{eq:supEgH}
\sup_{t\in [0,T]}\E\|s\mapsto
(t-s)^{-\alpha}
\Psi(s)\|_{\g(L^2(0,t;H),E_{-\theta_B})}^p < \infty.
\end{equation}
We estimate
the $\V$-norm of $S\diamond \Psi$.

From Proposition \ref{prop:path-hcont2} we obtain an
$\e>0$  such that
\[\begin{aligned}
\Big(\E\|S\diamond \Psi\|_{C([0,T];E_\eta)}^p \Big)^{\frac1p} \le
CT^\e \sup_{t\in [0,T]}\Big(\E\|s\mapsto (t-s)^{-\alpha}
\Psi(s)\|_{\g(L^2(0,t;H),E_{-\theta_B})}^p\big)^\frac1p.
\end{aligned}\]

For the other part of the norm, by Proposition
\ref{prop:gammanormconv} we obtain that
\[
\begin{aligned}
\Big( & \E\|s\mapsto (t-s)^{-\alpha} S\diamond
\Psi(s))\|_{\g(L^2(0,t;H),E_\eta)}^p\Big)^{\frac{1}{p}}
\\ & \qquad\qquad \le C
T^{\frac12-\eta-\theta_B}
 \Big(\E \|s\mapsto
(t-s)^{-\alpha}
\Psi\|_{\g(L^2(0,t;H),E_{-\theta_B})}^p\Big)^{\frac{1}{p}}.
\end{aligned}
\]
Combining things we conclude that
\begin{equation}\label{eq:stochLPhi}
\begin{aligned}
\ & \|S\diamond \Psi\|_{\V}
\\ & \quad \le
CT^{\min\{\frac12-\eta-\theta_B, \e\}} \Big( \sup_{t\in [0,T]}
\Big(\E\|s\mapsto
(t-s)^{-\alpha}\Psi(s)\|^p_{\g(L^2(0,t;H),E_{-\theta_B})}\Big)^{\frac{1}{p}}.
\end{aligned}
\end{equation}

(b): \ For $t\in [0,T]$
let $\mu_{t,\alpha}$ be the finite measure on $((0,t), \mathcal{B}_{(0,t)})$ defined
by
\[\mu_{t,\alpha}(B) =  \int_0^t (t-s)^{-2\alpha}\one_{B}(s) \, ds.\]
Notice that for a function $\phi\in C([0,t];E)$ we have
\[\phi\in
\g(L^2((0,t),\mu_{t,\alpha}), E) \Longleftrightarrow s\mapsto
(t-s)^{-\alpha}\phi(s)\in \g(L^2(0,t),E).\] Trivially,
\[
\|\phi\|_{L^2((0,t),\mu_{t,\alpha};E)}
 =\|(t-\cdot)^{-\a}\phi(\cdot)\|_{L^2(0,t;E)}
 \le C t^{\frac12-\a}\|\phi\|_{C([0,T];E)}.
\]

Now let $\phi_1, \phi_2\in \V$. Since $B$ is $\gL$-Lipschitz and of
linear growth and $\phi_1$ and $\phi_2$ belong to
$L_\g^2((0,t),\mu_{t,\a};E_\eta)$ uniformly, $B(\cdot,\phi_1)$ and $B(\cdot,\phi_2)$ satisfy
\eqref{eq:supEgH}. Since $B(\cdot,\phi_1)$ and $B(\cdot,\phi_2)$ are
$H$-strongly measurable and adapted, it follows from
\eqref{eq:stochLPhi} that $B(\cdot,\phi_1)$, $B(\cdot,\phi_2)\in
\V$ and
\begin{equation}\label{eq:stochLBphi}
\begin{aligned}
\ & \|S\diamond (B(\cdot,\phi_1) - B(\cdot,\phi_2))\|_{\V}
\\ & \ \  \lesssim T^{\min\{\frac12-\eta-\theta_B, \e\}}
\\ & \ \  \ \  \times  \Big(
\sup_{t\in [0,T]} \Big(\E\|s\mapsto (t-s)^{-\alpha} [B(s,\phi_1(s))
-B(s,\phi_2(s))]\|^p_{\g(L^2(0,t;H),E_{-\theta_B})}\Big)^{\frac{1}{p}}
\\ & \ \  = T^{\min\{\frac12-\eta-\theta_B, \e\}}
\sup_{t\in [0,T]} \Big(\E\|B(\cdot,\phi_1)
-B(\cdot,\phi_2)\|^p_{\g(L^2((0,t), \mu_{t,\a};H),E_{-\theta_B})}\Big)^{\frac{1}{p}}
\\ & \ \  \lesssim L_B^{\g}T^{\min\{\frac12-\eta-\theta_B, \e\}}
\sup_{t\in [0,T]} \Big(\E\|\phi_1-\phi_2\|^p
_{L_\g^2((0,t), \mu_{t,\a};E_{\eta})}\Big)^{\frac{1}{p}}
\\ & \ \  \lesssim L_B^{\g}T^{\min\{\frac12-\eta-\theta_B, \e\}}
\Big[\sup_{t\in [0,T]} \Big(\E\|s\mapsto (t-s)^{-\alpha}[\phi_1
-\phi_2]\|^p_{\g(L^2(0,t),E_{\eta})}\Big)^{\frac{1}{p}}
\\ & \hskip5cm + T^{\frac{p}{2}-\a p}\E\n \phi_1
-\phi_2\|^p_{C([0,T]; E_{\eta})}\Big)^{\frac{1}{p}}\Big]
\\ & \ \  \lesssim
L_B^{\g}T^{\min\{\frac12-\eta-\theta_B, \e\}} \|\phi_1-\phi_2\|_{\V}.
\end{aligned}
\end{equation}

{\em Step 4: Collecting the estimates.}\ It follows from the above
considerations that $L_T$ is well-defined on $\V$ and there exist
constants $C\ge 0$ and $\beta>0$ such that for all $\phi_1,\phi_2\in
\V$ we have
\begin{equation}\label{eq:fixpointestprecize}
\|L_T(\phi_1)-L_T(\phi_2)\|_{\V}\le C T^{\beta}
\|\phi_1-\phi_2\|_{\V}.
\end{equation}

The estimate \eqref{eq:operatorest} follows from \eqref{eq:fixpointestprecize}
and
\[\|L_T(0)\|_{\V}\le C(1+(\E\|u_0\|_{E_\eta}^p)^{\frac1p}.\]
\end{proof}

\begin{theorem}[Existence and uniqueness]\label{thm:mainexistenceL}
Let $E$ be a UMD space with type $\tau\in [1, 2]$. Suppose that
{\rm \ref{as:semigroup}-\ref{as:initial_value}} are
satisfied and assume that $0\le
\eta+\theta_F<\frac32 -\frac1\tau$ and $0\le \eta+\theta_B<\frac12$.
Let $p>2$ and $\alpha\in (0,\frac12)$ be such that
$\eta+\theta_B<\a-\frac1p$. If $u_0\in L^p(\O,\F_0;E_\eta)$, then
there exists a mild solution $U$ in $\VVinf$ of \eqref{SE}. As a mild solution
in $V_{\a,p}^p([0,T]\times\O;E_\eta)$, this solution $U$ is unique.
Moreover, there exists a constant $C\ge 0$, independent of $u_0$, such that
\begin{equation}\label{eq:mainthmestimateL}
\|U\|_{\VVinf} \le C ( 1 + (\E\|u_0\|^p_{E_{\eta}})^{\frac{1}{p}}).
\end{equation}
\end{theorem}

\begin{proof}
By Proposition \ref{prop:contr} we can find $T\in (0,\Tend]$,
independent of $u_0$, such that $C_T<\frac12$.
It follows from \eqref{eq:fixpointest} and the Banach fixed point theorem that $L_T$ has a unique fixed point
$U\in \Vinf$. This gives a continuous adapted
process $U:[0,T]\times\O\to E_\eta$ such that almost surely for all $t\in [0,T]$,
\begin{equation}\label{eq:fixedpointsol}
U(t) = S(t) u_0 + S*F(\cdot,U)(t)+S\diamond B(\cdot,U)(t).
\end{equation}
Noting that $U = \limn L_T^n(0)$ in $\V$, \eqref{eq:operatorest}
implies the inequality
\[\n U\n_{\V} \le C(1+(\E\n
u_0\n_{E_\eta}^p)^\frac1p) + C_T \n U\n_{\V},\]
and then $C_T<
\frac12$ implies
\begin{equation}\label{eq:mainthmestimateLT}
\|U\|_{\V} \le C ( 1 + (\E\|u_0\|^p_{E_{\eta}})^{\frac{1}{p}}).
\end{equation}

Via a standard induction argument one may construct a mild solution
on each of the intervals $[T,2T], \ldots, [(n-1)T,nT], [nT,\Tend]$
for an appropriate integer $n$. The induced solution $U$ on
$[0,\Tend]$ is the mild solution of \eqref{SE}. Moreover, by
\eqref{eq:mainthmestimateLT} and induction we deduce
\eqref{eq:mainthmestimateL}.

For small $T\in (0,\Tend]$, uniqueness on $[0,T]$ follows from the uniqueness
of the fixed point of $L_T$ in $V_{\a,p}^p([0,T]\times\O;E_\eta)$. Uniqueness on $[0,T_0]$ follows by induction.
\end{proof}

In the next theorem we deduce regularity properties of the solution.
They are formulated for $U-S u_0$; if $u_0$ is regular enough,
regularity of $U$ can be deduced.

\begin{theorem}[Regularity]\label{thm:Holdercont}
Let $E$ be a UMD space with type $\tau\in [1, 2]$ and suppose that {\rm \ref{as:semigroup}-\ref{as:initial_value}} are
satisfied. Assume that $0\le
\eta+\theta_F<\frac32 -\frac1\tau$ and $0\le \eta+\theta_B<\frac12-\frac1p$ with $p>2$.
 Let $\lambda\ge 0$ and $\delta\geq \eta$
satisfy
$\lambda+\delta<\min\{\frac12-\frac1p-\theta_B, 1- \theta_F\}$.
Then there exists a constant $C\ge 0$ such that for all $u_0\in
L^p(\O;E_{\eta})$,
\begin{equation}\label{eq:mainthmestimate}
\Big(\E\|U-Su_0\|_{C^\lambda([0,\Tend];E_{\delta})}^p\Big)^{\frac{1}{p}}
\le C ( 1 + (\E\|u_0\|^p_{E_{\eta}})^{\frac{1}{p}}.
\end{equation}
\end{theorem}
\begin{proof}
Choose $r\ge 1$ and  $0<\alpha<\frac{1}{2}$
such that $\lambda+\delta< 1- \frac1{r} - \theta_F$,
 $\eta+\theta_B< \a
-\frac{1}{p}$, and $\lambda+\delta+\theta_B<\a-\frac1p$. Let
$\tilde{U}\in \VVinf$ be the mild solution from Theorem
\ref{thm:mainexistenceL}. It follows from Lemma \ref{lem:detconv} (with $\a=1$)
that we may take a version of $S*F(\cdot,\tilde{U})$ with
\[
\begin{aligned}\E\|S*F(\cdot,\tilde U)\|_{C^{\lambda}([0,\Tend];E_{\delta})}^p
& \le C\E\| F(\cdot,\tilde U)\|_{L^{r}(0,\Tend;E_{-\theta_F})}^p
\\ & \le C\E\| F(\cdot,\tilde U)\|_{C([0,\Tend];E_{-\theta_F})}^p.
\end{aligned}\] Similarly, via Proposition \ref{prop:path-hcont2} we may
take a version of $S*B(\cdot,\tilde{U})$ with
\[\begin{aligned} \E\|&S\diamond
B(\cdot,\tilde U)\|_{C^{\lambda}([0,\Tend]; E_{\delta})}^p
\\ & \le C
 \sup_{t\in [0,\Tend]}\E\|s\mapsto
(t-s)^{-\alpha}
B(\cdot,\tilde U(s))\|_{\g(L^2(0,t;H),E_{-\theta_B})}^p.
\end{aligned}
\]
Define $U:[0,\Tend]\times\O\to E_\eta$ as
\[U(t) = S(t) u_0 + S*F(\cdot,\tilde{U})(t)+S\diamond B(\cdot,\tilde{U})(t),\]
where we take the versions of the convolutions as above. By uniqueness we have
almost surely $U \equiv \tilde{U}$. Arguing as in \eqref{eq:stochLBphi}
deduce that
\[\E\|U-Su_0\|_{C^{\lambda}([0,\Tend];
E_{\delta})}^p\le C(1+\|U\|_{\VVinf}^p).
\]
Now \eqref{eq:mainthmestimate} follows from
\eqref{eq:mainthmestimateL}.
\end{proof}

\section{Stochastic evolution equations II: measurable initial values}
\label{sec:meas_initial}

So far we have solved the problem \eqref{SE} for initial values
$u_0\in L^p(\O,\F_0;E_\eta)$. In this section we discuss the case of
initial values $u_0\in L^0(\O,\F_0;E_\eta)$.

Fix $T\in (0,\Tend]$. For $p\in [1, \infty)$ and  $\alpha\in
(0,\frac12)$ we define $V^{0}_{\a,p}([0,T]\times\O;E)$ as the linear
space of continuous adapted processes $\phi:[0,T]\times \O\to E$
such that almost surely,
\[\|\phi\|_{C([0,T];E)}+ \Big(\int_0^T \|s\mapsto (t-s)^{-\alpha} \phi(s)\|_{\g(L^2(0,t),E)}^p \,dt\Big)^{\frac1p} <\infty.\]
As usual we identify indistinguishable processes.

\begin{theorem}[Existence and uniqueness]\label{thm:mainexistenceLloc}
Let $E$ be a UMD space of type $\tau\in [1, 2]$ and suppose
that {\rm \ref{as:semigroup}}-{\rm \ref{as:initial_value}} are
satisfied. Assume that $0\le\eta+\theta_F<\frac32 -\frac1\tau$ and
$\eta+\theta_B<\frac12$.
If $\alpha\in (0,\frac12)$ and $p>2$ are such that
$\eta+\theta_B<\a-\frac1p$, then
 there exists a unique mild solution
$U\in \VVo$ of \eqref{SE}.
\end{theorem}

For the proof we need the following uniqueness result.

\begin{lemma}\label{lem:localityexistence}
Under the conditions of Theorem \ref{thm:mainexistenceL} let $U_1$ and $U_2$ in
$\Vinf$ be the mild solutions of \eqref{SE} with initial values
$u_1$ and $u_2$ in $L^p(\O,\F_0;E_\eta)$. Then almost surely on the set $\{u_1 =
u_2\}$ we have $U_1 \equiv U_2$.
\end{lemma}
\begin{proof}
Let $\Gamma = \{u_1 = u_2\}$. First consider small $T\in (0,\Tend]$ as in Step
1 in the proof of Theorem \ref{thm:mainexistenceL}. Since $\Gamma$ is
$\F_0$-measurable we have
\[\begin{aligned}\|U_1\one_{\Gamma}-U_2\one_{\Gamma}\|_{\Vinf}
&=  \|L_T(U_1)\one_{\Gamma}-L_T(U_2)\one_{\Gamma}\|_{\Vinf}
\\ & = \|(L_T(U_1\one_{\Gamma})-L_T(U_2\one_{\Gamma}))\one_{\Gamma}\|_{\Vinf}
\\ & \le \tfrac12
\|U_1\one_{\Gamma}-U_2\one_{\Gamma}\|_{\Vinf},
\end{aligned}\]
hence almost surely $U_1|_{[0,T]\times\Gamma}\equiv U_2|_{[0,T]\times\Gamma}$.

To obtain uniqueness on the interval $[0,\Tend]$
one may proceed as in the proof of Theorem
\ref{thm:mainexistenceL}.
\end{proof}

\begin{proof}[Proof of Theorem \ref{thm:mainexistenceLloc}]
(Existence): \ Define $(u_n)_{n\geq 1}$ in $L^p(\O,\F_0;E_\eta)$ as
\[u_n := \one_{\{\|u_0\|_{E_\eta}\le n\}}u_0.\]
By Theorem \ref{thm:mainexistenceL}, for each $n\geq 1$ there is a
unique solution $U_n\in \V$ of \eqref{SE} with initial value $u_n$.
By Lemma \ref{lem:localityexistence} we may define
$U:(0,\Tend)\times\O\to E_\eta$ as $U(t) = \limn U_n(t)$ if this
limit exists and $0$ otherwise. Then, $U$ is strongly measurable and
adapted, and almost surely on $\{\|u_0\|_{E_\eta}\le n\}$, for all
$t\in (0,\Tend)$ we have $U(t) = U_n(t)$. Hence, $U\in \Vo$. It is
routine to check
that $U$ is a solution of \eqref{SE}.

(Uniqueness): \  The argument is more or less standard, but there are
some subleties due to the presence of the radonifying norms.

Let $U,V\in \VVo$ be mild solutions of \eqref{SE}.
For each $n\geq 1$ let the stopping times $\mu_n^U$ and $\nu_n^U$ be
defined as
\[\begin{aligned}
\mu_n^U &= \inf\Big\{r\in[0,T_0]:\int_0^{T_0} \|s\mapsto
(t-s)^{-\alpha}U(s) \one_{[0,r]}(s)\|_{\g(L^2(0,t),E_\eta)}^p
\,dt\geq n\Big\},
\\ \nu_n^U &= \inf\Big\{r\in[0,T]: \|U(r)\n_{E_\eta}
\geq n\Big\}.
\end{aligned}\] This is well-defined since
\[r\mapsto \int_0^{T_0} \|s\mapsto (t-s)^{-\alpha}U(s)
\one_{[0,r]}(s)\|_{\g(L^2(0,t),E_\eta)}^p \,dt\] is a continuous
adapted process by \cite[Proposition 2.4]{NVW} and the dominated
convergence theorem. The stopping times $\mu_n^V$ and $\nu_n^V$ are
defined in a similarly. For each $n\geq 1$ let
\[\tau_n = \mu_n^U \wedge \nu_n^U \wedge \mu_n^V \wedge \nu_n^V,\]
and let $U_n = U\one_{[0,\tau_n]}$ and $V_n = V\one_{[0,\tau_n]}$.
Then for all $n\geq 1$, $U_n$ and $V_n$ are in $\VVp$. One easily
checks that
\[
U_n = \one_{[0,\tau_n]} (L_T (U_n))^{\tau_n} \ \ \text{and} \ \ V_n =
\one_{[0,\tau_n]} (L_T (V_n))^{\tau_n},
\]
where $L_T$ is the map introduced preceding Proposition
\ref{prop:contr} and  $(L_T (U_n))^{\tau_n}(t) := (L_T
(U_n))(t\wedge \tau_n)$. By Proposition \ref{prop:contr} we can find
$T\in (0,\Tend]$ such that $C_T\le \frac12$. A routine computation
then implies
\[
\|U_n - V_n\|_{\Vp} \leq \tfrac12 \|U_n - V_n\|_{\Vp}.
\]
We obtain that $U_n = V_n$ in $\Vp$, hence $\P$-almost surely,
$U_n\equiv V_n$. Letting $n$ tend to infinity, we may conclude that
almost surely, $U\equiv V$ on $[0,T]$. This gives the uniqueness on
the interval $[0,T]$. Uniqueness on $[0,\Tend]$ can obtained by the
usual induction argument.
\end{proof}

Note that in the last paragraph of the proof we needed
to work in the space $\Vp$ rather than
in $\Vpwithouttilde$ because the truncation with the stopping time
destroys the pathwise continuity.

By applying Theorem \ref{thm:Holdercont} to the unique solution $U_n$ with initial value
$u_n := \one_{\{\|u_0\|_{E_\eta}\le n\}} u_0$,
the solution $U:= \limn U_n$ constructed in Theorem \ref{thm:mainexistenceLloc}
enjoys the following regularity property.

\begin{theorem}[H\"older regularity]\label{thm:Holdercontloc}
Let $E$ be a UMD space and type $\tau\in [1, 2]$ and suppose that
{\rm \ref{as:semigroup}}-{\rm \ref{as:initial_value}} are
satisfied. Assume that $0\le
\eta+\theta_F<\frac32 -\frac1\tau$ and $0\le \eta+\theta_B<\frac12$.
Let $\lambda\ge 0$ and $\delta\geq \eta$ satisfy
$\lambda+\delta<\min\{\frac12-\theta_B, 1-\theta_F\}$.
Then the mild solution $U$ of \eqref{SE} has a version such that
almost all paths satisfy $U-Su_0\in
C^\lambda([0,\Tend];E_{\delta})$.
\end{theorem}

\begin{proof}[Proof of Theorem \ref{thm:main}]
Part (1) is a the special case of Theorem \ref{thm:mainexistenceL} corresponding to
$\tau=1$ and $\theta_F = \theta_B = 0$. For part (2) we apply
Theorem \ref{thm:Holdercontloc}, again with $\tau=1$ and $\theta_F = \theta_B = 0$.
\end{proof}

\section{Stochastic evolution equations  III: the locally Lipschitz case}
\label{sec:loc_Lipschitz}

Consider the following assumptions on $F$ and $B$.
\let\ALTERWERTA\theenumi
\let\ALTERWERTB\labelenumi
\def\theenumi{(A2)$'$}
\def\labelenumi{(A2)$'$}
\begin{enumerate}
\item \label{as:locLipschitzF}
The function $F:[0,T_0]\times\O\times E_{\eta}\to E_{-\theta_F}$ is
locally Lipschitz, uniformly in $[0,T_0]\times\O$, i.e., for all
$R>0$ there exists a constant $L_F^R$ such that for all $t\in
[0,T_0], \ \omega\in \O \ \text{and} \
\|x\|_{E_\eta},\|y\|_{E_\eta}\leq R$,
\[\begin{aligned}
\|F(t,\omega,x) - F(t,\omega,y)\|_{E_{-\theta_F}} & \leq
L_F^R\|x-y\|_{E_\eta}.
\end{aligned}\]
Moreover, for all $x\in E_\eta$, $(t,\omega)\mapsto F(t,\omega,
x)\in E_{-\theta_F}$ is strongly measurable and adapted, and there exists
a constant $C_{F,0}$ such that for all $t\in [0,T_0] \ \text{and} \
\omega\in \O$,
\[\|F(t,\omega,0)\|_{E_{-\theta_F}}\leq C_{F,0}.\]
\end{enumerate}
\let\theenumi\ALTERWERTA
\let\labelenumi\ALTERWERTB

\let\ALTERWERTA\theenumi
\let\ALTERWERTB\labelenumi
\def\theenumi{(A3)$'$}
\def\labelenumi{(A3)$'$}
\begin{enumerate}
\item \label{as:locLipschitzB}
The function $B:[0,T_0]\times\O \times E_\eta\to
\calL(H,E_{-\theta_B})$ is locally $\gL$-Lipschitz, uniformly in
$\O$, i.e., there exists a sequence of $\gL$-Lipschitz functions
 $B_n:[0,T_0]\times\O\times E_\eta\to
\calL(H,E_{-\theta_B})$ such that $B(\cdot, x) = B_n(\cdot, x)$ for
all $\|x\|_{E_\eta}<n$. Moreover, for all $x\in E_\eta$,
$(t,\omega)\mapsto B(t,\omega, x)\in E_{-\theta_B}$ is $H$-strongly
measurable and adapted, and there exists a constant $C_{B,0}$ such
that for all finite measures $\mu$ on $([0,T_0],
\mathcal{B}_{[0,T_0]})$ and all $\omega\in\Omega$,
\[\|t\mapsto B(t,\omega,0)\|_{\g(L^2((0,T_0),\mu;H),E_{-\theta_B})}\leq C_{B,0}.\]

\end{enumerate}
\let\theenumi\ALTERWERTA
\let\labelenumi\ALTERWERTB
One may check that the locally Lipschitz version of Lemma
\ref{lem:type2Lipschitz} holds as well. This gives an easy way to
check \ref{as:locLipschitzB} for type $2$ spaces $E$.

Let $\varrho$ be a stopping time with values in $[0,T_0]$. For $t\in
[0,T_0]$ let
\[\O_t(\varrho) = \{\omega\in \O: t<\varrho(\omega)\},\]
\[[0,\varrho)\times\O=\{(t,\omega)\in [0,T_0]\times\O: 0\leq
t<\varrho(\omega)\},\]
\[[0,\varrho]\times\O=\{(t,\omega)\in [0,T_0]\times\O: 0\leq t\leq \varrho(\omega)\}.\]
A process $\zeta:[0,\varrho)\times \O \to E$ (or $(\zeta(t))_{t\in
[0,\varrho)})$ is called {\em admissible} if for all $t\in [0,T_0]$,
$\O_t(\varrho)\ni\omega\to \zeta(t,\omega)$ is $\F_t$-measurable and
for almost all $\omega\in \O$, $[0,\varrho(\omega))\ni t\mapsto
\zeta(t,\omega)$ is continuous.

Let $E$ be a UMD space. An admissible $E_\eta$-valued process
$(U(t))_{t\in [0,\varrho)}$ is called a {\em local solution} of
\eqref{SE} if $\varrho\in(0,\Tend]$ almost surely and there exists
an increasing sequence of stopping times $(\varrho_n)_{n\geq 1}$
with $\varrho = \limn \varrho_n$ such that
\begin{enumerate}
\item[(i)] \!for all $t\in [0,T_0]$, $s\mapsto S(t-s) F(\cdot,U(s))
\one_{[0,\varrho_n]}(s) \in L^0(\O;L^1(0,t;E_\eta))$,
\item[(ii)] \!for all $t\in [0,T_0]$,
$s\mapsto S(t-s) B(\cdot,U(s)) \one_{[0,\varrho_n]}(s)\!\in\! L^0(\O;\g(L^2(0,t;H),E_\eta))$,
\item[(iii)] \!almost surely for all $t\in [0,\rho_n]$,
\[U(t) = S(t) u_0 + S*F(\cdot,U)(t) + S\diamond B(\cdot,U)(t).\]
\end{enumerate}
By (i) the deterministic convolution is defined pathwise as a
Bochner integral. Since $E$ is a UMD space, by (ii) and Proposition
\ref{prop:NVW} we may define the stochastic convolution as
\[S\diamond B(\cdot,U)(t) = \int_0^{t} S({t-s}) B(s,U(s))
\one_{[0,\varrho_n]}(s) \, d W_H(s), \ \ t\in [0,\rho_n].\]

A local solution $(U(t))_{t\in [0,\varrho)}$ is called {\em maximal}
for a certain space $V$ of $E_\eta$-valued admissible processes if
for any other local solution
$(\tilde{U}(t))_{t\in[0,\tilde{\varrho})}$ in $V$, almost surely we
have $\tilde \varrho\leq \varrho$ and $\tilde{U}\equiv
U|_{[0,\tilde \varrho)}$. Clearly, a maximal local solution for such
a space $V$ is always unique in $V$. We say that a local solution
$(U(t))_{t\in [0,\varrho)}$ of \eqref{SE} is a {\em global solution}
of \eqref{SE} if $\varrho = \Tend$ almost surely and $U$ has an
extension to a solution $\hat{U}:[0,\Tend]\times\O\to E_\eta$ of
\eqref{SE}. In particular, almost surely ``no blow" up occurs at
$t=\Tend$.

We say that $\varrho$ is an {\em explosion time} if for almost all
$\omega\in \O$ with $\varrho(\omega)<\Tend$,
\[
\limsup_{t\uparrow \varrho(\omega)} \|U(t,\omega)\|_{E_\eta} =
\infty.
\]
Notice that if $\varrho=\Tend$ almost surely, then $\varrho$ is
always an explosion time in this definition. However, there need not
be any ``blow up" in this case.

Let $\varrho$ be a stopping time with values in $[0,T_0]$. For $p\in
[1, \infty)$, $\alpha\in [0,\frac12)$ and $\eta\in [0,1]$ we define
$V^{0,{\rm loc}}_{\alpha,p}([0,\varrho)\times \O;E)$ as all
$E$-valued admissible processes $(\phi(t))_{t\in [0,\varrho)}$ such
that there exists an increasing sequence of stopping times
$(\varrho_n)_{n\geq 1}$ with $\varrho = \limn \varrho_n$ and almost
surely
\[\|\phi \|_{C([0,\varrho_n];E)} + \Big(\int_0^{T} \|s\mapsto
(t-s)^{-\alpha}\phi(s) \one_{[0,\varrho_n]}(s)\|_{\g(0,t;E)}^p \,dt\Big)^{\frac1p}<\infty.\]

In the case that for almost all $\omega$, $\varrho_n(\omega)=T$ for
$n$ large enough,
\[V^{0,{\rm loc}}_{\alpha,p}([0,\varrho)\times \O;E) = V^{0}_{\a,{p}}([0,T]\times\O;E).\]

\begin{theorem}\label{thm:mainexistenceLocallyLipschitzCase}
Let $E$ be a UMD space with type $\tau\in [1, 2]$ and suppose that
{\rm\ref{as:semigroup}}, {\rm\ref{as:locLipschitzF}},
{\rm\ref{as:locLipschitzB}}, {\rm \ref{as:initial_value}} are
satisfied, and assume that $0\leq \eta+\theta_F<\frac32 -
\frac1\tau$.
\begin{enumerate}
\item For all $\alpha\in (0,\frac12)$
and $p>2$ such that $\eta+\theta_B<\a-\frac1p$ there exists a unique
maximal local solution $(U(t))_{[0,\varrho)}$ in $\Vadm$ of
\eqref{SE}.
\item For all $\lambda>0$ and $\delta\geq \eta$ such that
$\lambda+\delta<\min\{\frac{1}{2}-\theta_B,
1-\theta_F\}$,
 $U$ has a version such that for almost all $\omega\in \O$,
\[t\mapsto U(t,\omega) - S(t) u_0(\omega)\in C^\lambda_{\rm loc}([0,\varrho(\omega));E_{\delta}),\]
\end{enumerate}
If in addition the linear growth conditions of
{\rm\ref{as:LipschitzF}} and  {\rm\ref{as:LipschitzB}} hold, then
the above function $U$ is the unique global solution of \eqref{SE}
in $\VVo$ and the following assertions hold:
\begin{enumerate}
\item[(3)]\label{it:lingrowth1} The solution $U$ satisfies the statements of
Theorems \ref{thm:mainexistenceLloc} and \ref{thm:Holdercontloc}.
\item[(4)]\label{it:lingrowth2} If $\alpha\in (0,\frac12)$ and $p>2$ are such that $\alpha>\eta+\theta_B+\frac1p$
and $u_0\in L^p(\O,\F_0;E_\eta)$, then the solution $U$ is in
$\VVinf$ and \eqref{eq:mainthmestimateL} and the statements of
Theorem \ref{thm:Holdercont} hold.
\end{enumerate}
\end{theorem}

Before we proceed, we prove the following local uniqueness result.
\begin{lemma}\label{lem:localuniqueness2}
Suppose that the conditions of Theorem
\ref{thm:mainexistenceLocallyLipschitzCase} are satisfied and let
 $(U_1(t))_{t\in [0,\varrho_1)}$ in $\Vadmone$ and
$(U_2(t))_{t\in [0,\varrho_2)}$ in $\Vadmtwo$ be local solutions of
\eqref{SE} with initial values $u_0^1$ and $u_0^2$. Let $\Gamma =
\{u_0^1=u_0^2\}$. Then almost surely on $\Gamma$,
$U_1|_{[0,\varrho_1\wedge \varrho_2)} \equiv
U_2|_{[0,\varrho_1\wedge \varrho_2)}$. Moreover, if $\varrho_1$ is
an explosion time for $U_1$, then almost surely on $\Gamma$,
$\varrho_1 \geq \varrho_2$. If $\varrho_1$ and $\varrho_2$ are
explosion times for $U_1$ and $U_2$, then almost surely on $\Gamma$,
$\varrho_1=\varrho_2$ and $U_1\equiv U_2$.
\end{lemma}
\begin{proof}
Let $\varrho = \varrho_1\wedge \varrho_2$. Let $(\mu_n)_{n\geq 1}$
be an increasing sequences of bounded stopping times such that
$\limn \mu_n = \varrho$ and for all $n\geq 1$, $U_1
\one_{[0,\mu_n]}$ and $U_2 \one_{[0,\mu_n]}$ are in $\VVp$. Let
\[\nu^1_n = \inf\{t\in [0,T_0]: \|U_1(t)\|_{E_\eta}\geq n\} \ \text{and} \ \nu^2_n = \inf\{t\in [0,T_0]: \|U_2(t)\|_{E_\eta}\geq n\}\]
and let $\sigma^i_n = \mu_n\wedge \nu^i_n$ and let $\sigma_n =
\sigma^1_n \wedge\sigma^2_n$. On $[0,T_0]\times\O\times \{x\in
E_\eta: \|x\|_{E_\eta}\leq n\}$ we may replace $F$ and $B$ by $F_n$
(for a possible definition of $F_n$, see the proof of Theorem
\ref{thm:mainexistenceLocallyLipschitzCase}) and $B_n$ which satisfy
\ref{as:LipschitzF} and \ref{as:LipschitzB}. As in the proof of
Theorem \ref{thm:mainexistenceLloc} it follows that for all $0<
T\leq T_0$,
\[\begin{aligned}
\|U_1^{\sigma_n}&\one_{[0,\sigma_n]\times\Gamma} -
U_2^{\sigma_n}\one_{[0,\sigma_n]\times\Gamma}\|_{\Vp} \\ &=
\|(L_T(U_1^{\sigma_n}\one_{[0,\sigma_n]\times\Gamma}) -
L_T(U_2^{\sigma_n}\one_{[0,\sigma_n]\times\Gamma}))\one_{[0,\sigma_n]\times\Gamma}\|_{\Vp}
\\ &\leq \|L_T(U_1^{\sigma_n}\one_{[0,\sigma_n]\times\Gamma}) -
L_T(U_2^{\sigma_n}\one_{[0,\sigma_n]\times\Gamma})\|_{\Vp}
\\ &\leq C_T\|U_1^{\sigma_n}\one_{[0,\sigma_n]\times\Gamma} -
U_2^{\sigma_n}\one_{[0,\sigma_n]\times\Gamma}\|_{\Vp},
\end{aligned}
\]
where $C_T$ satisfies $\lim_{T\downarrow 0} C_T = 0$. Here
$\one_{[0,\sigma_n]\times\Gamma}$ should be interpreted as the
process $(t,\omega)\mapsto \one_{[0,\sigma_n(\omega)]\times\Gamma}(t,\omega)$.
For $T$ small enough it follows that
$U_1^{\sigma_n}\one_{[0,\sigma_n]\times\Gamma} =
U_2^{\sigma_n}\one_{[0,\sigma_n]\times\Gamma}$ in $\Vp$. By an
induction argument this holds on $[0,T_0]$ as well. By path
continuity it follows that almost surely, $U_1 \equiv U_2$ on
$[0,\sigma_n]\times\Gamma$. Since $\varrho = \limn \sigma_n$ we may
conclude that almost surely, $U_1 \equiv U_2$ on
$[0,\varrho)\times\Gamma$.

If $\varrho_1$ is an explosion time, then as in \cite[Lemma
5.3]{Sei} this yields $\varrho_1 \geq \varrho_2$ on $\Gamma$ almost
surely. Indeed, if for some $\omega\in \Gamma$,
$\varrho_1(\omega)<\varrho_2(\omega)$, then we can find an $n$ such
that $\varrho_1(\omega)<\nu^2_n(\omega)$. We have $U_1(t,\omega) =
U_2(t,\omega)$ for all $0\leq t\leq \nu^1_{n+1}(\omega)<
\varrho_1(\omega)$. If we combine both assertions we obtain that
\[n+1=\|U_1(\nu^1_{n+1}(\omega),\omega)\|_{E_a} = \|U_2(\nu^1_{n+1}(\omega),\omega)\|_{E_a}\leq n.\]
This is a contradiction. The final assertion is now obvious.
\end{proof}

\begin{proof}[Proof of Theorem \ref{thm:mainexistenceLocallyLipschitzCase}]
We follow an argument of \cite{Brz2, Sei}.

For $n\geq 1$ let $\Gamma_n=\{\|u_0\|\leq \frac{n}{2}\}$ and $u_n =
u_0 \one_{\Gamma_n}$. Let $(B_n)_{n\geq 1}$ be the sequence of
$\gL$-Lipschitz functions from \ref{as:locLipschitzB}. Fix an
integer $n\geq 1$. Let $F_n:[0,T_0]\times \O\times E_\eta\to
E_{-\theta_F}$ be defined by
\[\begin{aligned}
F_n(\cdot,x) = F(\cdot,x) \ \ \text{for} \ \|x\|_{E_\eta}\leq n,
\end{aligned}\]
and $F_n(\cdot, x) = F\big(\cdot, \frac{n x}{\|x\|_{E_\eta}}\big)$ otherwise.
Clearly, $F_n$ and $B_n$ satisfy \ref{as:LipschitzF} and
\ref{as:LipschitzB}. It follows from Theorem
\ref{thm:mainexistenceL} that there exists a solution $U_n\in \VV$
of \eqref{SE} with $u_0$, $F$ and $B$ replaced by $u_n$, $F_n$ and
$B_n$. In particular, $U_n$ has a version with continuous paths. Let
$\varrho_n$ be the stopping time defined by
\[\varrho_n(\omega) = \inf\{t\in [0,T_0]: \|U_n(t,\omega)\|_{E_\eta}\geq n\}.\]
It follows from Lemma \ref{lem:localuniqueness2} that for all $1\leq
m\leq n$, almost surely, $U_m \equiv U_n$ on $[0,\varrho_m\wedge
\varrho_n]\times\Gamma_m$. By path continuity this implies
$\varrho_m\leq\varrho_n$. Therefore, we can define
$\varrho=\lim_{n\to \infty} \varrho_n$ and on $\Gamma_n$, $U(t) =
U_n(t)$ for $t\leq \varrho_n$. By approximation and Lemma
\ref{lem:localityloc} it is clear that $U\in \Vadm$ is a local
solution of \eqref{SE}. Moreover, $\varrho$ is an explosion time.
This proves the existence part of (1). Maximality is a consequence
of Lemma \ref{lem:localuniqueness2}. Therefore,
$(U(t))_{t\in[0,\varrho)}$ is a maximal local solution. This
concludes the proof of (1).

We continue with (2). By Corollary \ref{thm:Holdercont}, each
$U_n$ has the regularity as stated by (2). Therefore, the
construction yields the required pathwise regularity properties of $U$.

Turning to (4), let $(U_n)_{n\geq 1}$ be as before. As
in the proof of Proposition \ref{prop:contr} one can check that by
the linear growth assumption,
\[\begin{aligned}
\|U_n\|_{\V} &= \|L_T(U_n)\|_{\V} \\ & \leq C_T \|U_n\|_{\V} + C +
C\|u_n\|_{L^p(\O;E_\eta)},
\end{aligned}
\] where the constants do not
depend on $n$ and $u_0$ and we have $\lim_{T\downarrow 0}C_T =0$.
Since $\|u_n\|_{L^p(\O;E_\eta)}\leq \|u_0\|_{L^p(\O;E_\eta)}$, it
follows that for $T$ small we have
\[\|U_n\|_{\V}\leq C(1+\|u_0\|_{L^p(\O;E_\eta)}),\] where
$C$ is a constant independent of $n$ and $u_0$. Repeating this
inductively, we obtain a constant $C$ independent of $n$ and $u_0$
such that $\|U_n\|_{V_{\a,\infty}^p([0,T_0]\times\Omega;E_\eta)}\leq C(1+\|u_0\|_{L^p(\O;E_\eta)})$. In
particular,
\[\E\sup_{s\in [0,T_0]}\|U_n(s)\|^p_{E_\eta}\leq C^p (1+\|u_0\|_{L^p(\O;E_\eta)})^p.\]
It follows that
\[\P(\sup_{s\in [0,T_0]}\|U_n(s)\|_{E_\eta}\geq n)\leq  C^p n^{-p}.\]
Since $\sum_{n\geq 1}n^{-p}<\infty$, the Borel-Cantelli Lemma
implies that
\[\P\Big(\bigcap_{k\geq 1}\bigcup_{n\geq k}\Big\{\sup_{s\in [0,T_0]}\|U_n(s)\|_{E_\eta}\geq n\Big\}\Big) = 0.\]
This gives that almost surely, $\varrho_n=T_0$ for all $n$ large
enough, where $\varrho_n$ is as before. In particular, $\varrho=T_0$
and by Fatou's lemma
\[
\|U\|_{\VV}\leq \liminf_{n\to \infty} \|U_n\|_{\VV}\leq
C(1+\|u_0\|_{L^p(\O;E_\eta)}).
\]
Via an approximation argument one can check that $U$ is a global
solution. The final statement in (4) can be obtained as in Theorem
\ref{thm:Holdercont}.

For the proof of (3) one may repeat the construction of Theorem
\ref{thm:mainexistenceLloc}, using Lemma \ref{lem:localuniqueness2}
instead of Lemma \ref{lem:localityexistence}.
\end{proof}

\section{Generalizations to one-sided UMD spaces\label{sec:gen}}
In this section we explain how the theory of the preceding sections
 can be extended to a class of Banach spaces which contains, besides all UMD spaces,
the spaces $L^1$.

A Banach space $E$ is called a {\em UMD$^+$-space} if  for some
(equivalently, for all) $p\in (1,\infty)$ there exists a constant
$\beta_{p,E}^+\ge 1$ such that for all $E$-valued $L^p$-martingale
difference sequences $(d_j)_{j=1}^n$ we have
\[
\Bigl(\E\Bigl\|\sum_{j=1}^n r_j d_j \Bigr\|^p\Bigr)^\frac1p \le
\beta_{p,E}^+ \, \Bigl(\E\Bigl\|\sum_{j=1}^n
d_j\Bigr\|^p\Bigr)^\frac1p
\]
where $(r_j)_{j=1}^n$ is a Rademacher sequence independent of
$(d_j)_{j=1}^n$. The space $E$ is called a {\em UMD$^-$ space} if the
reverse inequality holds:
\[
\Bigl(\E\Bigl\|\sum_{j=1}^n d_j \Bigr\|^p\Bigr)^\frac1p \le
\beta_{p,E}^- \, \Bigl(\E\Bigl\|\sum_{j=1}^n r_j
d_j\Bigr\|^p\Bigr)^\frac1p
\]
Both classes of spaces were introduced and studied by Garling
\cite{Ga2}. By a standard randomization argument, every UMD spaces
is both UMD$^+$ and UMD$^-$, and conversely a Banach space which is
is both UMD$^+$ and UMD$^-$ is UMD. At present, no examples are
known of UMD$^+$-spaces which are not UMD. For the UMD$^-$property
the situation is different: if $E$ is UMD$^-$, then also $L^1(S;E)$
is UMD$^-$. In particular, every $L^1$-space is UMD$^-$ (cf.
\cite{NVfub}).

Assume that $(\mathcal{F}_t)_{t\geq 0}$ is the complete filtration
induced by $W_H$. If $E$ is a UMD$^-$-space, condition \eqref{NVW3}
still gives a {\em sufficient} condition for stochastic
integrability of $\Phi$, and instead of a norm equivalence one
obtains the one-sided estimate
\[ \E \Big\n \int_0^T \Phi\,dW_H\Big\n^p \lesssim_{p,E}
\E\n R\n_{\g(L^2(0,T;H),E)}^p\] for all $p\in (1, \infty)$, where we use the
notations of Proposition \ref{prop:NVW}.
The condition on the
filtration is needed for the approximation argument used in
\cite{Ga1}. By using Fubini's theorem it is obvious that the result also
holds if the probability space has the following product structure
$\O = \O_1\times \O_2$, $\F = \F\otimes \G$, $\P = \P_1\otimes
\P_2$, and the filtration is of the form $(\F_t\otimes \G)_{t\geq
0}$.

Mutatis mutandis, the theory presented in the previous sections
extends to UMD$^-$ spaces $E$, with two exceptions: (i) Proposition
\ref{prop:gammanormconv} relies, via the use of Lemma
\ref{lem:Stein}, on the fact that UMD spaces have property
$(\Delta)$; this property should now be included into the
assumptions. (ii) One needs the above assumption on the filtration.
We note that it follows from \cite{CV} that for $E=L^1$ the
assumption on the filtration is not needed.

\section{Applications to stochastic PDEs}
\label{sec:applications}

\paragraph{\bf Case of bounded $A$} We start with the case of a bounded
operator $A$. By putting $\tilde F :=  A + F$ it suffices to consider the case
$A=0$.

Let $E$ be a UMD$^-$ space with property $(\alpha)$ (see Section \ref{sec:g-Lipschitz}). Consider the
equation
\begin{equation}\label{eq:bounded}
\begin{aligned}
d U(t) & = F(t,U(t)) \, dt + B(U(t)) \, d W_E(t), \ \ t\in [0,T],
\\ U(0)& = u_0,
\end{aligned}
\end{equation}
where $W_E$ is an $E$-valued Brownian motion. With every $E$-valued
Brownian motion $W_E$ one can canonically associate an
$H$-cylindrical Brownian motion $W_H$, where $H$ is the so-called
reproducing kernel Hilbert space associated with $W_E(1)$ (see the
proof of Theorem \ref{thm:A0} below). Using this $H$-cylindrical
Brownian motion $W_H$, the problem \eqref{eq:bounded} can be
rewritten as a special instance of \eqref{SE}.

We make the following assumptions:
\begin{enumerate}
\item $F:[0,T]\times\O\times E\to E$ satisfies
\ref{as:LipschitzF} with $a=\theta_F=0$;
\item $B\in \calL(E,\calL(E))$;
\item $u_0:\O\to E$ is $\F_0$-measurable.
\end{enumerate}

\begin{theorem}\label{thm:A0}
Under these assumptions, for all $\alpha>0$ and $p>2$ such that $\alpha<\frac12-\frac1p$ there
exists a unique strong and mild solution $U:[0,T]\times\O\to E$ of
\eqref{eq:bounded} in $V^{0}_{\a,p}([0,T]\times\O;E)$. Moreover, for
all $0\le\lambda<\frac12$, $U$ has a version with paths in
$C^{\lambda}([0,T];E)$.
\end{theorem}
\begin{proof}
Let $H$ be the reproducing kernel Hilbert space associated with $W_E(1)$.
Then $H$ is a separable Hilbert space which is continuously embedded into $E$ by
means of an inclusion operator $i:H\embed E$ which belongs to $\g(H,E)$.
Putting $ W_H(t) i^*
x^* := \lb W_E(t), x^*\rb$ (cf. \cite[Example 3.2]{NW1}) we obtain an
$H$-cylindrical Brownian motion.

Assumption
\ref{as:semigroup} is trivially fulfilled, and (A2) and (A4) hold by assumption.
Let $\hat{B}\in \calL(E,\g(H,E))$ be given by
$\hat{B}(x) h = B(x) i h$.
Using Lemma \ref{lem:alphaGammaLipschitz}
one checks that $\hat{B}$ satisfies \ref{as:LipschitzB} with
$a=\theta_B=0$. Therefore, the result follows from Theorems
\ref{thm:mainexistenceLloc} and \ref{thm:Holdercontloc} (applied to
$\hat B$ and the $H$-cylindrical Brownian motion $W_H$).
Here we use the extension to UMD$^-$ space as explained in
Section \ref{sec:gen}.
\end{proof}

\paragraph{\bf Elliptic equations on bounded domains\label{sec:applGenEll}}

Below we will consider an elliptic equation of order $2m$ on a
domain $S\subseteq \R^d$. We will assume the noise is white in space
and time. The regularizing effect of the elliptic operator will be
used to be able to consider the white-noise in a suitable way.
Space-time white noise equations seem to be studied in the
literature in the case $m=1$ (cf. \cite{Brz2,DPZ}).

Let $S\subseteq \R^d$ be a bounded domain with $C^\infty$ boundary.
We consider the problem
\begin{equation}\label{eq:SEwhitenoise2}
\begin{aligned}
\frac{\partial u}{\partial t}(t,s) &= A(s,D) u(t,s)+ f(t,s,u(t,s))\\
& \qquad  + g(t,s,u(t,s)) \, \frac{\partial w}{\partial t}(t,s),
 && s\in S, \ t\in
(0,T],
\\  B_j(s,D) u(t,s) &= 0, && s\in \partial S, \ t\in (0,T],
\\ u(0,s) &= u_0(s), && s\in S.
\end{aligned}
\end{equation} Here $A$ is of the form
\[
A(s,D) = \sum_{|\alpha|\leq 2m} a_{\alpha}(s) D^{\alpha}
\]
where $D = -i(\partial_1, \ldots, \partial_d)$ and for $j=1, \ldots,
m$,
\[B_j(s,D) = \sum_{|\beta|\leq m_j} b_{j\beta}(s) D^{\beta}\]
where $1\leq m_j<2m$ is an integer. We assume that $a_{\alpha}\in
C(\overline{S})$ for all $|\alpha|=2m$. For $|\alpha|<2m$ the
coefficients $a_{\alpha}$ are in $L^\infty(S)$. For the principal
part $\sum_{|\alpha|= 2m} a_{\alpha}(s) D^{\alpha}$ of $A$ we assume
that there is a $\kappa>0$ such that
\[ (-1)^{m+1} \sum_{|\alpha|=2m} a_{\alpha}(s) \xi^{\alpha} \geq \kappa
|\xi|^{2m},  \ \ s\in S, \ \xi\in \R^d.\] For the coefficients of
the boundary value operator we assume that for $j=1, \ldots, m$ and
$|\beta|\leq m_j$ we have $b_{j\beta}\in C^{\infty}(\overline{S})$.
The boundary operators $(B_j)_{j=1}^m$ define a normal system of
Dirichlet type, i.e. $0\leq m_j<m$ (cf. \cite[Section 3.7]{Ta1}).
The $C^\infty$ assumption on the boundary of $S$ and on the
coefficients $b_{j\beta}$ is made for technical reasons. We will
need complex interpolation spaces for Sobolev spaces with boundary
conditions. It is well-known to experts that one can reduce the the
assumption to $S$ has a $C^{2m}$-boundary and $b_{j\beta}\in
C^{2m-m_j}(\overline{S})$. However, this seems not to be explicitly
contained in the literature.

The functions $f,g:[0,T]\times\O\times S\times \R\to \R$ are jointly
measurable, and adapted in the sense that for each $t\in [0,T]$,
$f(t,\cdot)$ and $g(t,\cdot)$ are $\F_t\otimes
\mathcal{B}_{S}\otimes \mathcal{B}_{\R}$-measurable. Finally, $w$ is
a space-time white noise (see, e.g., \cite{Walsh}) and $u_0:S\times
\O \to \R$ is an $\mathcal{B}_{S}\otimes\F_0$-measurable initial
value condition. We say that $u:[0,T]\times\O\times S\to \R$ is a
{\em solution} of \eqref{eq:SEwhitenoise2} if the corresponding
functional analytic model \eqref{SE} has a mild solution $U$ and
$u(t,s,\omega) = U(t,\omega)(s)$.

Consider the following conditions:

\let\ALTERWERTA\theenumi
\let\ALTERWERTB\labelenumi
\def\theenumi{(C1)}
\def\labelenumi{(C1)}
\begin{enumerate}
\item \label{cond:localLipfgb}
The functions $f$ and $g$ are locally Lipschitz in the fourth
variable, uniformly on $[0,T]\times\O\times S$, i.e., for all $R>0$
the exist constants $L^R_f$ and $L^R_g$ such that
\[\begin{aligned}
|f(t,\omega,s,x) - f(t,\omega,s,y)| & \leq L^R_f |x-y|, \\
|g(t,\omega,s,x) - g(t,\omega,s,y)| & \leq L^R_g |x-y|,
\end{aligned}
\] for all $t\in [0,T]$, $\omega\in \O$, $s\in S$, and  $|x|,|y|<R.$
Furthermore, $f$ and $g$ satisfy the boundedness conditions
\[\sup |f(t,\omega,s,0)|<\infty, \qquad  \sup |g(t,\omega,s,0)|<\infty,\]
where the suprema are taken over $t\in [0,T]$, $\omega\in\O$, and
$s\in S$.

\end{enumerate}
\let\theenumi\ALTERWERTA
\let\labelenumi\ALTERWERTB

\let\ALTERWERTA\theenumi
\let\ALTERWERTB\labelenumi
\def\theenumi{(C2)}
\def\labelenumi{(C2)}
\begin{enumerate}
\item \label{cond:lingrowthfgb}
The functions $f$ and $g$ are of linear growth in the fourth
variable, uniformly in $[0,T]\times\O\times S$, i.e., there exist
constants $C_f$ and $C_g$ such that
\[
|f(t,\omega,s,x)|  \leq C_f (1+|x|),\qquad |g(t,\omega,s,x)|  \leq
C_g (1+|x|), \] for all $t\in [0,T]$, $\omega\in \O$, $s\in S$, and
$x\in \R$.
\end{enumerate}
\let\theenumi\ALTERWERTA
\let\labelenumi\ALTERWERTB

Obviously, if $f$ and $g$ are Lipschitz and $f(\cdot, 0)$ and
$g(\cdot, 0)$ are bounded, i.e., if \ref{cond:localLipfgb} holds with
constants $L_f$ and $L_g$ not depending on $R$, then
\ref{cond:lingrowthfgb} is automatically fulfilled.

The main theorem of this section will be formulated in the terms of
the spaces $B^{s}_{p,1, \{B_j\}}(S)$. For their definition and further
properties
we refer to \cite[Section 4.3.3]{Tr} and references therein. For
$p\in [1, \infty]$, $q\in [1, \infty]$ and $s>0$, let
\[
\begin{aligned}
H^{s,p}_{\{B_j\}}(S) & := \big\{f\in H^{s,p}(S):\ B_j f = 0 \
\text{for $m_j<s-\frac1p$}, \ j=1, \ldots,
m\},\\
C^{s}_{\{B_j\}}(\overline{S}) & := \big\{f\in C^s(S):\ B_j f = 0 \
\text{for $m_j\leq s$}, \ j=1, \ldots, m\big\}.
\end{aligned}
\]
For $p\in (1, \infty)$ let $A_p$ be the realization of $A$ on the
space $L^p(S)$ with domain $H^{2m, p}_{\{B_j\}}(S)$. In this way
$-A_p$ is the generator of an analytic $C_0$-semigroup
$(S_p(t))_{t\geq 0}$. Since we may replace $A$ and $f$ in
\eqref{eq:SEwhitenoise2} by $A-w$ and $w+f$, we may assume that
$(S_p(t))_{t\geq 0}$ is uniformly exponentially stable. From
\cite[Theorem 4.1]{Se} and \cite[Theorem 1.15.3]{Tr} (also see
\cite{DDHPV}) we deduce that if $\theta\in (0,1)$ and $p\in (1,
\infty)$ are such that
\begin{equation}\label{eq:HBinterpolationcond}
\text{$2m \theta - \frac1p \neq m_j$, for all $j=1, \ldots, m$, }
\end{equation}
then
\[[L^p(S), D(A_p)]_{\theta} = [L^p(S), H^{2m,p}_{\{B_j\}}(S)]_{\theta} = H^{2m\theta,p}_{\{B_j\}}(S)\]
isomorphically.

\begin{theorem}\label{thm:mainexdm}
Assume that {\rm \ref{cond:localLipfgb}} holds, let  $\frac{d}{m}<2$, and let $p\in
(1, \infty)$ be such that $\frac{d}{2mp}<\frac12-\frac{d}{4m}$.
\begin{enumerate}
\item\label{it:appl1b}
If $\eta\in (\frac{d}{2mp},\frac12-\frac{d}{4m})$ is such that
\eqref{eq:HBinterpolationcond} holds for the pair $(\eta,p)$ and if
$u_0\in H^{2 m \eta,p}_{\{B_j\}}(S)$ almost surely, then for all
$r>2$ and $\alpha\in (\eta+\frac{d}{4m},\frac12-\frac1r)$ there
exists a unique maximal solution $(u(t))_{t\in [0,\varrho)}$ of
\eqref{eq:SEwhitenoise2} in $V^{0,{\rm
loc}}_{\alpha,r}([0,\varrho)\times \O;H^{2m\eta,p}_{\{B_j\}}(S))$.
\item\label{it:appl2b}
Moreover, if $\delta>\frac{d}{2mp}$ and $\lambda\ge 0$ are such that
$\delta+\lambda<\frac12-\frac{d}{4m}$ and
\eqref{eq:HBinterpolationcond} holds for the pair $(\delta,p)$, and
if $u_0\in H^{m-\frac{d}{2},p}_{\{B_j\}}(S)$ almost surely, then $u$
has paths in $C_{\rm loc}^{\lambda}([0,\tau);H^{2
m\delta,p}_{\{B_j\}}(S))$ almost surely.
\end{enumerate}
Furthermore, if condition {\rm\ref{cond:lingrowthfgb}} holds as well,
then:
\begin{enumerate}
\setcounter{enumi}{2}
\item\label{it:appl3b} If $\eta\in (\frac{d}{2mp},\frac12-\frac{d}{4m})$
is such that \eqref{eq:HBinterpolationcond} holds for the pair
$(\eta,p)$ and if $u_0\in H^{2 m \eta,p}_{\{B_j\}}(S)$ almost
surely, then for all $r>2$ and $\alpha\in
(\eta+\frac{d}{4m},\frac12-\frac1r)$ there exists a unique global
solution $u$ of \eqref{eq:SEwhitenoise2} in
$V^{0}_{\a,r}([0,T]\times\O;H^{2m\eta,p}_{\{B_j\}}(S))$.
\item\label{it:globalub} Moreover, if $\delta> \frac{d}{2mp}$ and
$\lambda\ge 0$ are such that
$\delta+\lambda<\frac{1}{2}-\frac{d}{4m}$ and
\eqref{eq:HBinterpolationcond} holds for the pair $(\delta, p)$, and
if $u_0\in H^{m-\frac{d}{2},p}_{\{B_j\}}(S)$ almost surely, then $u$
has paths in $C^{\lambda}([0,T];H^{2m \delta,p}_{\{B_j\}}(S))$
almost surely.
\end{enumerate}
\end{theorem}

\begin{remark} \
\renewcommand{\labelenumi}{(\roman{enumi})}
\renewcommand{\theenumi}{(\roman{enumi})}
\begin{enumerate}
\item For $p\in [2, \infty)$ the uniqueness result in \eqref{it:appl1b} and
\eqref{it:appl3b} can be simplified. In that
case one obtains a unique solution in
\[L^0(\O;C([0,T];H^{2m\eta,p}_{\{B_j\}}(S)))\subseteq V^{0}_{\a,r}([0,T]\times\O;H^{2m\eta,p}_{\{B_j\}}(S)).\]
For this case on could also apply martingale type $2$ integration
theory from \cite{Brz2} to obtain the result.
\item By the Sobolev embedding theorem one obtains H\"older continuous solutions in time
and space. For instance, assume in (4) that $u_0\in
C^{m-\frac{d}{2}}_{\{B_j\}}(\overline{S})$ almost surely. It follows
from
\[C^{m-\frac{d}{2}}_{\{B_j\}}(\overline{S})\hookrightarrow H^{\eta,p}_{\{B_j\}}(S) = [E,D(-A_p)]_{\frac{\eta}{2m}}
\hookrightarrow  D((-A_p)^{\frac{\eta-\e}{2m}})\] for all $p\in (1,
\infty)$ and $\eta<m-\frac{d}{2}$ and $\e>0$, that $t\mapsto S(t)
u_0$ is in $C^{\lambda}([0,T]; D((-A_p)^{\delta})$ for all
$\delta,\lambda>0$ that satisfy
$\delta+\lambda<\frac{1}{2}-\frac{d}{4m}$. Since
\[D((-A_p)^{\delta})\hookrightarrow [E,D(-A_p)]_{\delta-\e} =
H^{2m(\delta-\e),p}_{\{B_j\}}(S)\] for all $p\in (1, \infty)$ and
$\e>0$, by Sobolev embedding we obtain that the solution $u$ has
paths in $C^{\lambda}([0,T];C^{2m \delta}_B(\overline{S}))$ for all
$\delta,\lambda>0$ that satisfy
$\delta+\lambda<\frac{1}{2}-\frac{d}{4m}$.
\end{enumerate}
\end{remark}
\renewcommand{\labelenumi}{(\arabic{enumi})}
\renewcommand{\theenumi}{(\arabic{enumi})}

\begin{proof}[Proof of Theorem \ref{thm:mainexdm}]
Let $p\in (1, \infty)$ be as in the theorem and take $E:= L^p(S)$.
For $b\in (0,1)$ let $E_b$ denote the complex interpolation space
$[E,\D(A_p)]_{b}$. Note that we use the notation $E_b$ for complex
interpolation spaces instead of fractional domain spaces as we did
before. This will be more convenient, since we do not assume that
$A_p$ has bounded imaginary powers, and therefore we do not know the
fractional domain spaces explicitly. Recall (cf. \cite{Lun}) that
$E_{a}\embed\D((-A)^{b})$ and that $\D((-A)^{a}) \embed E_b$ for all
$a\in (b,1)$ for all $b\in (0,1)$.

If $b> \frac{d}{2mp}$, then by \cite[Theorem 4.6.1]{Tr} we have
\[[E,\D(A_p)]_{b} \hookrightarrow C(\overline{S}).\]
Assume now that $\eta\in (\frac{d}{2mp},\frac12-\frac{d}{4m})$. Let
$F,G:[0,T]\times\O\times E_\eta\to L^\infty(S)$ be defined as
\[(F(t,\omega,x))(s) = f(t,\omega,s,x(s)) \ \text{and} \ G(t,\omega,x))(s) = g(t,\omega,s,x(s)).\]
We show that $F$ and $G$ are
well-defined and locally Lipschitz. Fix $x,y\in
E_\eta$ and let
\[R:=\max\{{\rm ess}\sup_{s\in S}|x(s)|, {\rm ess}\sup_{s\in
S}|y(s)|\}<\infty.\] From the measurability of $x,y$ and $f$ it is
clear that $s\mapsto (F(t,\omega,x))(s)$ and $s\mapsto
(F(t,\omega,y))(s)$ are measurable. By
\ref{cond:localLipfgb} it follows that for almost all $s\in S$, for
all $t\in [0,T]$ and $\omega\in \O$ we have
\[\begin{aligned}
|(F(t,\omega,x))(s) - (F(t,\omega,y))(s)| &= |f(t,\omega,s,x(s)) -
f(t,\omega,s,y(s))|
\\ &\leq L_f^R |x(s) -
y(s)|
\\ &\leq L_f^R \|x-y\|_{L^\infty(S)}
\\ & \lesssim_{\eta} L_f^R \|x-y\|_{E_\eta}.
\end{aligned}\]
Also, by the second part of \ref{cond:localLipfgb}, for almost all
$s\in S$, for all $t\in [0,T]$ and $\omega\in \O$ we have
\[|(F(t,\omega,0))(s)| = |f(t,\omega,s,0)|<\sup_{t,s,\omega} |f(t,\omega,s,0)|<\infty.\]
Combing the above results we see that $F$ is well-defined and
locally Lipschitz. In a similar way one shows that $F$ has linear
growth (see (A2)) if \ref{cond:lingrowthfgb} holds.
The same arguments work for $G$.

Since $L^\infty(S)\hookrightarrow L^p(S)=E$ we may consider $F$ as an
$E$-valued mapping. It follows from the Pettis measurability theorem
that for all $x\in E_\eta$, $(t,\omega)\mapsto F(t,\omega,x)$
is strongly measurable in $E$ and adapted.

To model the term $g(t,x,u(t,s)) \, \frac{\partial w(t,s)}{\partial t}$, let
$H := L^2(S)$ and let $W_H$ be a cylindrical Brownian motion. Define
the multiplication operator function $\Gamma:[0,T]\times\O\times
E_\eta\to \calL(H)$ as
\[(\Gamma(t,\omega, x) h)(s) := (G(t,\omega, x))(s) h(s), \ \ s\in S.\]
Then $\Gamma$ is well-defined, because for all $t\in [0,T]$, $\omega\in \O$ we
have $G(t,\omega, x)\in L^\infty(S)$.

Now let $\theta_B>\theta_B'>\frac{d}{4m}$ be such that
$\theta_B+\eta<\frac12$ and \eqref{eq:HBinterpolationcond} holds for
the pair $(\theta_B, 2)$. Define $(-A)^{-\theta_B}B:[0,T]\times\O\times
E_\eta\to \g(H,E)$ as
\[(-A)^{-\theta_B} B(t,\omega, x) h = i (-A)^{-\theta_B} G(t,\omega, x) h,\]
where $i:H^{2 m \theta_B', 2}(S)\to L^p(S)$ is the inclusion
operator. This is well-defined, because  $(-A)^{-\theta_B}:H\to H^{2
m \theta_B', 2}(S)$ is a bounded operator and therefore
by the right-ideal
property and Corollary \ref{cor:gammaSobolevembed} it follows that
\[\|i (-A)^{-\theta_B}\|_{\g(H,E)} \leq \|(-A)^{-\theta_B}\|_{\calL(L^2(S),H^{2 m \theta_B', 2}(S))} \|i\|_{\g(H^{2 m \theta_B', 2}(S),L^p(S))} <\infty.\]
Moreover, $B$ is locally Lipschitz. Indeed, fix $x,y\in E_\eta$ and
let
\[R:=\max\{{\rm ess}\sup_{s\in S}|x(s)|, {\rm ess}\sup_{s\in
S}|y(s)|\}<\infty.\] It follows from the right-ideal property that
\[\begin{aligned}
\|i(-A)^{\theta_B}&(B(t,\omega, x) - B(t,\omega, y))\|_{\g(H,E)}
\\ &\leq \|i(-A)^{-\theta_B}\|_{\g(H,E)} \|\Gamma(t,\omega,x) -
\Gamma(t,\omega,y)\|_{\calL(H)}
\\ & \leq \|i(-A)^{\theta_B}\|_{\g(H,E)} \|G(t,\omega,x) - G(t,\omega,y)\|_{L^\infty(S)}
\\ & \leq \|i(-A)^{\theta_B}\|_{\g(H,E)} L_g^R \|x-y\|_{L^\infty(S)}
\\ & \lesssim_{a,p} \|i(-A)^{\theta_B}\|_{\g(H,E)} L_g^R \|x-y\|_{E_\eta}.
\end{aligned}\]
In a similarly way one shows that $B$ has linear growth. Notice that $B$
is $H$-strongly measurable and adapted by the Pettis measurability
theorem.

If $p\in [2, \infty)$, then $E$ has type $2$ and it follows from
Lemma \ref{lem:type2Lipschitz} that $(-A)^{-\theta_B}B$ is locally
$\gL$-Lipschitz and $B$ has linear growth in the sense of
(A3) if \ref{cond:lingrowthfgb} holds.

In case $p\in (1, 2)$ the above result holds as well. This may be
deduced from the previous case. Indeed, for each $n$ define
$(-A)^{-\theta_B} B_n:[0,T]\times\O\times E_\eta\to \g(H,E)$ as
$(-A)^{-\theta_B}B_n(x) = (-A)^{-\theta_B}B(x)$ for all
$\|x\|_{E_\eta}\leq n$ and $(-A)^{-\theta_B}B_n(x) =
(-A)^{-\theta_B}B_n\big(\frac{n x}{\|x\|}\big)$ otherwise. Define
$(-A)^{-\theta_B} B_n^\infty:[0,T]\times\O\times L^\infty(S)\to
\g(H,H)$ as $(-A)^{-\theta_B}B_n^\infty(x) =
(-A)^{-\theta_B}B_n(x)$. Replacing $E$ with $L^2(S)$ in the above
calculation it follows that $B_n^\infty$ is a Lipschitz function
uniformly on $[0,T]\times\O$. Since $H$ has type $2$,
$(-A)^{-\theta_B}B_n^\infty$ is $\gL$-Lipschitz. Fix a finite Borel
measure $\mu$ on $(0,T)$ and fix $\phi_1, \phi_2\in
L^2_{\g}((0,T),\mu;E_\eta)$. Since $H\hookrightarrow E$ continuously, it follows
that
\[\begin{aligned}
\|(-A)^{-\theta_B}&(B_n(t,\omega, \phi_1) - B_n(t,\omega,
\phi_2))\|_{\g(L^2((0,T),\mu;H),E)} \\ &\leq C
\|(-A)^{-\theta_B}(B_n^\infty(t,\omega, \phi_1) -
B_n^\infty(t,\omega, \phi_2))\|_{\g(L^2((0,T),\mu;H),H)}
\\ & \leq C \|\phi_1-\phi_2\|_{L^2((0,T),\mu; L^\infty(S))}
\\ & \leq C\|\phi_1-\phi_2\|_{L^2((0,T),\mu;E_{\eta})},
\end{aligned}
\]
where $C$ also depends on $n$. In a similarly way one shows that $B$
has linear growth in the sense of (A3)
 if $g$ has linear growth.

If $u_0\in H^{2m\beta,p}_{\{B_j\}}(S)$ almost surely, where
$\beta\in (\frac{d}{2mp},\frac12 - \frac{d}{4m}]$ is such that
\eqref{eq:HBinterpolationcond} holds for the pair $(\beta, p)$, then
$\omega\mapsto u_0(\cdot,\omega)\in
H^{2m\beta,p}_{\{B_j\}}(S)=E_{\beta}$ is strongly $\F_0$-measurable.
This follows from the Pettis measurability theorem.

(\ref{it:appl1b}): It follows from Theorem
\ref{thm:mainexistenceLocallyLipschitzCase} with $\eta$, $\theta_B$
as above and with $\eta+\theta_B<\frac12$ and $\theta_F=0$,
$\tau=p\wedge 2$ that there is a unique maximal local mild solution
$(U(t))_{t\in [0,\varrho)}$ in $V^{0,{\rm
loc}}_{\alpha,r}([0,\varrho)\times \O;E_\eta)$ for all $\a>0$ and
$r>2$ satisfying $\eta+\theta_B<\alpha<\frac12 - \frac1r$. In
particular $U$ has almost all paths in $C([0,\varrho),E_\eta)$. Now
take $u(t,\omega, s) := U(t,\omega)(s)$ to finish the proof of
(\ref{it:appl1b}).

(\ref{it:appl2b}): Let $\delta=\eta>\frac{d}{2mp}$ and $\lambda\ge
0$ be such that $\lambda+\delta<\frac12-\frac{d}{4m}$. Choose
$\theta_B>\frac{d}{4m}$ such that $\lambda+\delta<\frac12-\theta_B$.
It follows from Theorem \ref{thm:mainexistenceLocallyLipschitzCase}
that almost surely, $U-S u_0\in C^\lambda_{\rm
loc}([0,\varrho(\omega));H^{2m\delta,p}_{\{B_j\}}(S))$. First
consider the case that $(\frac12-\frac{d}{4m},p)$ satisfies
\eqref{eq:HBinterpolationcond}. Since $u_0\in
H^{m-\frac{d}{2},p}_{\{B_j\}}(S)=E_{\frac12-\frac{d}{4m}}\subseteq
E_\delta$ almost surely and $\lambda+\delta<\frac12-\frac{d}{4m}$ we
have $S u_0\in C^{\lambda}([0,T];H^{2m\delta,p}_{\{B_j\}}(S))$
almost surely. Therefore, almost all paths of $U$ are in
$C^\lambda_{\rm
loc}([0,\varrho(\omega));H^{2m\delta,p}_{\{B_j\}}(S))$. In the case
$(\frac12-\frac{d}{4m},p)$ does not satisfy
\eqref{eq:HBinterpolationcond} one can redo above argument with
$\frac12-\frac{d}{4m}-\epsilon$ for $\epsilon>0$ small. This proves
(\ref{it:appl2b}).

(\ref{it:appl3b}), (\ref{it:globalub}): This follows from Theorems
\ref{thm:mainexistenceLocallyLipschitzCase} and parts (3), (4) of
\ref{thm:mainexistenceLocallyLipschitzCase}.
\end{proof}

\begin{remark}
The above approach also works for systems of equations.
\end{remark}

\paragraph{\bf Laplacian in $L^p$}

Let $S$ be a open subset (not necessarily bounded) of $\R^d$.
Consider the following perturbed heat equation with Dirichlet
boundary values:

\begin{equation} \label{eq:SELpS}
\begin{aligned}
\frac{\partial u}{\partial t}(t,s) & = \Delta u(t,s) +
f(t,s,u(t,s))
\\ \nonumber  & \ \ \ + \sum_{n\geq 1} b_n(t,s,u(t,s)) \, \frac{\partial W_n(t)}{\partial t},
&& s\in S, \ t\in
(0,T],
\\  u(t,s) &= 0, && s\in \partial S, \ t\in (0,T],
\\ u(0,s) &= u_0(s), && s\in S.
\end{aligned}
\end{equation}
The functions $f,b_n:[0,T]\times\O\times S\times \R\to \R$ are
jointly measurable, and adapted in the sense that for each $t\in
[0,T]$, $f(t,\cdot)$ and $b_n(t,\cdot)$ are $\F_t\otimes
\mathcal{B}_{S}\otimes \mathcal{B}_{\R}$-measurable. We assume that
$(W_n)_{n\geq 1}$ is a sequence of independent standard Brownian
motions on $\O$ and $u_0:S\times \O \to \R$ is an
$\mathcal{B}_{S}\otimes\F_0$-measurable initial value condition. We
say that $u:[0,T]\times\O\times S\to \R$ is a {\em solution} of
\eqref{eq:SELpS} if the corresponding functional analytic model
\eqref{SE} has a mild solution $U$ and $u(t,s,\omega) =
U(t,\omega)(s)$.

Let $p\in [1, \infty)$ be fixed and let $E:=L^p(S)$. It is well-known
that the Dirichlet Laplacian $\Delta_p$ generates a uniformly exponentially
stable and analytic $C_0$-semigroup $(S_p(t))_{t\geq 0}$ on $L^p(S)$, and
under a regularity assumption on $\partial S$ one can identify
$\D(\Delta_p)$ as $W^{2,p}(S)\cap W^{1,p}_0(S)$. Consider the
following $p$-dependent condition:
\let\ALTERWERTA\theenumi
\let\ALTERWERTB\labelenumi
\def\theenumi{(C)}
\def\labelenumi{(C)}
\begin{enumerate}
\item \label{cond:Lipandbounded}
There exist constants $L_f$ and $L_{b_n}$ such that
\[
\begin{aligned}
|f(t,\omega,s,x) - f(t,\omega,s,y)| & \leq L_f |x-y|, \\
|b_n(t,\omega,s,x) - b_n (t,\omega,s,y)| & \leq L_{b_n} |x-y|,
\end{aligned}
\]
for all $t\in [0,T]$, $\omega\in \O$, $s\in S$, and $x,y\in \R$.
Furthermore, $f$ satisfies the boundedness condition
\[
\sup \|f(t,\omega,\cdot,0)\|_{L^p(S)}<\infty,
\]
where the supremum is taken over all $t\in [0,T]$ and $\omega\in\O$,
and the $b_n$ satisfy the following boundedness
condition: for all finite measures $\mu$ on $(0,T)$,
\[\sup \Big\|\Big(\int_0^T
\sum_{n\geq 1}|b_n(t,\omega,\cdot,0)|^2\,
d\mu(t)\Big)^{\frac12}\Big\|_{L^p(S)}<\infty,
\]
where the supremum is taken over all $\omega\in \O$.

\end{enumerate}
\let\theenumi\ALTERWERTA
\let\labelenumi\ALTERWERTB

\begin{theorem}
Let $S$ be an open subset of $\R^d$ and let $p\in [1, \infty)$.
Assume that condition {\rm\ref{cond:Lipandbounded}} holds with
$\sum_{n\geq 1} L_{b_n}^2<\infty$. If $u_0\in L^p(S)$ almost surely,
then for all $\alpha>0$ and $r>2$ such that
$\alpha<\frac12-\frac1r$, the problem \eqref{eq:SELpS} has a unique
solution $U\in V^{0}_{\a,r}([0,T]\times\O;L^p(S))$. Moreover, for
all $\l\ge 0$ and $\d\ge 0$ such that $\lambda+\delta<\frac12$ there
is a version of $U$ such that almost surely, $t\mapsto U(t) - S_p(t)
u_0$ belongs to $C^\lambda([0,T];[L^p(S),\D(\Delta_p)]_{\delta}).$
\end{theorem}

\begin{remark}
Under regularity conditions on $\partial S$ and for $p\in (1, \infty)$ one
has
\[[L^p(S),D(\Delta_p)]_{\delta}
=\big\{x\in H^{2\delta,p}(S): x=0 \ \ \text{on $\partial S$ if
$2\delta-\frac1p>0$}\big\}\] provided $\delta\in (0,1)$ is such that
$2\delta-\frac1p\neq 0$.
\end{remark}

\begin{proof}
We check the conditions of Theorem \ref{thm:mainexistenceLloc} (for
$p=1$ we use the extensions of our results to UMD$^-$ spaces described in Section
\ref{sec:gen}, keeping in mind that $L^1$-spaces have this property).
It was already noted that \ref{as:semigroup} is
fulfilled. Let $E:= L^p(S)$ and define $F:E\to E$ as $F(t,x)(s) := f(t,s,x(s))$.
One easily checks that $F$ satisfies \ref{as:LipschitzF} with
$\theta_F=\eta=0$. Let $H:=l^2$ with standard unit basis $(e_n)_{n\ge 1}$,
and let $B:[0,T]\times\O\times E\to
\calL(H,E)$ be defined as $(B(t,\omega,x) e_n)(s) :=
b_n(t,\omega,s,x(s))$. Then for all finite measures $\mu$ on $(0,T)$ and all
$\phi_1,\phi_2\in \g(L^2((0,T),\mu;H),E)$,
\[\begin{aligned}
\|B&(\cdot, \phi_1) - B(\cdot, \phi_2)\|_{\g(L^2((0,T),\mu;H),E)} \\ &
\eqsim_p \Big\| \Big(\int_0^T \sum_{n\geq 1}|b_n(t,
\cdot,\phi_1(t)(\cdot)) - b_n(t,\cdot, \phi_2(t)(\cdot))|^2\,
d\mu(t)\Big)^{\frac12}\Big\|_{E}
\\ & \leq \Big\|\Big(\int_0^T \sum_{n\geq 1} L_{b_n}^2 |\phi_1(t)(\cdot) -
\phi_2(t)(\cdot)|^2\, d \mu(t)\Big)^{\frac12}\Big\|_{E}
\\ & \eqsim_p L \|\phi_1 - \phi_2\|_{\g(L^2(0,T),\mu),E)},
\end{aligned}\]
where $L=(\sum_{n\geq 1} L_{b_n}^2)^{\frac12}$. Moreover,
\[\|B(\cdot, 0)\|_{\g(L^2(0,T),\mu;H),E))} \eqsim_p \Big\|
\Big(\int_0^T \sum_{n\geq 1}|b_n(t,\cdot,0)|^2\, d \mu(t)\Big)^{\frac12}\Big\|_{E} <\infty. \]
From these two estimates one can obtain \ref{as:LipschitzB}.
\end{proof}

\appendix
\section{Measurability of stochastic convolutions}\label{sec:measconv}

In this appendix we study progressive
measurability properties of processes of
the form \[t\mapsto \int_0^t \Phi(t,s) \, d W_H(s)\] where $\Phi$ is
a two-parameter process with values in $\calL(H,E)$.

\begin{proposition}\label{prop:convprogr}
Let $E$ be a UMD space. Assume that $\Phi: \R_+\times\R_+\times \O
\to \calL(H,E)$ is $H$-strongly measurable and for each $t\in \R_+$,
$\Phi(t, \cdot)$ is adapted and has paths in $\g(L^2(\R_+;H),E)$
almost surely. Then the process $\zeta:\R_+\times\O\to E$,
\[\zeta(t) = \int_0^t \Phi(t,s) \, d W_H(s),\]
has a version which is adapted and strongly measurable.
\end{proposition}

\begin{proof}
It suffices to show that $\zeta$ has a strongly measurable version
$\tilde{\zeta}$, the adaptedness of $\tilde{\zeta}$ being clear.
Below we use strong measurability for metric spaces as in
\cite{VTC}.

Let $L^0_{\bbF}(\O;\g(L^2(\R_+;H),E))$ denote the closure of all
adapted strongly measurable processes which are almost surely in
$\g(L^2(\R_+;H),E)$. Note that by \cite{NVW} the stochastic integral
mapping extends to $L^0_{\bbF}(\O;\g(L^2(\R_+;H),E))$.

Let $G\subseteq \R_+\times\O$ be the set of all $(t,\omega)$ such
that $\Phi(t,\cdot,\omega)\in \g(L^2(\R_+;H),E)$. Since $\Phi$ is
$H$-strongly measurable, we have $G\in {\mathcal
B}_{\R_+}\otimes\mathcal{A}$. Moreover, letting $G_t=\{\omega\in \O:
(t,\omega)\in G\}$ for $t\in \R_+$, we have $\P(G_t) =1$ and
therefore $G_t\in \F_0$. Define the $H$-strongly measurable function
$\Psi:\R_+\times\R_+\times\O\to\mathcal{B}(H,E)$ as
$\Psi(t,s,\omega) := \Phi(t,s,\omega) \one_{[0,t]}(s)
\one_{G}(t,\omega)$. It follows from \cite[Remark 2.8]{NVW} that the
map $\R_+\times\O\ni (t,\omega)\to \Psi(t,\cdot,\omega) \in
\g(L^2(\R_+;H),E)$ is strongly measurable. Hence, the map $\R_+\ni
t\to \Psi(t,\cdot) \in L^0(\O;\g(L^2(\R_+;H),E))$ is strongly
measurable. Since it takes values in
$L^0_{\bbF}(\O;\g(L^2(\R_+;H),E))$ it follows from an approximation
argument that it is strongly measurable as an
$L^0_{\bbF}(\O;\g(L^2(\R_+;H),E))$-valued map. Since the elements
which are represented by an adapted step process are dense in
$L^0_{\bbF}(\O;\g(L^2(\R_+;H), E))$, it follows from
\cite[Proposition 1.9]{VTC} that we can find a sequence of processes
$(\Psi_n)_{n\geq 1}$, where each $\Psi_n:\R_+\to
L^0_{\bbF}(\O;\g(L^2(\R_+;H),E))$ is a countably valued simple
function of the form
\[\Psi_n = \sum_{k\geq 1} \one_{B_k^n} \Phi^n_{k}, \ \ \text{with} \ B^n_k\in \mathcal{B}_{\R_+} \ \text{and} \ \Phi^n_{k}\in L^0_{\bbF}(\O;\g(L^2(\R_+;H),E)),\]
such that for all $t\in \R_+$ we have
$\|\Psi(t)-\Psi_n(t)\|_{L^0(\O;\g(L^2(\R_+;H),E))}\le 2^{-n}$, where
with a slight abuse of notation we write
\[ \|\xi\|_{L^0(\O;F)} :=
\E (\n \xi\n_F\wedge 1)\] keeping in mind that this is not a norm.
Notice that by the Chebyshev inequality, for a random variable
$\xi:\O\to F$, where $F$ is a normed space, and $\e\in (0,1]$, we
have
\[\P(\|\xi\|_F>\e)= \P((\|\xi\|_F\wedge 1)>\e) \le \e^{-1}\|\xi\|_{L^0(\O;F)}.\]
It follows from \cite[Theorems 5.5 and 5.9]{NVW} that for all $t\in
\R_+$, for all $n\geq 1$ and for all $\e,\d\in (0,1]$,
\[
\P\Big(\Big\|\int_{\R_+}\Psi(t,s)-\Psi_n(t,s)\, d
W_H(s)\Big\|>\e\Big)\le \frac{C \d^2}{\e^2} + \frac{1}{\d 2^n}.
\]
Taking $\e\in (0,1]$ arbitrary and $\d=\frac1n$, it follows from the
Borel-Cantelli lemma that for all $t\in \R_+$,
\[\P\Big(\bigcap_{N\geq 1} \bigcup_{n\geq N} \Big\{\Big\|\int_{\R_+}\Psi(t,s)-\Psi_n(t,s)\, d W_H(s)\Big\|>\e\Big\}\Big) =0.\]
Since $\e\in (0,1]$, was arbitrary, we may conclude that for all $t\in \R_+$,
almost surely,
\[\zeta(t, \cdot) = \int_{\R_+}\Psi(t,s)\, d W_H(s) = \limn \int_{\R_+} \Psi_n(t,s)\, d W_H(s).\]
Clearly,
\[\int_{\R_+}\Psi_n(\cdot,s)\, d W_H(s) = \sum_{k\geq 1} \one_{B_k^n}(\cdot) \int_{\R_+} \Phi^n_{k}(s) \, d W_H(s)\]
has a strongly $\mathcal{B}_{\R_+} \otimes\mathcal{F}_{\infty}$-measurable
version, say $\zeta_n:\R_+\times\O\to E$. Let $C\subseteq \R_+\times\O$ be the
set of all points $(t,\omega)$ such that $(\zeta_n(t, \omega))_{n\geq 1}$
converges in $E$. Then $C\in \mathcal{B}_{\R_+} \otimes\mathcal{F}_{\infty}$
and we may define the process $\tilde{\zeta}$ as $\tilde{\zeta} = \limn
\zeta_n \one_{C}$. It follows that $\tilde{\zeta}$ is strongly
$\mathcal{B}_{\R_+}\otimes\mathcal{F}_{\infty}$-measurable and for all $t\in
\R_+$, almost surely, $\tilde{\zeta}(t, \cdot) = \zeta(t, \cdot)$.
\end{proof}

{\em Acknowledgment} -- We thank Tuomas Hyt\"onen for suggesting an improvement
in Proposition \ref{prop:gammanormconv}.

\providecommand{\bysame}{\leavevmode\hbox
to3em{\hrulefill}\thinspace}

\end{document}